\documentclass[preprint,12pt,letterpaper]{elsarticle}

\usepackage{amsmath,amssymb,amsfonts,amsthm}
\usepackage{enumitem}
\usepackage{array}
\usepackage[expansion=false]{microtype}
\usepackage[hidelinks]{hyperref}
\usepackage{mathrsfs}
\usepackage{stmaryrd}
\input{Ustmry.fd}
\DeclareFontShape{U}{stmry}{b}{n}{<->ssub*stmry/m/n}{}
\usepackage[all]{xy}
\usepackage{tikz}
\usetikzlibrary{backgrounds}
\usepackage{etoolbox}

\makeatletter
\let\ps@pprintTitle\ps@empty
\patchcmd{\pprintMaketitle}
  {\elsaddress\par\vskip36pt}
  {\elsaddress\par\vskip24pt}
  {}{\PackageError{whittaker}{Could not adjust title spacing}{Check the elsarticle title definition.}}
\patchcmd{\pprintMaketitle}
  {\ifvoid\absbox\else\unvbox\absbox\par\vskip10pt\fi}
  {\ifvoid\absbox\else\unvbox\absbox\par\vskip12pt\hrule height0.4pt\vskip12pt\fi}
  {}{\PackageError{whittaker}{Could not separate the abstract and keywords}{Check the elsarticle title definition.}}
\patchcmd{\pprintMaketitle}
  {\ifvoid\keybox\else\unvbox\keybox\par\vskip10pt\fi\hrule\vskip12pt}
  {\ifvoid\keybox\else\unvbox\keybox\par\fi\vskip18pt}
  {}{\PackageError{whittaker}{Could not adjust the keyword separator}{Check the elsarticle title definition.}}
\makeatother

\newcommand{\ydraw}[3]{%
 \tikz[x=#1cm,y=#1cm,baseline=(current bounding box.center)]{%
  \foreach \len [count=\r] in {#2}{\foreach \c in {1,...,\len}{%
    \ifnum#3=1 \ifnum\r=1 \fill[.!18] (\c-1,-\r) rectangle (\c,1-\r);\fi\fi
    \ifnum#3=2 \ifnum\c=1 \fill[.!18] (\c-1,-\r) rectangle (\c,1-\r);\fi\fi
    \draw[line width=0.4pt] (\c-1,-\r) rectangle (\c,1-\r);}}}}

\allowdisplaybreaks
\numberwithin{equation}{section}
\theoremstyle{plain}
\newtheorem{theorem}{Theorem}[section]
\newtheorem{proposition}[theorem]{Proposition}
\newtheorem{lemma}[theorem]{Lemma}
\newtheorem{corollary}[theorem]{Corollary}
\theoremstyle{remark}
\newtheorem{remark}[theorem]{Remark}
\newtheorem{example}[theorem]{Example}

\newcommand{\Fq}{\mathbb F_q}
\newcommand{\Sp}{\operatorname{Sp}}
\newcommand{\SO}{\operatorname{SO}}
\newcommand{\OO}{\operatorname{O}}
\newcommand{\GL}{\operatorname{GL}}
\newcommand{\U}{\operatorname{U}}
\newcommand{\bSp}{\mathbf{Sp}}
\newcommand{\bO}{\mathbf{O}}
\newcommand{\bSO}{\mathbf{SO}}
\newcommand{\Irr}{\operatorname{Irr}}
\newcommand{\Hom}{\operatorname{Hom}}
\newcommand{\Wh}{\operatorname{Wh}}
\newcommand{\rank}{\operatorname{rank}}
\newcommand{\defect}{\operatorname{def}}
\newcommand{\LC}{\operatorname{LC}}
\newcommand{\CO}{\mathcal O}
\newcommand{\Ost}[1]{\mathcal O^{\rm st}\!\left(#1\right)}
\newcommand{\ScO}{\mathscr O}
\newcommand{\bG}{\mathbf{G}}
\newcommand{\bH}{\mathbf{H}}
\newcommand{\Fg}{\mathfrak g}
\newcommand{\CN}{\mathcal N}
\newcommand{\CL}{\mathcal L}
\newcommand{\cD}{\mathcal D}
\newcommand{\cE}{\mathcal E}
\newcommand{\bA}{\overline A}
\newcommand{\fd}{\mathrm{fd}}
\newcommand{\fdsearrow}{\overset{\fd}{\searrow}}
\newcommand{\FJ}{\mathrm{FJ}}
\newcommand{\disc}{\operatorname{disc}}

\begin{document}
\begin{frontmatter}

\title{Polynomial Degrees and Leading Terms of Generalized Whittaker Multiplicities for Finite Symplectic and Orthogonal Groups}

\author[hku]{Pengcheng Li}
\ead{pcli25@hku.hk}
\author[jlu]{Zhicheng Wang}
\ead{wangzhicheng@jlu.edu.cn}

\address[hku]{Department of Mathematics, The University of Hong Kong,
Hong Kong, China}
\address[jlu]{Department of Mathematics, Jilin University,
Changchun 130012, Jilin, P. R. China}
\hypersetup{pdfauthor={Pengcheng Li and Zhicheng Wang}}

\begin{abstract}
We study generalized Whittaker multiplicities for irreducible representations
of finite symplectic and full orthogonal groups over $\mathbb F_q$.
For fixed combinatorial data in Lusztig's parametrization and a fixed
rational nilpotent orbit, we prove that generalized Whittaker multiplicities are given by
polynomials in $q$ for all sufficiently large $q$.
Moreover, we determine the degree of this polynomial in terms of stable
nilpotent orbit dimensions and compute its leading coefficient.
\end{abstract}


\end{frontmatter}
\enlargethispage{1pt}

\section{Introduction}\label{sec:intro}

Let $\Fq$ be a finite field of odd characteristic $p$, with algebraic
closure $\overline{\Fq}$. Let $V$ be a nondegenerate symplectic or quadratic
space over $\Fq$, and let $\bG$ be its algebraic isometry group.
For the corresponding Frobenius endomorphism $F$, let
$G=\bG(\Fq)=\bG^F$;
thus $G=\Sp(V)$ or $\OO(V)$, respectively.
\par
If the root system of $\bG^\circ$ over $\overline{\Fq}$ is nonempty,
we denote its Coxeter number by $h_{\bG}$ and assume, as in \cite{wang},
that
\begin{equation}\label{eq:char-hyp}
 p>3(h_{\bG}-1).
\end{equation}
\par
Whenever we use finite Gan--Gross--Prasad (GGP)
branching or the chosen Jordan decomposition, we assume that $q\ge q_0$,
where $q_0$ is large enough for the results of
\cite{Wang21,WangGGP,WangJordanCorrected} to apply.
We fix a nontrivial additive character $\psi:\Fq\to\mathbb C^\times$.

Write $\Fg=\operatorname{Lie}(\bG)(\overline{\Fq})$ and let $\CN\subset\Fg$ be its nilpotent cone.
For $x\in\CN^F$, write $\CO=G\cdot x$ for its \emph{rational nilpotent orbit}
and $\CO^{\rm st}=\bG(\overline{\Fq})\cdot x$ for the corresponding
\emph{stable nilpotent orbit}. The orbit $\CO^{\rm st}$ is stable under
$\operatorname{Gal}(\overline{\Fq}/\Fq)$. Dimensions of
stable orbits and algebraic groups are taken over $\overline{\Fq}$.
Let $\Gamma_\CO$ be Kawanaka's generalized Gelfand--Graev representation of $G$
attached to $\CO$ \cite{Kawanaka85,Kawanaka87}.
Let $\pi$ be an irreducible complex representation of $G$.
We denote by
$\Wh_\CO(\pi)$ the generalized Whittaker model of $\pi$ attached to $\CO$
and write
\begin{equation*}
 m_\CO(\pi):=\dim\Wh_\CO(\pi)=\dim\Hom_{G}(\pi,\Gamma_\CO).
\end{equation*}

For a fixed $\pi$, we consider two basic questions  about its
generalized Whittaker models:
\begin{itemize}
 \item[(1)] For which rational nilpotent orbits $\CO$ is $\Wh_\CO(\pi)$ nonzero?
 \item[(2)] If $\Wh_\CO(\pi)\ne0$, what is its dimension?
\end{itemize}
The first question amounts to determining the rational wavefront set
\begin{equation*}
 \ScO(\pi):=\{\CO:m_\CO(\pi)\ne0\}.
\end{equation*}
Let $\le_{\rm rat}$ denote the rational orbit order defined by the second
author in \cite{wang} and recalled in Section~\ref{subsec:rational-orbit-order}.
We write
$ \ScO(\pi)^{\max}$ for the set of maximal elements of $\ScO(\pi)$ with respect to this order.
This set may contain more than one rational orbit.

On stable nilpotent orbits, we use the Zariski closure order:
$\CO_1^{\rm st}\le\CO^{\rm st}_2$ means
$\CO^{\rm st}_1\subseteq\overline{\CO^{\rm st}_2}$, with closure taken over $\overline{\Fq}$.
Under our standing assumptions, the set of underlying stable orbits
$\{\CO^{\rm st}:\CO\in\ScO(\pi)\}$ has a unique maximal element
\cite{Lusztig92,GM00,AA07,Taylor16}. We call it the \emph{stable wavefront orbit}
of $\pi$ and denote it by
$
 \CO_\pi^{\rm st}.$
Thus every $\CO\in\ScO(\pi)$ satisfies
$\CO^{\rm st}\subseteq\overline{\CO_\pi^{\rm st}}$.

For the general linear group $\GL_n(\Fq)$ and the finite unitary group $\U_n(\Fq)$,
the computation of generalized Whittaker
multiplicities goes back to Kawanaka and Lusztig. By Kawanaka's character
formula \cite[(3.2.14)]{Kawanaka85}, the generalized Gelfand--Graev
characters can be expressed in terms of Green polynomials, and their irreducible
multiplicities can be computed.
Lusztig subsequently obtained a general geometric formula in terms of
intersection cohomology complexes \cite[Theorem~7.3]{Lusztig92}.
Further developments include a combinatorial treatment for $\GL_n(\Fq)$
\cite{AndrewsThiem17} and a descent approach for finite unitary groups
\cite{PengWangU23}. For certain degenerate Whittaker models of cuspidal
representations, see also \cite{Prasad00}.

For $\GL_n(\Fq)$ and $\U_n(\Fq)$, the $\Fq$-points of each stable
nilpotent orbit form a single rational orbit. For symplectic and
orthogonal groups, however, the $F$-rational points of a single stable
nilpotent orbit may form several rational orbits; see
Section~\ref{subsec:wavefront-def}. The stable wavefront orbit alone therefore does
not suffice to determine the nonvanishing or the dimension of a Whittaker
model. 

For these groups, the second author \cite{wang}
answered the first question for sufficiently large $p$ and $q$ by determining
the maximal rational wavefront orbits and proving the following support formula
\cite[Theorem~1.1]{wang}:
\[
 \ScO(\pi)=\{\CO:\CO\le_{\rm rat}\CO_0
                   \text{ for some }\CO_0\in\ScO(\pi)^{\max}\}.
\]
The multiplicities at the maximal rational wavefront orbits are given in
\cite[Theorem~1.9]{wang}. The proofs there use the descent
method, which goes back to the automorphic descent of Ginzburg, Rallis and
Soudry \cite{GRS11}. For related work on unipotently supported class
functions, unipotent supports and generalized Gelfand--Graev
representations, see \cite{DLM03,DLM14,GH08,BDT20}; for Fourier--Jacobi
models of Deligne--Lusztig characters of unitary groups, see
\cite{LiuMaShi26}.

In this paper we address the second question for every orbit in the rational
wavefront set, in the symplectic and orthogonal cases. Assume that $p$ and $q$ satisfy the standing assumptions above, and fix the combinatorial data in Lusztig's parametrization and the rational nilpotent orbit. We show that the multiplicities $m_\CO(\pi)$, $\CO\in\ScO(\pi)$, are then polynomials in $q$, and we compute their degrees and leading coefficients.

Let $\varepsilon_q$ be the quadratic character of $\Fq^\times$, and set
$\delta=\varepsilon_q(-1)=(-1)^{(q-1)/2}$. Let $q$ run through the odd prime powers
satisfying our standing assumptions. Under Lusztig's Jordan decomposition
(Section~\ref{sec:classification}), each $\pi\in\Irr(G)$ lies in a Lusztig series
$\cE(G,s)$ for a semisimple element $s$ of the dual group and corresponds to a unipotent
representation of its centralizer $C_{G^*}(s)$. The eigenvalues of $s$ in
$\overline{\Fq}^\times\setminus\{\pm1\}$ partition into 
$[a_i]:=\{a_i^{q^j},a_i^{-q^j}:j\in\mathbb Z\}$ ($i\in I$). Each class $[a_i]$ has
even cardinality, with \emph{degree} $d_i:=\#[a_i]/2$ and \emph{kind}
$\kappa_i\in\{\GL,\U\}$ given by
\[
 \kappa_i=\begin{cases}
  \U,&\text{if }a_i^{-1}=a_i^{q^j}\text{ for some }j\ge0,\\
  \GL,&\text{otherwise};
 \end{cases}
\]
the corresponding factor $G^*_{[a_i]}(s)$ of $C_{G^*}(s)$ is isomorphic to
$\GL_{m_i}(\mathbb F_{q^{d_i}})$ (if $\kappa_i=\GL$) or
$\U_{m_i}(\mathbb F_{q^{d_i}})$ (if $\kappa_i=\U$), and the unipotent factor of $\pi$ at
$[a_i]$ is parametrized by a nonempty partition $\lambda_i\vdash m_i$. The unipotent factors at the eigenvalues $1$ and $-1$ are
parametrized by Lusztig symbols $\Lambda_1$ and $\Lambda_{-1}$
(Section~\ref{subsec:symbols}), together with a sign $\epsilon_0\in\{\pm\}$ occurring
only in the odd orthogonal case. Recording this combinatorial data (without the specific
eigenvalues $a_i$) defines the \emph{type} of $\pi$ (Section~\ref{sec:qsystems}):
\[
 \mathfrak t=\bigl(\{(d_i,\kappa_i,\lambda_i)\}_{i\in I};
                         \Lambda_1,\Lambda_{-1};\epsilon_0\bigr),
\]
where the triples $(d_i,\kappa_i,\lambda_i)$ form a multiset. 

As usual, $\Irr(K)$ denotes the set of isomorphism classes of irreducible complex representations of a finite group $K$, and we write
$\Irr_{\mathfrak t}(\bG(\Fq))$ for the set of irreducible representations of $\bG(\Fq)$ of type $\mathfrak t$.
A $q$-system of type $\mathfrak t$ is a choice of representations
$\boldsymbol\pi=\{\pi_q\}_q$ with $\pi_q\in\Irr_{\mathfrak t}(\bG(\Fq))$
for sufficiently large $q$; see Section~\ref{sec:qsystems} for details.
By the degree formula in Lusztig's Jordan decomposition and the
unipotent character-degree formulas (see Section~\ref{sec:qsystems}), we have
\begin{equation*}
 \pi_q(1)=D_{\mathfrak t}(q),\qquad D_{\mathfrak t}(X)\in\mathbb Q[X].
\end{equation*}

\subsection{Main results}
We retain our standing assumptions on $p$ and $q$.
Fix a type $\mathfrak t$, and let $\CO$ be a rational nilpotent orbit (see Section~\ref{subsec:rational-orbit-order}).

\begin{theorem}\label{thm:mainA}
There is a  polynomial
$m_{\CO,\mathfrak t}(X)\in\mathbb Q[X]$ such that, for every sufficiently large $q$,
\[
 m_\CO(\pi)=m_{\CO,\mathfrak t}(q)
\]
for all $\pi\in\Irr_{\mathfrak t}(\bG(\Fq))$.
In particular, the polynomial is independent of the choice of $\pi$
within the type $\mathfrak t$.
Moreover,
\begin{equation}\label{eq:main-degree}
 \deg m_{\CO,\mathfrak t}(X)
 =\tfrac12\bigl(\dim\CO_{\pi_q}^{\rm st}-\dim\CO^{\rm st}\bigr).
\end{equation}
\end{theorem}
By Theorem~\ref{thm:mainA}, for $q$ sufficiently large, the rational wavefront set $\ScO(\pi)$ of $\pi\in\Irr_{\mathfrak t}(\bG(\Fq))$ depends only on $\mathfrak t$. Hence so do $\ScO(\pi)^{\max}$ and $\CO_\pi^{\rm st}$.

For $\CO\in\ScO(\pi)$, define
\begin{equation*}
 Q_\pi(\CO):=\#\{\CO_0\in\ScO(\pi)^{\max}:\CO\le_{\rm rat}\CO_0\}.
\end{equation*}
For a stable nilpotent orbit $\CO^{\rm st}$ of an algebraic group $\bH$,
choose $x\in\CO^{\rm st}$ and write $A(\CO^{\rm st})=C_\bH(x)/C_\bH(x)^\circ$
for its component group.
For $\pi_q$ in the Lusztig series $\cE(G,s_q)$, define
$\pi_q^{\rm J}:=\CL_{s_q}(\pi_q)$, where $\CL_{s_q}$ is Lusztig's
Jordan decomposition; see Section~\ref{sec:classification}. Its stable
wavefront orbit lies in the Lie algebra of $C_{\bG^*}(s_q)$, where $\bG^*$ is the dual group of $\bG$. Write
$\bA(\CO_{\pi_q^{\rm J}}^{\rm st})$ for Lusztig's canonical quotient of
this orbit. We write
$[X^d]P$ for the coefficient of $X^d$ in $P\in\mathbb Q[X]$, and $\LC(P)$ for the leading coefficient of a nonzero polynomial $P$. Set
\[
 c_G=\begin{cases}2,&G\text{ odd orthogonal},\\
                  1,&G\text{ symplectic or even orthogonal},
       \end{cases}
\]
where the factor $c_G=2$ in odd orthogonal dimension reflects the direct product
decomposition $\OO_{2n+1}(\Fq)\simeq\SO_{2n+1}(\Fq)\times\{\pm1\}$.

The leading coefficient involves an integer
$\Delta_*(\CO,\CO_\pi^{\rm st})$ determined by the partitions of
$\CO^{\rm st}$ and $\CO_\pi^{\rm st}$, given by $\Delta(\CO,\CO_\pi^{\rm st})$ in
\eqref{eq:Delta-stable-partitions} for symplectic groups and by
$\Delta_{\OO}(\CO,\CO_\pi^{\rm st})$ in \eqref{eq:Delta-orth-combinatorial} for
orthogonal groups.
\begin{theorem}\label{thm:mainB}
If $m_{\CO,\mathfrak t}(X)\ne0$, then, for every sufficiently large $q$ and every
$\pi_q\in\Irr_{\mathfrak t}(\bG(\Fq))$,
\begin{align}
 \LC(m_{\CO,\mathfrak t}(X))
 &=2^{-\Delta_*(\CO,\CO_{\pi_q}^{\rm st})}
 \sum_{\substack{\CO_0\in\ScO(\pi_q)^{\max}\\\CO\le_{\rm rat}\CO_0}}
 \dim\Wh_{\CO_0}(\pi_q)
 \label{eq:intro-leading-direct}\\
 &=2^{-\Delta_*(\CO,\CO_{\pi_q}^{\rm st})}
 \frac{Q_{\pi_q}(\CO)|A(\CO_{\pi_q}^{\rm st})|}
 {c_G|\ScO(\pi_q)^{\max}|\,|\bA(\CO_{\pi_q^{\rm J}}^{\rm st})|}.
 \notag
\end{align}
\end{theorem}
By Theorem~\ref{thm:mainA} and the remark following it,
both right-hand sides are independent of $q$ and of $\pi_q$.
For the rational orbit $\CO=\{0\}$ in type $C$, this recovers the classical degree and leading
coefficient of $D_{\mathfrak t}(X)$; see Remark~\ref{rem:geck-taylor-endpoints}.

\subsection{Strategy of the proof}\label{subsec:proof-strategy}
We first treat $G=\Sp_{2N}(\Fq)$. Our strategy is analogous to the classical branching rule for symmetric groups. Let $V_\lambda$
be the irreducible complex representation of $S_n$ indexed by a
partition $\lambda\vdash n$. By Young's branching rule,
\[
 \operatorname{Res}^{S_n}_{S_{n-1}}V_\lambda
 \simeq\bigoplus_{\lambda\searrow\nu}V_\nu,
\]
where $\lambda\searrow\nu$ means that the Young diagram of $\nu$ is
obtained from that of $\lambda$ by removing a box. Thus each constituent
occurs with multiplicity one. The Young graph has partitions as vertices,
with an edge between two diagrams that differ by one box
\cite[Section~5]{VershikOkounkov05}. Iterating along $S_n\supset S_{n-1}\supset\cdots\supset S_m$, we see that for
$\mu\vdash m$ the multiplicity of $V_\mu$ in
$\operatorname{Res}^{S_n}_{S_m}V_\lambda$ equals the number of downward paths in the
Young graph from $\lambda$ to $\mu$. In particular, restricting all the way to the
trivial group $S_0$ (whose trivial representation $1_{S_0}=V_\varnothing$ is indexed by
the empty partition) yields
\[
\begin{aligned}
 \dim V_\lambda
 &=\dim\Hom_{S_0}\!\left(
    1_{S_0},\operatorname{Res}^{S_n}_{S_0}V_\lambda\right)\\
 &=\#\left\{\lambda=\lambda^{(0)}\searrow\lambda^{(1)}
       \searrow\cdots\searrow\lambda^{(n)}=\varnothing\right\}.
\end{aligned}
\]
These paths are in bijection with the standard Young tableaux of shape
$\lambda$.

For symplectic groups, Fourier--Jacobi descent plays a role analogous to
restriction in Young's branching rule. Let
$\sigma\in\Irr(\Sp_{2M}(\Fq))$. For $0<\ell\le M$ and $a\in\{\pm1\}$, consider the
rational nilpotent orbit (see Section~\ref{subsec:GGP})
\[
 \CO_{2\ell,a}=[(2\ell,a),1^{2M-2\ell}].
\]
The \emph{Fourier--Jacobi descent} of $\sigma$ is the Whittaker model
\cite[Section~4.1]{wang}
\[
 \cD^{\FJ}_{2\ell,a}(\sigma)
 :=\Wh_{\CO_{2\ell,a}}(\sigma).
\]
Since the stabilizer of the Whittaker datum contains a symplectic factor
$\Sp_{2M-2\ell}(\Fq)$ acting naturally on $\cD^{\FJ}_{2\ell,a}(\sigma)$, this
descent is a representation of $\Sp_{2M-2\ell}(\Fq)$, and we may iterate the procedure
on each of its irreducible constituents.
For $\sigma'\in\Irr(\Sp_{2M-2\ell}(\Fq))$, let
$m_{2\ell,a}(\sigma,\sigma')$ be the multiplicity of $\sigma'$ in
$\cD^{\FJ}_{2\ell,a}(\sigma)$. Unlike the multiplicities in Young's rule,
$m_{2\ell,a}(\sigma,\sigma')$ can be greater than one. By the GGP branching
formulas, $m_{2\ell,a}(\sigma,\sigma')$ is determined by the parameters of
$\sigma$ and $\sigma'$ in Lusztig's parametrization.

Fix $\pi\in\Irr(\Sp_{2N}(\Fq))$ and a list of descent parameters
\[
 \Gamma=\llbracket(2\ell_1,a_1),\ldots,(2\ell_r,a_r)\rrbracket,
 \qquad \sum_j\ell_j=N.
\]
Here $\ell_j\ge0$ and $a_j\in\{\pm1\}$. The analogue of a downward path
in the Young graph is a descent path of irreducible representations
\[
 {\mathcal P:\quad}
 \pi=\pi^{(0)}\searrow_{2\ell_1,a_1}\pi^{(1)}
 \searrow_{2\ell_2,a_2}\cdots
 \searrow_{2\ell_r,a_r}\pi^{(r)}=\mathbf{1}_{\Sp_0},
\]
where $\mathbf{1}_{\Sp_0}$ is the trivial representation,
$\pi^{(j)}\in\Irr(\Sp_{2N_j}(\Fq))$,
$N_j=N-\sum_{i\le j}\ell_i$, and the $j$th arrow means that
$\pi^{(j)}$ is an irreducible constituent of
$\cD^{\FJ}_{2\ell_j,a_j}(\pi^{(j-1)})$.
The list $\Gamma$ is the \emph{descent sequence index} of
$\mathcal P$; see Section~\ref{subsec:good-indexes}.   Decomposing at each step, we set
\[
 m_\Gamma(\pi):=\dim\cD^{\FJ}_{2\ell_r,a_r}\cdots
   \cD^{\FJ}_{2\ell_1,a_1}(\pi)
 =\sum_{\mathcal P}\prod_{j=1}^r
   m_{2\ell_j,a_j}\bigl(\pi^{(j-1)},\pi^{(j)}\bigr),
\]
where $\mathcal P$ runs over the descent paths with index
$\Gamma$. It is the analogue of the formula for $\dim V_\lambda$.

The composition laws of \cite[Section~5]{wang} reduce Whittaker multiplicities
of rational orbits to iterated descents. Represent a rational orbit $\CO$ in
$\mathfrak{sp}_{2N}(\Fq)$ by its signed Young diagram
$\CO=[(\lambda_1^{n_1},a_1),\ldots,(\lambda_t^{n_t},a_t)]$ with
$\lambda_1>\cdots>\lambda_t>0$ (see Section~\ref{subsec:wavefront-def}, \eqref{eq:signed-YD}),
where each even part carries a sign in $\{\pm1\}$ and each odd part has even
multiplicity. If the largest part $\lambda_1=2\ell$ is even, splitting the block
$((2\ell)^{n_1},a_1)$ via \eqref{eq:equal-even-combination} if necessary and removing a
row $(2\ell,a)$ leaves a rational orbit $\CO'$ in $\mathfrak{sp}_{2N-2\ell}(\Fq)$.
By Proposition~\ref{prop:even-composition-law},
\begin{equation*}
 m_\CO(\pi)
 =\dim\Wh_{\CO'}\!\left(\cD^{\FJ}_{2\ell,a}(\pi)\right).
\end{equation*}
If the largest part $\lambda_1=2\ell+1$ is odd, removing one pair
$(2\ell+1)^2$ leaves a rational orbit $\CO'$ in $\mathfrak{sp}_{2N-4\ell-2}(\Fq)$. By Proposition~\ref{prop:odd-pair-identities},
for each $a\in\{\pm1\}$,
\begin{equation}\label{eq:intro-odd-composition}
\begin{aligned}
 2m_\CO(\pi)
 ={}&\dim\Wh_{\CO'}\!\left(
       \cD^{\FJ}_{2\ell+2,\delta a}
       \cD^{\FJ}_{2\ell,a}(\pi)\right)\\
 &+\frac{q+1}{2}\dim\Wh_{\CO'}\!\left(
       \cD^{\FJ}_{2\ell,\delta a}
       \cD^{\FJ}_{2\ell+2,a}(\pi)\right)\\
 &+\frac{q-1}{2}\dim\Wh_{\CO'}\!\left(
       \cD^{\FJ}_{2\ell,-\delta a}
       \cD^{\FJ}_{2\ell+2,-a}(\pi)\right).
\end{aligned}
\end{equation}
Iterating these two identities, we obtain
\begin{equation}\label{eq:intro-orbit-descent-expansion}
\begin{aligned}
 m_\CO(\pi)
 &=2^{-s_\CO}\sum_{\mathbf c\in\{0,+,-\}^{s_\CO}}
   \left(\prod_{j=1}^{s_\CO}w_q(c_j)\right)
   m_{\Gamma_{\CO,\mathbf c}}(\pi),\\
 w_q(0)&=1,\qquad
 w_q(+)=\frac{q+1}{2},\qquad
 w_q(-)=\frac{q-1}{2}.
\end{aligned}
\end{equation}
Here $s_\CO$ is the number of pairs of equal odd parts of
$\CO$, and $\Gamma_{\CO,\mathbf c}$ is obtained by replacing the $j$th such
pair by the two descents selected by $c_j$ in
\eqref{eq:intro-odd-composition}; see \eqref{eq:odd-channel-table}. It
therefore suffices to study the descent multiplicities $m_\Gamma(\pi_q)$.

To study
$m_\Gamma(\pi_q)=\sum_{\mathcal P}\prod_{j=1}^r m_{2\ell_j,a_j}(\pi^{(j-1)},\pi^{(j)})$
as a function of $q$, we evaluate this sum over descent paths
\[\mathcal P\colon \pi_q=\pi^{(0)}\searrow_{2\ell_1,a_1}\pi^{(1)}\searrow_{2\ell_2,a_2}\cdots\searrow_{2\ell_r,a_r}\pi^{(r)}=\mathbf{1}_{\Sp_0}\]
to determine its degree and leading coefficient, and then relate them to stable orbit
dimensions. The key mechanism is to separate the Lusztig parameters of each $\pi^{(j)}$
into \emph{inherited} and \emph{new} parts. By the finite GGP branching rule
(\cite{WangJordanCorrected}; see Theorem~\ref{thm:FJ-branching}), the multiplicity
\[
m_{2\ell,a}(\pi^{(j-1)},\pi^{(j)}) \textrm{ of the step $\pi^{(j-1)}\searrow_{2\ell_j,a_j}\pi^{(j)}$}
\]
is
computed by matching the primary partition of $\pi^{(j-1)}$ at $[x]$ with that of $\pi^{(j)}$
at $[-x]$ (together with the unipotent symbols at $\pm1$). Any class
$[y]\ne[\pm1]$ appearing in $\pi^{(j)}$ whose negative $[-y]$ does not occur in $\pi^{(j-1)}$
is called \emph{new}, while the remaining parameters of $\pi^{(j)}$ are called
\emph{inherited} (see \eqref{eq:inherited-primary-partition}). By the GGP branching rule,
the possible choices for the inherited parameters of $\pi^{(j)}$ and the resulting step
multiplicities $m_{2\ell_j,a_j}(\pi^{(j-1)},\pi^{(j)})$ are independent of $q$ and depend only on
the type $\mathfrak t$ of $\pi_q$ (Section~\ref{subsec:type-transitions}). By contrast,
for every new class $[y]$, the branching rule forces the partition of $\pi^{(j)}$ at $[y]$
to be a single column $(1^e)$. By Proposition~\ref{prop:primary-counting}
(see \eqref{eq:AB-leading}), the number of classes of fixed degree $d$ and kind
$\kappa\in\{\GL,\U\}$ is a polynomial in $q$ of the form $q^d/(2d)+O(q^{d-1})$, so the
number of choices of new classes with given degrees, kinds and partitions is also a
polynomial in $q$. Grouping descent paths $\mathcal P$ by the types of $\pi^{(j)}$ and the
new/inherited data therefore yields a polynomial $m_{\Gamma,\mathfrak t}(X)$ with
$m_\Gamma(\pi_q)=m_{\Gamma,\mathfrak t}(q)$ (Proposition~\ref{prop:exact-poly}), and
\eqref{eq:intro-orbit-descent-expansion} gives the polynomiality of $m_\CO(\pi_q)$.

For each $\pi^{(j)}\in\Irr(\Sp_{2N_j}(\Fq))$, the group rank $N_j$ is the sum of the
ranks of its Lusztig parameters: a class $[x]$ of degree $d$ with partition $\lambda$
contributes $d|\lambda|$, and the unipotent symbols $\Lambda_{\pm1}$ contribute their
symbol ranks (Section~\ref{sec:parameters}). Separating the Lusztig parameters of
$\pi^{(j)}$ into inherited parameters and new classes directly determines the
degree and leading coefficient of $m_{\Gamma,\mathfrak t}(X)$: only the choices of new
classes contribute to the degree in $q$, while the rank occupied by the inherited
parameters of $\pi^{(j)}$ limits the rank available for new classes. Since a new class of degree $d$ with partition $(1^e)$ uses rank $de$ while contributing
$d$ to the degree in $q$, maximizing the degree in $q$ forces two conditions at every
intermediate stage $1\le j<r$: every new class in $\pi^{(j)}$ carries the partition
$(1)$, and the inherited parameters of $\pi^{(j)}$ have the smallest possible rank,
which is attained along successive \emph{first descents} (the maximal-step nonzero
descents at each stage; Section~\ref{subsec:GGP} and Proposition~\ref{prop:stagewise}).
Writing $N_j=N-\sum_{i\le j}\ell_i$ and letting $r_j^{\rm fd}(\mathfrak t)$ be the
rank after $j$ successive first descents, we obtain
\[
 \deg m_{\Gamma,\mathfrak t}(X)\le B_{\Gamma,\mathfrak t}
 :=\sum_{j=1}^{r-1}\bigl(N_j-r_j^{\rm fd}(\mathfrak t)\bigr)
 \qquad(m_{\Gamma,\mathfrak t}(X)\ne0),
\]
where the contribution of paths achieves degree $B_{\Gamma,\mathfrak t}$ precisely when
these two conditions hold
(Proposition~\ref{prop:degree-budget}). For the leading coefficient, each $1\le j<r$
with $N_j>r_j^{\rm fd}(\mathfrak t)$ contributes a factor of $1/2$ because the sign $a_j$ fixes
the parity of the unitary rank, reducing the leading coefficient to a count of
first-descent paths (Proposition~\ref{cor:cGamma-lc-explicit}).

To pass from $B_{\Gamma,\mathfrak t}$ to rational orbits, we compare good descent
sequence indices $\Gamma$ (each determining a rational orbit $\CO(\Gamma)$;
Section~\ref{subsec:good-indexes}) via elementary downward mutations (shifting $2$ from
$2\ell_i$ to $2\ell_{i+1}$; Section~\ref{subsec:mutations}), identifying
$B_{\Gamma,\mathfrak t}$ with
$\tfrac12(\dim\CO_{\pi_q}^{\rm st}-\dim\CO(\Gamma)^{\rm st})$
(Proposition~\ref{prop:endpoint-codim}).

In \eqref{eq:intro-orbit-descent-expansion},
every summand together with its weight $\prod_j w_q(c_j)$ satisfies the same degree
bound; for $\CO\in\ScO(\pi_q)$ this bound is attained at $\mathbf c=\mathbf0$ without
cancellation by positivity, proving \eqref{eq:main-degree}
(Theorem~\ref{prop:orbit-degree}). Grouping the highest-degree paths by the maximal
rational orbit $\CO(\Gamma')\in\ScO(\pi_q)^{\max}$ attached to their first-descent
index $\Gamma'$ then yields \eqref{eq:intro-leading-direct}
(Theorem~\ref{prop:orbit-leading-coefficient}).

Finally, for $G=\OO(V)$ with $n=\dim V$, the finite-field Gomez--Zhu theorem
\cite[Theorem~8.6]{wang} (Section~\ref{sec:orth-theta}) expresses $m_\CO(\pi_q)$ as a
finite sum of symplectic Whittaker multiplicities of the theta lift of $\pi_q^\vee$ (or
$\pi_q^\vee\otimes\det_G$) at the lifted orbit
$\CO^\uparrow\subset\mathfrak{sp}_{2n}(\Fq)$ (see \eqref{eq:orth-first-column-lift}).
By the symplectic case, this sum is a polynomial in $q$ whose degree and leading
coefficient are governed by a single theta constituent of maximal stable wavefront orbit,
yielding Theorems~\ref{thm:mainA} and~\ref{thm:mainB} for orthogonal groups.

The paper is organized as follows. In Section~\ref{sec:parameters}, we introduce
the combinatorial data and Lusztig's parametrization as $q$ varies. In
Section~\ref{sec:descent}, we recall rational nilpotent orbits, the finite GGP branching
rule, and the composition laws. In Section~\ref{sec:minimality-degree}, we establish the
rank inequalities comparing arbitrary descents with successive first descents. In
Section~\ref{sec:poly}, we prove the polynomiality of $m_\Gamma(\pi_q)$, determine its
degree and leading coefficient, and relate the degree to stable orbit dimensions via
sequence mutations. In Section~\ref{sec:orbit-multiplicities}, we combine these results
with the composition laws to prove Theorems~\ref{thm:mainA} and~\ref{thm:mainB} for
symplectic groups, and in Section~\ref{sec:orth-theta} we treat orthogonal groups via the
theta correspondence. Finally, in \ref{app:symbol}, we prove the symbol interlacing
inequalities used in Section~\ref{sec:minimality-degree}.

\section*{Acknowledgements}
Z.W. acknowledges support from the National Natural Science Foundation of
China (Grant No.~12571012) and the Jilin Provincial Postdoctoral Foundation.
P.L. acknowledges support from the New Cornerstone Science Foundation through
the New Cornerstone Investigator Program awarded to Professor Xuhua He.

\noindent\emph{Use of AI tools.} We used GPT and Gemini to polish the exposition, simplify several proofs, and find and fix sign errors in an early draft. The mathematical ideas and proofs are our own, and we
take full responsibility for the final manuscript.

\section{Combinatorial data and Lusztig's parametrization}\label{sec:parameters}
Throughout, for an $F$-stable algebraic group $\bH$ we put $H=\bH(\Fq)=\bH^F$. The
discriminant of a quadratic space $V$ of dimension $m$ is
\begin{equation}\label{eq:disc-def}
 \disc(V)=(-1)^{m(m-1)/2}\det(V)
 \in \Fq^\times/\Fq^{\times2}\simeq\{\pm\},
\end{equation}
where $\det(V)$ is the determinant of a Gram matrix of $V$.
We write $\OO(V)=\OO_m^{\disc(V)}$; with this convention $\OO_{2n}^{+}$ is the
split even orthogonal group.

\subsection{Partitions and Lusztig symbols}\label{subsec:symbols}
By a \emph{partition} we mean a sequence $\lambda=(\lambda_1,\lambda_2,\ldots)$ of
nonnegative integers with $\lambda_1\ge\lambda_2\ge\cdots$ and only finitely many
nonzero terms; sequences that differ only by trailing zeros are identified. We
write $|\lambda|=\sum_i\lambda_i$ for the size of $\lambda$,
$\ell(\lambda)=\#\{i:\lambda_i>0\}$ for its length, and $\lambda\vdash n$ if
$|\lambda|=n$. The transpose $\lambda^t$ is the partition with
$\lambda^t_i=\#\{j:\lambda_j\ge i\}$. If
$|\lambda|=|\mu|$, we write $\lambda\le\mu$ when
$\lambda_1+\cdots+\lambda_i\le\mu_1+\cdots+\mu_i$ for all $i\ge1$ (the dominance
order). We also use the exponential
notation $(r_1^{m_1},\ldots,r_s^{m_s})$, with $r_1>\cdots>r_s>0$, for the
partition in which $r_j$ occurs $m_j$ times.
Set
\begin{equation}\label{eq:first-column}
 \lambda_-:=
 \bigl(\max\{\lambda_1-1,0\},\max\{\lambda_2-1,0\},\ldots\bigr).
\end{equation}
Two partitions $\lambda,\mu$ are \emph{close} if $|\lambda_i-\mu_i|\le1$ for every $i$. The
common partition $\lambda\cap\mu$ consists of the positive entries for which
$\lambda_i=\mu_i$, counted with multiplicity. For
$\nu=(r_1^{m_1},\ldots,r_s^{m_s})$ with $r_1>\cdots>r_s>0$, define
\begin{equation*}
f(\nu):=\prod_{j=1}^s(m_j+1),\qquad f(\varnothing):=1.
\end{equation*}
The partitions $\lambda$ and $\mu$ are \emph{$2$-transverse} if they are close and every part
of $\lambda\cap\mu$ occurs with even multiplicity. For partitions $P,Q$, write
\begin{equation*}
P\preccurlyeq Q
 \quad\Longleftrightarrow\quad
Q_i-1\le P_i\le Q_i
\end{equation*}
for every $i$.
These definitions imply
\begin{equation*}
\begin{aligned}
 \lambda,\mu\text{ close}
 &\ \Longrightarrow\ \lambda_-,\mu_-\text{ close},\\
 \lambda,\mu\text{ $2$-transverse}
 &\ \Longrightarrow\ \lambda_-,\mu_-\text{ $2$-transverse},\\
P\preccurlyeq Q
 &\ \Longrightarrow\ P_-\preccurlyeq Q_-.
\end{aligned}
\end{equation*}

A Lusztig symbol is an equivalence class represented by
\[
 \Lambda=
 \begin{pmatrix}
a_1>\cdots>a_r\\
b_1>\cdots>b_s
 \end{pmatrix}
 =
 \begin{pmatrix}A\\ B\end{pmatrix},
 \qquad
A=\{a_1,\ldots,a_r\},\quad B=\{b_1,\ldots,b_s\}.
\]
Here $A,B$ are finite subsets of $\mathbb Z_{\ge0}$, and the equivalence relation is generated by the simultaneous shift
\[
 \begin{pmatrix}A\\ B\end{pmatrix}
 \sim
 \begin{pmatrix}
 \{a+1:a\in A\}\cup\{0\}\\
 \{b+1:b\in B\}\cup\{0\}
 \end{pmatrix}.
\]
The transpose of a symbol interchanges its rows:
$\Lambda^t=\begin{pmatrix}B\\A\end{pmatrix}$.
For $\Lambda=\begin{pmatrix}A\\ B\end{pmatrix}$ with $|A|=r$ and $|B|=s$, the defect and rank are
\begin{equation}\label{eq:symbol-rank}
 \defect(\Lambda)=r-s,
 \qquad
 \rank(\Lambda)=\sum_i a_i+\sum_j b_j-
 \left\lfloor\left(\frac{r+s-1}{2}\right)^2\right\rfloor.
\end{equation}
The associated bipartition is
\begin{equation}\label{eq:upsilon}
 \Upsilon(\Lambda)=
 \begin{pmatrix}
 \Upsilon(\Lambda)^+\\
 \Upsilon(\Lambda)^-
 \end{pmatrix},
 \qquad
 \Upsilon(\Lambda)^+_i=a_i-(r-i),\qquad
 \Upsilon(\Lambda)^-_j=b_j-(s-j).
\end{equation}
The rank, defect, and associated bipartition $\Upsilon$ are invariant under shifts.
In row indices, $+$ and $-$ are identified with $+1$ and $-1$, respectively. Write $|\Upsilon(\Lambda)|:=|\Upsilon(\Lambda)^+|+|\Upsilon(\Lambda)^-|$. Then
\begin{equation*}
 \rank(\Lambda)=|\Upsilon(\Lambda)|+c({\defect\Lambda}),
 \qquad
c({x})=\begin{cases}({x}^2-1)/4,&{x}\text{ odd},\\
                         {x}^2/4,&{x}\text{ even}.
          \end{cases}
\end{equation*}
The Alvis--Curtis dual symbol ${}^d\Lambda$ satisfies, for $\epsilon\in\{+,-\}$,
\begin{equation}\label{eq:AC-symbol}
 \defect({}^d\Lambda)=\defect(\Lambda),
 \qquad
 \Upsilon({}^d\Lambda)^\epsilon
   =\bigl(\Upsilon(\Lambda)^{-\epsilon}\bigr)^t.
\end{equation}

\subsection{Lusztig's parametrization}
\label{sec:classification}
Following the conventions of \cite{wang}, we take the dual groups to be
\begin{equation*}
\begin{array}{c|c}
\bG & \bG^* \\\hline
 \bSp_{2n} & \bSO(V^*),\quad \dim V^*=2n+1,\\
 \bO(V),\quad \dim V=2n+1 & \bSp(V^*)\times\{\pm1\},\quad \dim V^*=2n,\\
 \bO(V),\quad \dim V=2n & \bO(V^*),\quad V^*=V.
\end{array}
\end{equation*}
Let $G^*=(\bG^*)^F$. For a semisimple $s\in(\bG^*)^{\circ F}$, let $\cE(G,s)\subseteq\Irr(G)$ be the
Lusztig series attached to the $G^*$-conjugacy class of $s$. By Deligne--Lusztig theory and Lusztig's
classification \cite{DL76,Lusztig77,Lusztig81}, we have
\begin{equation*}
 \Irr(G)
 =\coprod_{(s)\in (\bG^*)_{ss}^{\circ F}/\operatorname{Ad}(G^*)}\cE(G,s).
\end{equation*}
Here $(\bG^*)_{ss}^{\circ F}$ denotes the set of semisimple elements of $(\bG^*)^{\circ F}$.

We use Lusztig's Jordan decomposition as in
\cite[Section~3.1]{WangJordanCorrected}. For a semisimple
$s\in(\bG^*)^{\circ F}$, write

 \begin{equation*}
 \CL_s:\cE(G,s)\xrightarrow{\sim}\cE(C_{G^*}(s),1)
\end{equation*}
for this correspondence, where $C_{G^*}(s):=C_{\bG^*}(s)^F$.
For $G=\OO_{2n+1}(\Fq)$, the group $C_{G^*}(s)$ has an extra factor $\{\pm1\}$ (see \eqref{eq:centralizer-primary}), and the right-hand side also involves a character $\epsilon_0\in\{\pm\}$ of this factor.
We call the semisimple parameter together with the labelled unipotent factors the \emph{{Lusztig parameters}} of $\pi$.

For $a\in\overline{\Fq}^\times$, the \emph{class} of $a$ is its orbit under
Frobenius and inversion:
\begin{equation*}
 [a]:=\{a^{q^j},a^{-q^j}:j\in\mathbb Z\}.
\end{equation*}
A class $[a]$ occurs in $s$ if it contains an eigenvalue of $s$
on the standard representation.
Let $G^*_{[a]}(s)$ be the corresponding factor of $C_{G^*}(s)$, called
the \emph{$[a]$-primary factor}, and define
$G^*_{[\ne\pm1]}(s):=\prod_{[a]\ne[\pm1]}G^*_{[a]}(s)$. Then
\begin{equation}\label{eq:centralizer-primary}
C_{G^*}(s)\simeq
G^*_{[\ne\pm1]}(s)\times G^*_{[1]}(s)\times G^*_{[-1]}(s)
 \times\begin{cases}\{\pm1\},&G=\OO_{2n+1}(\Fq),\\1,&\text{otherwise}.
 \end{cases}
\end{equation}
For $[a]\ne[\pm1]$, the finite group $G^*_{[a]}(s)$ is of general linear or
unitary type. Its unipotent irreducible representations are parametrized by
partitions \cite{LS77}. We denote the
corresponding unipotent factor by $\pi[a]$, its partition by $\lambda[a]$, and write
\begin{equation}\label{eq:rho-def}
 \rho=\bigotimes_{[a]\ne[\pm1]}\pi[a].
\end{equation}
If $[a]$ does not occur in the semisimple parameter, set
$\lambda[a]=\varnothing$.

For a finite symplectic, special odd orthogonal or full even orthogonal
group $H$, write $\mathcal S(H)$ for the symbols parametrizing
$\cE(H,1)$. For $\pi\in\cE(G,s)$ and $\xi=\pm1$, denote the symbol of the
$G^*_{[\xi]}(s)$-factor of
$\CL_s(\pi)$ by $\Lambda_\xi\in\mathcal S(G^*_{[\xi]}(s))$.
If $G^*_{[\xi]}(s)\simeq\Sp_{2m}(\Fq)$ or $\SO_{2m+1}(\Fq)$, then
\[
 \rank(\Lambda_\xi)=m,\qquad \defect(\Lambda_\xi)\equiv1\pmod4,
\]
whereas, if $G^*_{[\xi]}(s)\simeq\OO_{2m}^{{\delta_\xi}}(\Fq)$, then
\[
 \rank(\Lambda_\xi)=m,\qquad
 \defect(\Lambda_\xi)\equiv
 \begin{cases}0\pmod4,&{\delta_\xi}=+,\\2\pmod4,&{\delta_\xi}=-.
 \end{cases}
\]
Thus
\begin{equation*}
 \CL_s(\pi)=(\rho,\Lambda_1,\Lambda_{-1},\epsilon),
\end{equation*}
where $\epsilon$ is the scalar by which $-I$ acts in the odd orthogonal case
and is omitted otherwise. We write
\begin{equation*}
 \pi=\pi^s_{\rho,\Lambda_1,\Lambda_{-1},\epsilon},
\end{equation*}
omitting $s$ when clear and $\epsilon$ outside the odd orthogonal case.

For every class $[a]\ne[\pm1]$, the cardinality $\#[a]$ is even. Define
\begin{equation*}
d([a]):=\frac{\#[a]}{2}.
\end{equation*}
Set $\kappa([a])=\U$ if $a^{-1}=a^{q^j}$ for some $j\ge0$, and set
$\kappa([a])=\GL$ otherwise; we call $d([a])$ the \emph{degree} and
$\kappa([a])$ the \emph{kind} of $[a]$. Both $d([a])$ and $\kappa([a])$ are unchanged under
$[a]\mapsto[-a]$. If $[a]$ occurs in $s$, there is a unique integer $m_a\ge1$ such that
\begin{equation*}
G^*_{[a]}(s)\simeq
 \begin{cases}
  \GL_{m_a}(\mathbb F_{q^{d([a])}}),&\kappa([a])=\GL,\\
  \U_{m_a}(\mathbb F_{q^{d([a])}}),&\kappa([a])=\U.
 \end{cases}
\end{equation*}
The partition $\lambda[a]\vdash m_a$ contributes
\begin{equation*}
d([a])|\lambda[a]|
\end{equation*}
to the rank, and the sum of the eigenspaces indexed by $[a]$ has dimension
$2d([a])|\lambda[a]|$.
The primary sign used below is
\begin{equation}\label{eq:epsilon-rho}
 \epsilon_\rho:=(-1)^{\sum_{[a]:\,\kappa([a])=\U}|\lambda[a]|},
\end{equation}
as in \cite[Section~3.1]{WangJordanCorrected}; an empty sum is zero.

\subsection{Lusztig's parametrization as \texorpdfstring{$q$}{q} varies}\label{sec:qsystems}
Fix a symplectic or orthogonal group type and a rank $N$. Throughout, $p$ and $q$ satisfy our standing assumptions (see \eqref{eq:char-hyp} and $q\ge q_0$). For each such $q$, let $\bG$ be the corresponding algebraic
group over $\Fq$, with Frobenius endomorphism $F$. We fix one of the following
families of finite groups:
\[
\bG(\Fq)=\Sp_{2N}(\Fq),\qquad
\bG(\Fq)=\OO_{2N+1}(\Fq),\qquad
\bG(\Fq)=\OO_{2N}^{\epsilon_G}(\Fq),
\]
where $\epsilon_G\in\{\pm\}$ is fixed in the even orthogonal case.
For odd orthogonal groups we also fix the
discriminant of the defining quadratic space as $q$ varies.
Choose a finite set $I$ and
combinatorial data
\begin{equation}\label{eq:qtype}
 \mathfrak t=
 \left(
   \{(d_i,\kappa_i,\lambda_i)\}_{i\in I};
   \Lambda_1,\Lambda_{-1};\epsilon_0
 \right),
 \qquad \kappa_i\in\{\GL,\U\}.
\end{equation}
Here $d_i\ge1$ and each $\lambda_i$ is a nonempty partition. We regard
$\{(d_i,\kappa_i,\lambda_i)\}_{i\in I}$ as a multiset.
Let $r_{\pm1}=\rank(\Lambda_{\pm1})$. The symbols $\Lambda_1$ and
$\Lambda_{-1}$ parametrize unipotent representations of the following groups:
\[
\begin{array}{c|cc}
\bG(\Fq) & \Lambda_1 & \Lambda_{-1}\\ \hline
 \Sp_{2N}(\mathbb F_q)
   & \SO_{2r_1+1}(\Fq) & \OO_{2r_{-1}}^{\delta_{-1}}(\Fq)\\
 \OO_{2N+1}(\mathbb F_q)
   & \Sp_{2r_1}(\Fq) & \Sp_{2r_{-1}}(\Fq)\\
 \OO_{2N}^{\epsilon_G}(\mathbb F_q)
   & \OO_{2r_1}^{\delta_1}(\Fq) & \OO_{2r_{-1}}^{\delta_{-1}}(\Fq),
\end{array}
\]
where
\[
\begin{aligned}
 \delta_{\pm1}=+&\Longleftrightarrow\defect(\Lambda_{\pm1})\equiv0\pmod4,\\
 \delta_{\pm1}=-&\Longleftrightarrow\defect(\Lambda_{\pm1})\equiv2\pmod4,
\end{aligned}
\]
and \begin{equation*}
 \epsilon_G=\delta_1\delta_{-1}
 (-1)^{\sum_{i:\,\kappa_i=\U}|\lambda_i|}.
\end{equation*}
The sign $\epsilon_0\in\{\pm\}$ occurs only in the odd orthogonal case. We impose
\begin{equation*}
N=\sum_{i\in I}d_i|\lambda_i|+\rank(\Lambda_1)+\rank(\Lambda_{-1}).
\end{equation*}

A \emph{realization of $\mathfrak t$} over $\mathbb F_q$ is a pair
$(s,\pi)$ with
$\pi=\pi^s_{\rho,\Lambda_1,\Lambda_{-1},\epsilon_0}\in\cE(G,s)$ satisfying
the following condition. The classes $[a]\ne[\pm1]$ occurring in $s$
admit an indexing $[a_i]$, $i\in I$, such that
$\rho=\bigotimes_{i\in I}\pi[a_i]$ and
\begin{equation}\label{eq:qtype-primary}
 \bigl(d([a_i]),\kappa([a_i]),\lambda[a_i]\bigr)
   =(d_i,\kappa_i,\lambda_i)\qquad(i\in I).
\end{equation}

Let
\begin{equation*}
 \Irr_{\mathfrak t}(\bG(\Fq))
 :=\bigl\{\pi\in\Irr(\bG(\Fq)):(s,\pi)\text{ is a realization of }\mathfrak t
 \text{ for some }s\bigr\}.
\end{equation*}
By Proposition~\ref{prop:primary-counting}, the number of classes $[a]\ne[\pm1]$ of a given degree and kind tends to infinity with $q$, so $\Irr_{\mathfrak t}(\bG(\Fq))$ is nonempty for all sufficiently large $q$. A \emph{$q$-system of type $\mathfrak t$} is a family $\boldsymbol\pi=\{\pi_q\}$, indexed by sufficiently large $q$, with $\pi_q\in\Irr_{\mathfrak t}(\bG(\Fq))$.

By the character-degree formulas in the Jordan decomposition and
the degree formulas for unipotent characters, for each type $\mathfrak t$ there is a
polynomial $D_{\mathfrak t}(X)\in\mathbb Q[X]$ such that
\begin{equation*}
 \pi(1)=D_{\mathfrak t}(q)
 \qquad\bigl(\pi\in\Irr_{\mathfrak t}(\bG(\Fq))\bigr).
\end{equation*}
For the connected case, see \cite[Section~7.A]{MalleSpath16};
for the even orthogonal case, see \cite[Proposition~6.2 and Section~7]{Vinroot20}.
For odd orthogonal groups, $\OO_{2N+1}(\Fq)
=\SO_{2N+1}(\Fq)\times\{\pm I\}$. The polynomial is independent of the individual
classes $[a]$ and of the component label in odd orthogonal dimension; see
\cite{Lusztig84,LS77,WangJordanCorrected}.

\section{Rational wavefront sets and symplectic descent}\label{sec:descent}

\subsection{Rational nilpotent orbits and generalized Whittaker models}
\label{subsec:wavefront-def}
Recall that $\Fg=\operatorname{Lie}(\bG)(\overline{\Fq})$ and that $\CN\subset\Fg$ is
its nilpotent cone; thus $\CN^F$ is the set of nilpotent elements of
$\operatorname{Lie}(\bG)(\Fq)$.
Fix a $\bG$-invariant nondegenerate symmetric pairing
$\langle\,,\,\rangle$ on $\Fg$, defined over $\Fq$.

\begingroup
Let $\CO\subset\CN^F$ be a rational nilpotent orbit and choose $X\in\CO$.
For $X\ne0$,
choose a homomorphism $\phi_X:\mathbf{SL}_2\to\bG^\circ$, defined over $\Fq$,
with $d\phi_X\bigl(\begin{smallmatrix}0&1\\0&0\end{smallmatrix}\bigr)=X$ \cite[Theorem~47]{McNinch05}.
Let $(X,H_0,Y)\subset\Fg^F$ be the associated $\mathfrak{sl}_2$-triple,
and let $\lambda_X(t):=\phi_X(\operatorname{diag}(t,t^{-1}))$ be the
associated cocharacter. For $X=0$, set $H_0=Y=0$ and take $\lambda_X$
to be the trivial cocharacter.

As in \cite[Section~3.1]{BDT20}, the cocharacter $\lambda_X$ determines a
$\mathbb Z$-grading
\[
 \Fg=\bigoplus_{i\in\mathbb Z}\Fg_i,\qquad
 \Fg_i:=\{Z\in\Fg:\operatorname{Ad}(\lambda_X(t))Z=t^iZ
 \text{ for all }t\in\overline{\Fq}^{\,\times}\}.
\]
Each $\Fg_i$ is $F$-stable. Moreover, $X\in\Fg_2$, $Y\in\Fg_{-2}$, and
$[H_0,Z]=iZ$ for $Z\in\Fg_i$.
\endgroup
For $x,y\in\Fg_1$, the pairing
\[
 (x,y)\longmapsto\langle Y,[x,y]\rangle
\]
is a nondegenerate symplectic form on $\Fg_1$. Choose an $F$-stable Lagrangian subspace
$\mathfrak l\subset\Fg_1$ and define
\[
 \Fg_{\geq1.5}:=\mathfrak l\oplus\bigoplus_{i\ge2}\Fg_i.
\]
Let $\bG_{\ge1}$ and $\bG_{\geq1.5}$ be the connected unipotent subgroups with Lie
algebras $\bigoplus_{i\ge1}\Fg_i$ and $\Fg_{\geq1.5}$, respectively, and put
$G_{\ge1}=\bG_{\ge1}^F$ and $G_{\geq1.5}=\bG_{\geq1.5}^F$. Let $\mathcal U$ be the
unipotent variety of $\bG$. We fix the $F$-equivariant Cayley isomorphism
\[
 \operatorname{Cay}:\mathcal U\longrightarrow\CN
\]
of \cite[Section~3]{wang}. Its restriction maps $\bG_{\geq1.5}$ into
$\Fg_{\geq1.5}$.

Define
\[
 \psi_{\CO}(u)=\psi\bigl(\langle Y,\operatorname{Cay}(u)\rangle\bigr),
 \qquad u\in G_{\geq1.5}.
\]
The generalized Gelfand--Graev representation and the generalized Whittaker model are
\begin{equation*}
 \Gamma_{\CO}:=\operatorname{Ind}_{G_{\geq1.5}}^{G}\psi_{\CO},
 \qquad
 \Wh_{\CO}(\pi):=\Hom_{G_{\geq1.5}}(\psi_{\CO},\pi).
\end{equation*}
For $\pi\in\Irr(G)$, we have
\begin{equation*}
m_{\CO}(\pi)=\dim\Hom_{G}(\pi,\Gamma_{\CO})=\dim\Wh_{\CO}(\pi).
\end{equation*}

Define $M_\CO:=C_G(H_0)$ and
$M_{\psi_\CO}:=\{m\in M_\CO:\operatorname{Ad}(m)Y=Y\}$.
As in \cite[Section~3]{wang}, the representation
$\operatorname{Ind}_{G_{\geq1.5}}^{G_{\ge1}}\psi_\CO$ has a canonical extension
$\omega_{\psi_\CO}$ to $M_{\psi_\CO}\ltimes G_{\ge1}$. This gives an
$M_{\psi_\CO}$-action on
\[
 \Wh_\CO(\pi)\simeq\Hom_{G_{\ge1}}(\omega_{\psi_\CO},\pi).
\]

The rational wavefront set is
\begin{equation*}
 \ScO(\pi):=\{\CO:m_{\CO}(\pi)\ne0\}.
\end{equation*}

\subsection{The rational orbit order}\label{subsec:rational-orbit-order}
The rational nilpotent orbits in $\CN^F$ are parametrized by the signed Young diagrams of
\cite[Section~3.1]{wang}. These have the form
\begin{equation}\label{eq:signed-YD}
 \begin{gathered}
 [ (\lambda_1^{n_1},a_1),\ldots,(\lambda_t^{n_t},a_t) ],
 \qquad \lambda_1>\cdots>\lambda_t>0,\\
n_i\ge1,\qquad a_i\in\Fq^\times/\Fq^{\times2}\simeq\{\pm\}.
 \end{gathered}
\end{equation}
They satisfy $\sum_i n_i\lambda_i=\dim V$ and the following conditions. If $V$ is
skew-symmetric, $a_i$ records the discriminant of the multiplicity space when $\lambda_i$ is
even, while for odd $\lambda_i$ one has $n_i$ even and $a_i=+$. If $V$ is symmetric, $a_i$ is the discriminant of the multiplicity space for odd $\lambda_i$; for even $\lambda_i$, the multiplicity $n_i$ is even
and $a_i=+$. The signs for symmetric $V$ also satisfy
\begin{equation}\label{eq:YD-disc}
 \prod_{k<j}\varepsilon_q(-1)^{n_kn_j}\prod_i a_i=\disc(V).
\end{equation}
where the square class of $-1$ is identified with the sign $\varepsilon_q(-1)$;
thus $\delta=\varepsilon_q(-1)$.

We label a rational orbit $\CO$ by its signed Young diagram
\eqref{eq:signed-YD}, omitting signs that are forced to be $+$. Its underlying
stable orbit is denoted $\CO^{\rm st}$, and $\Ost{\lambda}$ is the
stable orbit with partition $\lambda$.
In particular, odd blocks carry no sign in the symplectic case; in the
orthogonal case the even blocks carry no sign, and the odd blocks keep their
signs unless these are determined by \eqref{eq:YD-disc}.

Equal even blocks in a symplectic signed Young diagram are combined by
\begin{equation}\label{eq:equal-even-combination}
 [(\lambda,a_1),\ldots,(\lambda,a_m)]
 =[(\lambda^m,\varepsilon_q(-1)^{m(m-1)/2}a_1\cdots a_m)].
\end{equation}

We use the rational order of \cite[Section~3.2]{wang}. For a rational orbit
$\CO=[(\lambda_i^{n_i},a_i)]$, let $c_j(\CO)$ be the number of boxes in the $j$-th column of its
Young diagram and define
\begin{equation}\label{eq:rational-order-data}
 \begin{aligned}
C_r(\CO)&=\sum_{j=1}^r c_j(\CO),\\
 \operatorname{sgn}_r(\CO)
 &=\prod_{\lambda_i\ge r}a_i
   \prod_{\substack{i<j\\ \lambda_i,\lambda_j\ge r}}
   \varepsilon_q(-1)^{n_in_j}.
 \end{aligned}
\end{equation}
When $V$ is skew-symmetric, define $\CO\le_{\rm rat}\CO'$ for two rational orbits $\CO$ and $\CO'$ by
\begin{equation}\label{eq:rational-order-sp}
 \begin{aligned}
C_r(\CO)&\ge C_r(\CO') &&(r\ge1),\\
C_{2s-1}(\CO)=C_{2s-1}(\CO')
 &\Longrightarrow
 \operatorname{sgn}_{2s-1}(\CO)=\operatorname{sgn}_{2s-1}(\CO') &&(s\ge1).
 \end{aligned}
\end{equation}
For symmetric $V$, the first condition is unchanged, while the sign condition applies at even
columns $2s$ rather than odd columns $2s-1$. The set of maximal rational wavefront orbits is
\begin{equation*}
 \ScO(\pi)^{\max}
 :=\left\{\CO\in\ScO(\pi):
 \nexists\,\CO'\in\ScO(\pi)\text{ with }\CO<_{\rm rat}\CO'\right\}.
\end{equation*}

On stable nilpotent orbits, we use the closure order
\[
 \CO_1^{\rm st}\le\CO_2^{\rm st}
 \quad\Longleftrightarrow\quad
 \CO_1^{\rm st}\subseteq\overline{\CO_2^{\rm st}}.
\]
For irreducible $\pi$, the set of stable nilpotent
orbits
$\{\CO^{\rm st}:\CO\in\ScO(\pi)\}$ has a unique maximal element for
the closure order \cite{Lusztig92,GM00,AA07,Taylor16}, denoted by
$\CO_\pi^{\rm st}$.

As $q$ varies, we identify rational nilpotent orbits in $\CN^F$ for different $q$ that have the same $(C_r(\CO),\operatorname{sgn}_r(\CO))$ in \eqref{eq:rational-order-data} for every $r\ge1$; by abuse of notation, we still denote such an orbit by $\CO$. By \eqref{eq:rational-order-data} and \eqref{eq:rational-order-sp}, these invariants uniquely determine $\CO$ for each $q$ and preserve $\le_{\rm rat}$.

\subsection{Fourier--Jacobi descent}\label{subsec:GGP}
Let $\pi\in\Irr(\Sp_{2n}(\Fq))$. For $0<\ell\le n$ and $a\in\{\pm1\}$, consider
the rational nilpotent orbit
\begin{equation*}
 \CO_{2\ell,a}:=[(2\ell,a),1^{2n-2\ell}].
\end{equation*}
For $\CO_{2\ell,a}$ we have
$M_{\psi_{\CO_{2\ell,a}}}\simeq\OO_1(\Fq)\times\Sp_{2n-2\ell}(\Fq)$. Restricting the natural
$M_{\psi_{\CO_{2\ell,a}}}$-action to the symplectic factor, define
\begin{equation*}
 \cD^{\FJ}_{2\ell,a}(\pi):=\Wh_{\CO_{2\ell,a}}(\pi)
\end{equation*}
as an $\Sp_{2n-2\ell}(\Fq)$-representation.

For $b\in\{\pm1\}$, let $\omega_{n,\psi_b}$ be the Weil
representation of $\Sp_{2n}(\Fq)$ associated with
$\psi_b(x)=\psi(cx)$, where $c\in\Fq^\times$ satisfies
$\varepsilon_q(c)=b$. Its isomorphism class depends only on $b$.
For the formal index $\ell=0$ we use the convention
\begin{equation}\label{eq:FJ-zero}
 \cD^{\FJ}_{0,a}(\pi):=\pi\otimes\omega_{n,\psi_{\varepsilon_q(-1)a}},
\end{equation}
as in \cite[(4.17)]{wang}.
For $0\le\ell\le n$ and $\pi'\in\Irr(\Sp_{2n-2\ell}(\Fq))$, define
\begin{equation}\label{eq:GGP-mult}
m_{2\ell,a}(\pi,\pi')
 :=\dim\Hom_{\Sp_{2n-2\ell}(\Fq)}
   \bigl(\pi',\cD^{\FJ}_{2\ell,a}(\pi)\bigr).
\end{equation}
This is the multiplicity in the finite Gan--Gross--Prasad problem; we call it the
\emph{GGP multiplicity} \cite[Section~1.2]{Wang21}.
The \emph{first descent index} of $\pi$ is the largest
$\ell\in\{0,\ldots,n\}$ for which $\cD^{\FJ}_{2\ell,a}(\pi)\ne0$ for some $a\in\{\pm1\}$; at this value of
$\ell$, each nonzero $\cD^{\FJ}_{2\ell,a}(\pi)$ is called a \emph{first descent}.

\subsection{The finite GGP branching rule}\label{subsec:branching-input}
Let $0\le\ell\le n$, $a\in\{\pm1\}$, and write
\[
 \pi=\pi_{\rho,\Lambda_1,\Lambda_{-1}}\in\Irr(\Sp_{2n}(\Fq)),\qquad
 \pi'=\pi_{\rho',\Lambda'_1,\Lambda'_{-1}}
       \in\Irr(\Sp_{2n-2\ell}(\Fq)).
\]
Let $\delta_{-1},\delta'_{-1}\in\{\pm1\}$ be the discriminants of the symmetric
spaces defining the even orthogonal factors labelled by
$\Lambda_{-1},\Lambda'_{-1}$, respectively, as in \eqref{eq:disc-def}.
Equivalently,
\[
 \delta_{-1}=(-1)^{\defect(\Lambda_{-1})/2},\qquad
 \delta'_{-1}=(-1)^{\defect(\Lambda'_{-1})/2}.
\]
If either factor has rank zero, its discriminant is $+1$.
For any symbol $\Lambda$, write
\[
 \Lambda^{\langle+1\rangle}=\Lambda,\qquad
 \Lambda^{\langle-1\rangle}=\Lambda^t.
\]

The factors $\pi[x]$ of $\rho$ and $\pi'[-x]$ of $\rho'$, with
partitions $\lambda[x]$ and $\lambda'[-x]$, are as in \eqref{eq:rho-def}.
For partitions $\lambda,\nu$, define
\begin{equation}\label{eq:exact-primary-pairings}
\begin{aligned}
 b_{\GL}(\lambda,\nu)
 &=\#\{\tau:\tau\preccurlyeq\lambda,\ \tau\preccurlyeq\nu\}\\
 &=\begin{cases}f(\lambda\cap\nu),&\lambda,\nu\text{ are close},\\
                0,&\text{otherwise},\end{cases}\\
 b_{\U}(\lambda,\nu)
 &=\begin{cases}1,&\lambda,\nu\text{ are $2$-transverse},\\
                0,&\text{otherwise}.
   \end{cases}
\end{aligned}
\end{equation}
Here $f(r_1^{m_1},\ldots,r_s^{m_s})=\prod_j(m_j+1)$, and the
partition $\tau=\varnothing$ is allowed.

For a symplectic symbol
$\Lambda$ and a symbol $\Xi$ of an even orthogonal group
$\OO_{2m}^{\epsilon}(\Fq)$,
define $b_\epsilon(\Lambda,\Xi)\in\{0,1\}$ by
\begin{equation}\label{eq:GGP-symbol-rule}
 b_\epsilon(\Lambda,\Xi)=
 \begin{cases}
 1,&\begin{aligned}
    &\defect(\Lambda)+\defect(\Xi)=\epsilon,\\
    &\Upsilon(\Xi)^\epsilon\preccurlyeq\Upsilon(\Lambda)^{-\epsilon},\\
    &\Upsilon(\Lambda)^\epsilon\preccurlyeq\Upsilon(\Xi)^{-\epsilon},
    \end{aligned}\\[
4pt]
 0,&\text{otherwise}.
 \end{cases}
\end{equation}
\begin{theorem}[{\cite[Theorems~7.4(1) and 7.6(1)]{WangJordanCorrected}}]
\label{thm:FJ-branching}
Assume $q\ge q_0$. Let
\begin{equation}\label{eq:GGP-rational-sign}
 {t=-a\epsilon_\rho\epsilon_{\rho'}\delta_{-1}\delta'_{-1}}.
\end{equation}
Then
\begin{equation}\label{eq:GGP-primary-numerical}
\begin{aligned}
 m_{2\ell,a}(\pi,\pi')
 &=\prod_{[x]\ne[\pm1]}
 b_{\kappa([x])}\bigl(\lambda[x],\lambda'[-x]\bigr)\\
 &\quad\times
 b_{\delta'_{-1}}\!\left(\Lambda_1,(\Lambda'_{-1})^{\langle\delta t\rangle}\right)
 b_{\delta_{-1}}\!\left(\Lambda'_1,\Lambda_{-1}^{\langle t\rangle}\right).
\end{aligned}
\end{equation}
\end{theorem}
\begingroup
For later use, we write the product of the last two factors in
\eqref{eq:GGP-primary-numerical} as
\begin{equation}\label{eq:symbol-pairing-factor}
 \mathbf b_t(\Lambda_{\pm1};\Lambda'_{\pm1})
 :=b_{\delta'_{-1}}\!\left(\Lambda_1,(\Lambda'_{-1})^{\langle\delta t\rangle}\right)
   b_{\delta_{-1}}\!\left(\Lambda'_1,\Lambda_{-1}^{\langle t\rangle}\right)\in\{0,1\},
\end{equation}
where $\Lambda_{\pm1}=(\Lambda_1,\Lambda_{-1})$ and
$\Lambda'_{\pm1}=(\Lambda'_1,\Lambda'_{-1})$, and
$\delta_{-1}=(-1)^{\defect(\Lambda_{-1})/2}$, $\delta'_{-1}=(-1)^{\defect(\Lambda'_{-1})/2}$
as above.
\endgroup
\begingroup
\begin{remark}\label{rem:GGP-theta}
Let $D$ denote Alvis--Curtis duality, and let
$m_\theta(\sigma,\sigma')$ be the multiplicity of $\sigma\otimes\sigma'$
in the Weil representation \cite{Ger}.
For $\kappa\in\{\GL,\U\}$ and $d\ge1$, let $\pi_\lambda$ be the unipotent
representation of $\kappa_{|\lambda|}(\mathbb F_{q^d})$ labelled by $\lambda$.
Similarly, let $\pi_\Lambda^{\Sp}$ and $\pi_\Xi^{\OO}$ be the unipotent
representations of $\Sp_{2\rank\Lambda}(\Fq)$ and
$\OO_{2\rank\Xi}^{\epsilon}(\Fq)$ with symbols $\Lambda$ and $\Xi$.
For $\kappa\in\{\GL,\U\}$, $m_\theta(\pi_\lambda,\pi_\nu)$ is given in
\cite[Theorem~5.5 and Proposition~5.6]{AMR96}.
For symplectic--even-orthogonal pairs, the theta correspondence was
determined by Pan \cite[Theorem~1.8]{PanUnipotent} and also follows from
Ma--Qiu--Zou \cite[Theorem~1.4]{MaQiuZou}.
The connection between finite Gan--Gross--Prasad multiplicities and theta
correspondence after Alvis--Curtis duality was observed by the second author
\cite{WangGGP,WangJordanCorrected}. In our conventions, the factors in
\eqref{eq:GGP-primary-numerical} are given by
\cite[Theorem~7.6(1)]{WangJordanCorrected}:
\begin{equation*}
\begin{aligned}
 b_\kappa(\lambda,\nu)&=m_\theta(D\pi_\lambda,D\pi_\nu)=m(\pi_\lambda,\pi_\nu),\\
 b_\epsilon(\Lambda,\Xi)
 &=m_\theta\!\left(D\pi_\Lambda^{\Sp},
 D\pi_{\Xi^{\langle\delta^{\rank\Lambda}\rangle}}^{\OO}\right).
\end{aligned}
\end{equation*}
Here $m(\pi_\lambda,\pi_\nu)$ is the GGP multiplicity for $\kappa=\U$
\cite[Theorem~1.1]{LiuWangGGP21}; for $\kappa=\GL$, it is the branching
multiplicity defined in \cite[Section~6.5]{WangGGP}.
\end{remark}
\endgroup
Choose representatives of the dual symbols with both rows nonempty,
applying shifts if necessary. For $\eta\in\{\pm1\}$, let
$\partial_\eta\Lambda$ be the symbol
obtained by deleting a largest entry of $\Lambda$ from its $\eta$-row.
This operation is defined only when the $\eta$-row contains
a largest entry of $\Lambda$.
If a largest entry occurs in both rows, both $\partial_{+1}\Lambda$
and $\partial_{-1}\Lambda$ are defined.
\begin{proposition}\label{cor:first-descent-rule}
Let $\cD^{\FJ}_{2\ell,a}(\pi)$ be a first descent of $\pi$, with $\ell\ge0$.
It is multiplicity-free, and its irreducible constituents are precisely the
$\pi'$ satisfying
\begin{equation}\label{eq:first-primary}
 \lambda'[-x]=\lambda[x]_-\qquad([x]\ne[\pm1])
\end{equation}
and
\[
\begin{aligned}
 {}^d\Lambda'_1&=\bigl(\partial_{t\delta_{-1}}({}^d\Lambda_{-1})\bigr)^{\langle-t\rangle},\\
 {}^d\Lambda'_{-1}&=\bigl(\partial_{\delta'_{-1}}({}^d\Lambda_1)\bigr)^{\langle-\delta t\rangle},
\end{aligned}
\]
with $t$ as in \eqref{eq:GGP-rational-sign} and both deletions defined.
\end{proposition}
\begin{proof}
The description of the irreducible constituents of $\cD^{\FJ}_{2\ell,a}(\pi)$ and their multiplicities
follows from \cite[Theorem~6.1(ii) and Section~6.1]{wang}, with the signs
translated using \eqref{eq:GGP-rational-sign}. The case $\ell=0$ is included
by \cite[Remark~6.3]{wang}. The symbol part can also be derived
directly from Theorem~\ref{thm:FJ-branching}; see the proof of
Lemma~\ref{lem:symbol-first-descent-comparison}.
\end{proof}

For an irreducible constituent $\pi'$ of $\cD^{\FJ}_{2\ell,a}(\pi)$ as in Proposition~\ref{cor:first-descent-rule}, let $h[x]=\ell(\lambda[x])$.
By \eqref{eq:first-primary} and \eqref{eq:epsilon-rho}, we have
\[
 \epsilon_\rho\epsilon_{\rho'}=(-1)^{\sum_{\kappa([x])=\U}h[x]}.
\]

If $m_{2\ell,a}(\pi,\pi')\ne0$, then for each $[x]\ne[\pm1]$ such
that exactly one of $\lambda[x]$ and $\lambda'[-x]$ is nonempty, the nonempty
partition $\lambda'[-x]$ is of the form $(1^e)$, because by \eqref{eq:exact-primary-pairings}, we have
{}
\begin{equation}\label{eq:GGP-one-sided-column}
b_\kappa(\lambda,\varnothing)
 =b_\kappa(\varnothing,\lambda)
 =\begin{cases}1,&\lambda=(1^e)\text{ for some }e\ge0,\\
               0,&\text{otherwise}.\end{cases}
\end{equation}

\begin{remark}\label{rem:generic-new-classes}
Let $B$ consist of the classes $[b]$ in $\rho'$ for which
$\lambda[-b]=\varnothing$. Suppose that $\lambda'[b]=(1^{e_b})$ for
each $[b]\in B$, as is necessary for $m_{2\ell,a}(\pi,\pi')\ne0$
by \eqref{eq:GGP-one-sided-column}.
Put $d_b=d([b])$. According to the kind of $[b]$, let $H_b=\GL_{e_b}(\mathbb F_{q^{d_b}})$ or $H_b=\U_{e_b}(\mathbb F_{q^{d_b}})$. The corresponding $F$-stable Levi subgroup $\mathbf L$ of $\bG'=\bSp_{2(n-\ell)}$ satisfies
\[
 \mathbf L^F\simeq\Bigl(\prod_{[b]\in B}H_b\Bigr)\times\Sp_{2r}(\Fq),
 \qquad r=n-\ell-\sum_{[b]\in B}d_b e_b.
\]
Put $\varepsilon_{\mathbf H}=(-1)^{\rank_{\Fq}\mathbf H}$. Since the dual of
$\mathbf L$ contains the full centralizer of the semisimple parameter of $\pi'$,
the map $\varepsilon_{\bG'}\varepsilon_{\mathbf L}R_{\mathbf L}^{\bG'}$ is a bijection
between the corresponding Lusztig series \cite[Theorem~13.25(ii)]{DM91}. On
$H_b$, the unipotent label $(1^{e_b})$ corresponds to the Steinberg
character, so there is a $\sigma\in\Irr(\Sp_{2r}(\Fq))$ with
\[
 \chi_{\pi'}=\varepsilon_{\bG'}\varepsilon_{\mathbf L}\,
 R_{\mathbf L}^{\bG'}\!\left(
   \Bigl(\mathop{\boxtimes}_{[b]\in B}
      (\theta_b\otimes\operatorname{St}_{H_b})\Bigr)
       \boxtimes\chi_\sigma\right).
\]
Here $\theta_b$ is the linear character of $H_b$ given by the chosen
representative of $[b]$ in the dual Levi, so that $\theta_b\otimes\operatorname{St}_{H_b}$
is the unique generic character in its Lusztig series.
\end{remark}

\subsection{Descent sequences and good indices}\label{subsec:good-indexes}

Let $\pi\in\Irr(\Sp_{2N}(\Fq))$. For $0\le\ell\le N$,
$a\in\{\pm1\}$ and $\pi'\in\Irr(\Sp_{2(N-\ell)}(\Fq))$, write
\[
 \pi\searrow_{2\ell,a}\pi'\quad\Longleftrightarrow\quad
 m_{2\ell,a}(\pi,\pi')>0.
\]
We write $\pi\fdsearrow_{2\ell,a}\pi'$ when, in addition, $\ell$ is the first
descent index of $\pi$ (Section~\ref{subsec:GGP}). A finite chain of arrows
$\searrow$ is a \emph{descent path}, and a \emph{first-descent path} if all its
arrows are of the form $\fdsearrow$.

A \emph{descent sequence} of $\pi$ is a descent path that terminates on first
reaching $\Sp_0$:
\[
 \mathcal P:\quad
 \pi=\pi^{(0)}\searrow_{2\ell_1,a_1}\pi^{(1)}
 \searrow_{2\ell_2,a_2}\cdots
 \searrow_{2\ell_r,a_r}\pi^{(r)}=\mathbf{1}_{\Sp_0},
\]
where
\[
 \begin{gathered}
 N_0=N,\qquad N_j=N_{j-1}-\ell_j,\\
 \pi^{(j)}\in\Irr(\Sp_{2N_j}(\Fq)),\qquad
 \ell_j\ge0,\quad a_j\in\{\pm1\}\quad(1\le j\le r).
 \end{gathered}
\]
Here $\mathbf{1}_{\Sp_0}$ is the trivial representation of $\Sp_0(\Fq)$. For $N=0$ the sequence is empty.
To compare indices of different lengths, we may append arrows
$\mathbf{1}_{\Sp_0}\searrow_{0,a}\mathbf{1}_{\Sp_0}$ as in Section~\ref{subsec:mutations}.

Following \cite[Section~4.1]{wang}, the associated \emph{descent sequence
index} is
\begin{equation*}
 \Gamma=\llbracket(2\ell_1,a_1),\ldots,(2\ell_r,a_r)\rrbracket,
 \qquad \sum_{j=1}^r\ell_j=N.
\end{equation*}
We also use this term for any finite list of this form, whether
or not it arises from a descent sequence. An index is \emph{realized by $\pi$} if it comes from a descent sequence of
$\pi$. Here appended zero entries are allowed. The \emph{numerical sequence} of
$\Gamma$ is $\Gamma^{\rm num}=(2\ell_1,\ldots,2\ell_r)$. It need not be a partition.

Let $\mathfrak t^{(j)}$ be the type of $\pi^{(j)}$, and let $\mathfrak t_\varnothing$ be the type of $\mathbf{1}_{\Sp_0}$, the unique type of rank zero. The path $\mathcal P$ determines a sequence of types
\begin{equation*}
 \mathfrak t=\mathfrak t^{(0)}
 \searrow_{2\ell_1,a_1}\mathfrak t^{(1)}
 \searrow\cdots\searrow_{2\ell_r,a_r}\mathfrak t^{(r)}={\mathfrak t_\varnothing}.
\end{equation*}
In Section~\ref{subsec:type-transitions} we sort the descent paths by these types and by their step data.

We now attach rational orbits to good indices. Set $\ell_0=+\infty$ and
$\ell_{r+1}=0$. Following \cite[Definition~4.2]{wang}, an index $\Gamma$ is
\emph{good} if, for every $1\le i\le r$, one of the following holds:
\begin{equation}\label{eq:good-cases}
\begin{aligned}
\textup{(i)}\quad&\ell_{i-1}\ge\ell_i\ge\ell_{i+1},\\
\textup{(ii)}\quad&\ell_i=\ell_{i+1}-1<\ell_{i-1}{\ \text{ and }\ a_i=\varepsilon_q(-1)a_{i+1}},\\
\textup{(iii)}\quad&\ell_i=\ell_{i-1}+1>\ell_{i+1}{\ \text{ and }\ a_{i}=\varepsilon_q(-1)a_{i-1}} .
\end{aligned}
\end{equation}
The sign condition also applies when $\ell_i=0$. For a good
index, define
\[
 p_i=\begin{cases}
 (2\ell_i,a_i),&\ell_{i-1}\ge\ell_i\ge\ell_{i+1},\\
 2\ell_i+1,&\ell_i=\ell_{i+1}-1<\ell_{i-1},\\
 2\ell_i-1,&\ell_i=\ell_{i-1}+1>\ell_{i+1}.
 \end{cases}
\]
As in \cite[Section~4]{wang}, we set $\CO(\Gamma):=[p_1,\ldots,p_r]$:
equal parts are combined by
\eqref{eq:equal-even-combination}. This is the signed Young diagram of a rational
nilpotent orbit in $\mathfrak{sp}_{2\sum_i\ell_i}(\Fq)$. Its odd parts occur in pairs
$(2\ell_i+1)^2$, which we call \emph{odd pairs}. For a rational orbit $\CO$, put
\[
 \mathfrak G(\CO):=\{\Gamma\text{ good}:\CO(\Gamma)=\CO\}.
\]

\subsection{Mutations of rational orbits}\label{subsec:orbit-mutations}
The rational orbit order is generated by certain elementary moves on
signed Young diagrams, which we now recall (the analogous order on
descent sequence indices appears in Section~\ref{subsec:mutations}).

For symplectic groups, an \emph{orbit mutation} replaces some
blocks of a diagram according to one of the rules
in \eqref{eq:orbit-mutations} below; here $\lambda\ge\mu$ are positive even
integers and $a,b\in\{\pm\}$.
Repeated blocks may have to be split first (for even blocks by
\eqref{eq:equal-even-combination}), and after the replacement zero parts are
deleted, the rows are reordered and equal parts are combined again.

 \begin{equation}\label{eq:orbit-mutations}
\begin{array}{ll}
\textup{(1)}&(\lambda-1)^2\longmapsto(\lambda,a),(\lambda-2,\varepsilon_q(-1)a),\\[
2pt]
\textup{(2)}&(\lambda,a),(\mu,b)\longmapsto
       (\lambda+2,a),(\mu-2,b),\\[2pt]
\textup{(3)}&(\lambda,a),(\mu,b)\longmapsto
       (\lambda+2,-a),(\mu-2,-b),\\[2pt]
\textup{(4)}&(\lambda^2,a),(\mu,b)\longmapsto
       (\lambda+1)^2,(\mu-2,ab),\\[2pt]
\textup{(5)}&(\lambda,a),(\mu^2,b)\longmapsto
       (\lambda+2,ab),(\mu-1)^2,\\[2pt]
\textup{(6)}&(\lambda^2,a),(\mu^2,a)\longmapsto
       (\lambda+1)^2,(\mu-1)^2.
\end{array}
\end{equation}
Write $\CO\rightsquigarrow_{\rm orb}\CO'$ if $\CO'$ is obtained from $\CO$ by one
replacement. Then $\CO<_{\rm rat}\CO'$.
\subsection{Elementary sequence mutations}\label{subsec:mutations}
To compare two indices of different lengths, we add
entries $(0,a)$, with specified signs, to the shorter one. This does not change
the iterated descent, since $\cD^{\FJ}_{0,a}(\mathbf{1}_{\Sp_0})=\mathbf{1}_{\Sp_0}$.
An \emph{elementary downward mutation} changes two consecutive entries
of an index as follows, where $\ell_i>0$:
\begin{equation}\label{eq:elementary-seq-mutation}
 \begin{gathered}
 (2\ell_i,a_i),(2\ell_{i+1},a_{i+1})
 \longmapsto(2\ell_i-2,b_i),(2\ell_{i+1}+2,b_{i+1}),\\
b_i b_{i+1}=a_i a_{i+1};
 \end{gathered}
\end{equation}
here $2\ell_{i+1}$ and $2\ell_i-2$ may be zero. They are the inverses of
the elementary mutations of \cite[Section~4.1]{wang}.
We write $\Gamma'\le_{\rm seq}\Gamma$ if $\Gamma'$ is obtained from $\Gamma$ by a
finite chain
\[
 \Gamma=\Gamma_0\rightsquigarrow\Gamma_1\rightsquigarrow\cdots
 \rightsquigarrow\Gamma_s=\Gamma'
\]
of elementary downward mutations. The intermediate indices
need not be good. By \cite[Lemma~4.3]{wang}, for good indices $\Gamma,\Gamma'$ one has
\begin{equation}\label{eq:good-order-compatibility}
 \Gamma'\le_{\rm seq}\Gamma
 \quad\Longrightarrow\quad
 \CO(\Gamma')\le_{\rm rat}\CO(\Gamma),
\end{equation}
and, conversely, if $\CO\le_{\rm rat}\CO'$ and $\Gamma\in\mathfrak G(\CO)$, then there exists
$\Gamma'\in\mathfrak G(\CO')$ such that
\begin{equation}\label{eq:orbit-to-sequence-order}
 \Gamma\le_{\rm seq}\Gamma'.
\end{equation}
\subsection{Composition laws}\label{subsec:composition-laws}
\begingroup
Let
$
 \CO'=[(d_1^{m_1},a_1),\ldots,(d_s^{m_s},a_s)]
$
be a rational nilpotent orbit in $\mathfrak{sp}_{2m}(\Fq)$.
For $\ell>0$, $a\in\{\pm1\}$ and $2\ell\ge d_1$, define
\begin{equation*}
 [(2\ell,a)]*\CO'
 :=[(2\ell,a),(d_1^{m_1},a_1),\ldots,(d_s^{m_s},a_s)]
 \quad{\text{in }\mathfrak{sp}_{2(m+\ell)}(\Fq)}.
\end{equation*}
For $r\ge0$ and $2r+1\ge d_1$, define
\begin{equation*}
 [(2r+1)^2]*\CO'
 :=[(2r+1)^2,(d_1^{m_1},a_1),\ldots,(d_s^{m_s},a_s)]
 \quad{\text{in }\mathfrak{sp}_{2(m+2r+1)}(\Fq)}.
\end{equation*}
Here, we combine equal odd parts by adding their multiplicities, and equal even
parts by \eqref{eq:equal-even-combination}.

We also write the two compositions as $\CO_{2\ell,a}*\CO'$ and
$\CO_{(2r+1)^2}*\CO'$, where $\CO_{2\ell,a}$ is an orbit in
$\mathfrak{sp}_{2(m+\ell)}(\Fq)$ and $\CO_{(2r+1)^2}$ is an orbit in
$\mathfrak{sp}_{2(m+2r+1)}(\Fq)$:
\[
 \CO_{2\ell,a}=[(2\ell,a),1^{2m}],
\]
\begin{equation*}
 \CO_{(2r+1)^2}=[(2r+1)^2,1^{2m}].
\end{equation*}
Every rational nilpotent orbit $\CO$ in $\mathfrak{sp}_{2n}(\Fq)$ has a
largest part that is even or occurs in an odd pair. Removing that row or pair
gives
\begin{equation*}
 \CO=\CO_{2\ell,a}*\CO'
 \qquad\text{or}\qquad
 \CO=\CO_{(2r+1)^2}*\CO',
\end{equation*}
where $\CO'$ is an orbit in $\mathfrak{sp}_{2n-2\ell}(\Fq)$ or
$\mathfrak{sp}_{2n-4r-2}(\Fq)$, respectively. For a repeated even part, split its
block using \eqref{eq:equal-even-combination}. In the second case, $r=0$
forces $\CO'=[1^{2n-2}]$ and $\CO=[1^{2n}]$.
\endgroup

\begin{proposition}\label{prop:even-composition-law}
Let $\pi$ be a finite-dimensional representation of $\Sp_{2n}(\Fq)$. If $\ell>0$ and $\CO'$ is
a rational nilpotent orbit in $\mathfrak{sp}_{2n-2\ell}(\Fq)$ whose largest part is at most $2\ell$, then
\begin{equation}\label{eq:even-composition-law}
 \dim\Wh_{\CO_{2\ell,a}*\CO'}(\pi)
 =\dim\Wh_{\CO'}\!\left(\cD^{\FJ}_{2\ell,a}(\pi)\right).
\end{equation}
\end{proposition}
\begin{proof}
This follows from \cite[Lemma~7.6]{PengWang}; see also \cite[Lemma~5.4]{wang}.
\end{proof}

\begin{proposition}\label{prop:odd-pair-identities}
Let $\pi$ be a finite-dimensional representation of $\Sp_{2n}(\Fq)$, let
$0\le\ell\le(n-1)/2$, and let $\CO'$ be a rational nilpotent orbit in
$\mathfrak{sp}_{2n-4\ell-2}(\Fq)$ whose largest part is at most $2\ell+1$.

\textup{(i)} If $\ell>0$, then
\begin{equation*}
 \dim\Wh_{\CO_{(2\ell+1)^2}*\CO'}(\pi)
 =\dim\Wh_{\CO'}\!\left(\Wh_{\CO_{(2\ell+1)^2}}(\pi)\right).
\end{equation*}

\textup{(ii)} For every $a\in\{\pm1\}$,
\begin{align}
 2\dim\Wh_{\CO_{(2\ell+1)^2}*\CO'}(\pi)
 ={}&\dim\Wh_{\CO'}\!\left(
       \cD^{\FJ}_{2\ell+2,\delta a}
       \cD^{\FJ}_{2\ell,a}(\pi)\right)\notag\\
 &+\frac{q+1}{2}\dim\Wh_{\CO'}\!\left(
       \cD^{\FJ}_{2\ell,\delta a}
       \cD^{\FJ}_{2\ell+2,a}(\pi)\right)\notag\\
 &+\frac{q-1}{2}\dim\Wh_{\CO'}\!\left(
       \cD^{\FJ}_{2\ell,-\delta a}
       \cD^{\FJ}_{2\ell+2,-a}(\pi)\right).
 \label{eq:odd-general-input}
\end{align}
\end{proposition}

For $\ell=0$, both $\CO'=[1^{2n-2}]$ and $\CO_{1^2}*\CO'=[1^{2n}]$
are zero orbits. Thus, by \eqref{eq:odd-general-input}, with $c=\delta a$, we have
\begin{equation}\label{eq:zero-index-odd-pair}
\begin{aligned}
 2\dim\pi={}&\dim\!\left(\cD^{\FJ}_{2,c}\cD^{\FJ}_{0,a}(\pi)\right)\\
 &+\frac{q+1}{2}\dim\!\left(\cD^{\FJ}_{0,c}\cD^{\FJ}_{2,a}(\pi)\right)
  +\frac{q-1}{2}\dim\!\left(\cD^{\FJ}_{0,-c}\cD^{\FJ}_{2,-a}(\pi)\right).
\end{aligned}
\end{equation}

\begin{proof}
\emph{Case $\ell=0$.} We prove \eqref{eq:zero-index-odd-pair}.
\begingroup
As in the proof of \cite[Theorem~5.6]{wang}, we compare Fourier expansions
for the two orders of descent. Let $U$ be the unipotent radical of the
stabilizer of an isotropic line in $\Sp_{2n}(\Fq)$, and let
\[
 Z=\{u(0,0,z):z\in\Fq\}
\]
be its centre, in the coordinates of \cite[(4.16)]{wang}. For a
representation $\sigma$ of $U$ and $t\in\Fq$, write
\[
 \sigma[t]:=\{v\in\sigma:\sigma(u(0,0,z))v=\psi(tz)v
 \text{ for every }z\in\Fq\}.
\]
By the Stone--von Neumann theorem, for $t\ne0$ the group $U$ has a
unique irreducible representation with central character $z\mapsto\psi(tz)$,
and its dimension is $q^{n-1}$. Let $\sigma$ be a representation of
$\Sp_{2n}(\Fq)$ and let $t\ne0$.
By the definition of the Fourier--Jacobi descent and conjugation by
$\GL_1(\Fq)$, we have
\begin{equation}\label{eq:heisenberg-isotypic}
 \dim\sigma[t]=q^{n-1}\dim\cD^{\FJ}_{2,\varepsilon_q(t)}(\sigma).
\end{equation}
Taking $\sigma=\pi$ in \eqref{eq:heisenberg-isotypic}
and summing over the nonzero square classes, each of size $(q-1)/2$, we obtain,
since $\pi[0]=\pi^Z$,
\begin{equation}\label{eq:ell-zero-restriction}
 \dim\pi=\dim\pi^Z+\frac{q-1}{2}q^{n-1}
 \left(\dim\cD^{\FJ}_{2,a}(\pi)
       +\dim\cD^{\FJ}_{2,-a}(\pi)\right).
\end{equation}

Choose $c_0\in\Fq^\times$ such that $\varepsilon_q(c_0)=c=\delta a$, and put
$\psi_c(z)=\psi(c_0z)$. By \eqref{eq:FJ-zero},
$\cD^{\FJ}_{0,a}(\pi)=\pi\otimes\omega_{n,\psi_c}$. Choose a normalization of the
Schr\"odinger model compatible with \cite[(4.16)]{wang}; then
$u(0,0,z)$ acts on $\mathbb C[\Fq\times\Fq^{n-1}]$ by
\[
 (\omega_{n,\psi_c}(u(0,0,z))f)(x,y)=\psi(c_0x^2z)f(x,y).
\]
The $c_0$-eigenspace of the centre in $\pi\otimes\omega_{n,\psi_c}$
therefore has dimension
\[
 q^{n-1}\sum_{x\in\Fq}\dim\pi[c_0(1-x^2)].
\]
Applying \eqref{eq:heisenberg-isotypic} to
$\sigma=\pi\otimes\omega_{n,\psi_c}$ and $t=c_0$, we obtain
\[
 \dim\!\left(\cD^{\FJ}_{2,c}\cD^{\FJ}_{0,a}(\pi)\right)
 =\sum_{x\in\Fq}\dim\pi[c_0(1-x^2)].
\]
For $x\ne\pm1$, we have
$\varepsilon_q(c_0(1-x^2))=a\,\varepsilon_q(x^2-1)$. Extending $\varepsilon_q$ to $\Fq$ by $\varepsilon_q(0)=0$, we have
$\#\{(x,y)\in\Fq^2:x^2-y^2=1\}=q-1$, so
$\sum_{x\in\Fq}\bigl(1+\varepsilon_q(x^2-1)\bigr)=q-1$, that is,
$\sum_{x\in\Fq}\varepsilon_q(x^2-1)=-1$. Hence
\[
 \begin{aligned}
 \#\{x:x^2=1\}&=2,\\
 \#\{x:\varepsilon_q(x^2-1)=1\}&=\frac{q-3}{2},\qquad
 \#\{x:\varepsilon_q(x^2-1)=-1\}=\frac{q-1}{2}.
 \end{aligned}
\]
It follows that
\begin{equation}\label{eq:ell-zero-schrodinger}
\begin{aligned}
 \dim\!\left(\cD^{\FJ}_{2,c}\cD^{\FJ}_{0,a}(\pi)\right)
 ={}&2\dim\pi^Z\\
 &+q^{n-1}\left(
 \frac{q-3}{2}\dim\cD^{\FJ}_{2,a}(\pi)
 +\frac{q-1}{2}\dim\cD^{\FJ}_{2,-a}(\pi)\right).
\end{aligned}
\end{equation}
By \eqref{eq:FJ-zero}, for every representation $\sigma$ of
$\Sp_{2n-2}(\Fq)$, we have
\[
 \dim\cD^{\FJ}_{0,\pm c}(\sigma)=q^{n-1}\dim\sigma.
\]
Write $A_\pm=\dim\cD^{\FJ}_{2,\pm a}(\pi)$. By
\eqref{eq:ell-zero-schrodinger} and the last identity, the right-hand side
of \eqref{eq:zero-index-odd-pair} equals
\[
 \begin{aligned}
 &2\dim\pi^Z+q^{n-1}\Bigl(\bigl(\tfrac{q-3}{2}+\tfrac{q+1}{2}\bigr)A_+
   +\bigl(\tfrac{q-1}{2}+\tfrac{q-1}{2}\bigr)A_-\Bigr)\\
 &\qquad=2\Bigl(\dim\pi^Z+\tfrac{q-1}{2}q^{n-1}(A_++A_-)\Bigr),
 \end{aligned}
\]
which is $2\dim\pi$ by \eqref{eq:ell-zero-restriction}. The
unknown dimension $\dim\pi^Z$ cancels.
\endgroup

\emph{Case $\ell>0$.}
\begingroup
Part \textup{(i)} is
\cite[Theorem~5.8]{wang}. For part \textup{(ii)}, set $2m=2n-4\ell-2$ and
\[
 \CO[u,v]:=[(2\ell+2,u),(2\ell,v),1^{2m}]
 \qquad(u,v\in\{\pm1\}).
\]
We follow the proof of \cite[Theorem~5.1]{wang}.
The root exchanges and Fourier expansions in the proof of
\cite[Theorem~5.6]{wang} are compatible with the action of
$\Sp_{2m}(\Fq)$. That calculation is therefore an identity of
virtual representations of $\Sp_{2m}(\Fq)$:
\[
 \begin{aligned}
 2\bigl[\Wh_{\CO_{(2\ell+1)^2}}(\pi)\bigr]
 ={}&\bigl[\cD^{\FJ}_{2\ell+2,\delta a}\cD^{\FJ}_{2\ell,a}(\pi)\bigr]
 +\frac{q+1}{2}\bigl[\Wh_{\CO[a,\delta a]}(\pi)\bigr]\\
 &+\frac{q-1}{2}\bigl[\Wh_{\CO[-a,-\delta a]}(\pi)\bigr].
 \end{aligned}
\]
Since $\Wh_{\CO'}$ is an exact functor on finite-dimensional
representations, taking dimensions of $\CO'$-Whittaker models gives
\[
\begin{aligned}
 2\dim\Wh_{\CO'}\!\left(\Wh_{\CO_{(2\ell+1)^2}}(\pi)\right)
 ={}&\dim\Wh_{\CO'}\!\left(
       \cD^{\FJ}_{2\ell+2,\delta a}\cD^{\FJ}_{2\ell,a}(\pi)\right)\\
 &+\frac{q+1}{2}\dim\Wh_{\CO'}
       \!\left(\Wh_{\CO[a,\delta a]}(\pi)\right)\\
 &+\frac{q-1}{2}\dim\Wh_{\CO'}
       \!\left(\Wh_{\CO[-a,-\delta a]}(\pi)\right).
\end{aligned}
\]
By \textup{(i)}, the left side is
$2\dim\Wh_{\CO_{(2\ell+1)^2}*\CO'}(\pi)$. To obtain
\eqref{eq:odd-general-input}, it remains to express the last
two terms by successive Fourier--Jacobi descents, that is, to
prove \eqref{eq:two-even-block-dimension} below for $(u,v)=(a,\delta a)$ and
$(-a,-\delta a)$.
\endgroup

\emph{The last two terms.}
Let $V$ be the $2n$-dimensional symplectic space over $\Fq$ on which $\Sp_{2n}(\Fq)$
acts, and fix an $\mathfrak{sl}_2$-triple attached to a representative of
$\CO[u,v]=[(2\ell+2,u),(2\ell,v),1^{2m}]$. As an $\mathfrak{sl}_2$-module,
$V$ decomposes orthogonally as $V=W_{2\ell+2}\oplus W_{2\ell}\oplus V_0$, where $W_{2\ell+2}$ and $W_{2\ell}$
correspond to the two even Jordan blocks and $V_0$ is the $2m$-dimensional trivial
$\mathfrak{sl}_2$-submodule. Viewing $\CO'$ in $\mathfrak{sp}(V_0)$, put
$U'=G_{\ge1.5}(\CO')\subset\Sp(V_0)$ and $\psi'=\psi_{\CO'}$.

To prove \eqref{eq:two-even-block-dimension}, we apply the root-exchange argument of
\cite[proof of Lemma~7.6]{PengWang} (see also \cite[Lemma~5.4]{wang}) and check that it
is $U'$-equivariant. That argument first replaces each Weil representation $\omega$ by a
$1$-dimensional character via a Lagrangian polarization so that root exchanges apply.
Although no Lagrangian subspace of $V_0$ is $\Sp(V_0)$-stable, the unipotent subgroup
$U'\subset\Sp(V_0)$ stabilizes a Lagrangian subspace $L_0\subset V_0$ and acts on $L_0$
with determinant $1$. Via $\mathfrak{sp}(V)\simeq\operatorname{Sym}^2 V$, the weight-one
space decomposes orthogonally as
\[
 \Fg_1(\CO[u,v])^F\simeq(W_{2\ell+2})_1\otimes V_0\oplus(W_{2\ell})_1\otimes V_0
 \simeq V_0\oplus V_0,
\]
so $L_0$ defines a $U'$-stable Lagrangian subspace $L_0\oplus L_0$ of
$\Fg_1(\CO[u,v])^F$. On each summand $V_0$ of $\Fg_1(\CO[u,v])^F\simeq V_0\oplus V_0$, writing $K$ for the Heisenberg group,
$K_{L_0}\subset K$ for the preimage of $L_0$, $\chi$ for the extension of the central
character that is trivial on $L_0$, and $\omega\simeq\operatorname{Ind}_{K_{L_0}}^K\chi$,
Frobenius reciprocity gives a $U'$-equivariant isomorphism
\[
 \Hom_K(\omega,T)\xrightarrow{\ \sim\ }
 \{w\in T:kw=\chi(k)w\text{ for all }k\in K_{L_0}\}
\]
for any representation $T$ of $U'\ltimes K$. Moreover, the root subgroups exchanged at each step of
\cite[proof of Lemma~7.6]{PengWang} pair only vectors of nonzero $\mathfrak{sl}_2$-weight
in $W_{2\ell+2}$ and $W_{2\ell}$; hence they lie in $\Sp(W_{2\ell+2}\oplus W_{2\ell})$
and commute with $U'\subset\Sp(V_0)$. Since $L_0\oplus L_0$ is $U'$-stable, every step of
the root exchange is $U'$-equivariant, yielding
$\Wh_{\CO[u,v]}(\pi)\simeq\cD^{\FJ}_{2\ell,v}\cD^{\FJ}_{2\ell+2,u}(\pi)$ as
$U'$-modules. Taking $\psi'$-eigenspaces gives
\begin{equation}\label{eq:two-even-block-dimension}
 \dim\Wh_{\CO'}\!\left(\Wh_{\CO[u,v]}(\pi)\right)
 =\dim\Wh_{\CO'}\!\left(
       \cD^{\FJ}_{2\ell,v}\cD^{\FJ}_{2\ell+2,u}(\pi)\right),
\end{equation}
and applying \eqref{eq:two-even-block-dimension} to $(u,v)=(a,\delta a)$ and
$(-a,-\delta a)$ yields \eqref{eq:odd-general-input}.
\end{proof}

\subsection{Iterated Fourier--Jacobi descents}\label{subsec:sequence-coefficients}
\begingroup
Let $\Gamma=\llbracket(2\ell_1,a_1),\ldots,(2\ell_r,a_r)\rrbracket$ be a descent sequence index of rank $N$; thus $\sum_j\ell_j=N$. For a representation $\pi$ of $\Sp_{2N}(\Fq)$, define
\begin{equation}\label{eq:cGamma-iterated}
 m_\Gamma(\pi):=\dim\!\left(
 \cD^{\FJ}_{2\ell_r,a_r}\circ\cdots\circ
 \cD^{\FJ}_{2\ell_1,a_1}(\pi)\right).
\end{equation}
As $q$ varies, let
$\Gamma=\llbracket(2\ell_1,a_1),\ldots,(2\ell_r,a_r)\rrbracket$ and
$\Gamma'=\llbracket(2\ell'_1,a'_1),\ldots,(2\ell'_r,a'_r)\rrbracket$ be descent
sequence indices for $\Sp_{2N}(\Fq)$ and $\Sp_{2N}(\mathbb F_{q'})$, respectively. We
identify $\Gamma$ with $\Gamma'$ if
\[
 \ell_j=\ell'_j
 \qquad\text{and}\qquad
 \varepsilon_q(-1)^{j-1}a_j=\varepsilon_{q'}(-1)^{j-1}a'_j
 \qquad(1\le j\le r);
\]
by abuse of notation, we still write $\Gamma$ for the identified index. This
identification preserves \eqref{eq:good-cases}, \eqref{eq:elementary-seq-mutation}, and
$\Gamma\mapsto\CO(\Gamma)$ (compatibly with Section~\ref{subsec:rational-orbit-order}).
For a $q$-system $\boldsymbol\pi=\{\pi_q\}$ of type
$\mathfrak t$, put
\begin{equation*}
 m_{\Gamma,\boldsymbol\pi}(q):=m_\Gamma(\pi_q).
\end{equation*}

Decomposing each descent into irreducible representations gives
\begin{equation}\label{eq:cGamma}
 m_{\Gamma,\boldsymbol\pi}(q)
 =\sum_{\mathcal P}\prod_{j=1}^r
   m_{2\ell_j,a_j}\bigl(\pi^{(j-1)},\pi^{(j)}\bigr),
\end{equation}
where the sum runs over descent paths with index $\Gamma$, allowing appended zero
entries as in Section~\ref{subsec:good-indexes}:
\[
 \begin{gathered}
 \mathcal P:\quad
 \pi_q=\pi^{(0)}\searrow_{2\ell_1,a_1}\pi^{(1)}
 \searrow\cdots\searrow_{2\ell_r,a_r}\pi^{(r)}=\mathbf{1}_{\Sp_0},\\
 \pi^{(j)}\in\Irr(\Sp_{2N_j}(\Fq)),\qquad
 N_j=N-\sum_{h\le j}\ell_h.
 \end{gathered}
\]
In particular, $m_{\Gamma,\boldsymbol\pi}(q)>0$ if and only if $\Gamma$ is
realized by $\pi_q$. Sorting the paths in \eqref{eq:cGamma} by their
combinatorial data will give the polynomial $m_{\Gamma,\mathfrak t}(X)$
(Section~\ref{sec:poly}), and the composition laws express the Whittaker
multiplicities of rational orbits through the $m_\Gamma$.
\endgroup

\medskip
From now until Section~\ref{sec:orth-theta}, we take $G=\Sp_{2N}(\Fq)$.

\section{First descents and rank inequalities}\label{sec:minimality-degree}
This section contains the rank inequalities used in
Section~\ref{sec:poly} to bound the degree of $m_{\Gamma,\mathfrak t}(X)$
(see \eqref{eq:weighted-descent-paths}). The main result is
Proposition~\ref{prop:stagewise}.

\subsection{Ranks after first descents}\label{subsec:primary-minimality}
Fix $\pi=\pi_{\rho,\Lambda_1,\Lambda_{-1}}$ and choose a
first-descent path
\[
 \mathcal P_{\rm fd}:\quad
 \pi=\pi^{(0)}\fdsearrow\pi^{(1)}\fdsearrow\cdots
 \fdsearrow\pi^{(s)}=\mathbf{1}_{\Sp_0}.
\]
Let
$\Lambda_1^{(j)},\Lambda_{-1}^{(j)}$ be the two Lusztig symbols of
$\pi^{(j)}$ \textup{(}so $\Lambda_\xi^{(0)}=\Lambda_\xi$\textup{)}.
For $\xi\in\{\pm1\}$ and $1\le j\le s$ we have, by
Proposition~\ref{cor:first-descent-rule},
\[
 {}^d\Lambda_{-\xi}^{(j)}
 =\bigl(\partial_\eta({}^d\Lambda_\xi^{(j-1)})\bigr)^{\langle\epsilon\rangle},
\]
where $\partial_\eta$ removes one copy of the largest entry and
$\eta,\epsilon$ are the signs in the proposition. The key point is that
the rank of the right-hand side does not depend on $\eta$ or $\epsilon$:
Alvis--Curtis duality and row interchange preserve the rank, and if the largest
entry occurs in both rows, deleting either copy leaves the same multiset of
entries. Hence, after $j$ steps, $\Lambda_\xi$ has become
$\Lambda_{(-1)^j\xi}^{(j)}$, and its rank depends only on $\Lambda_\xi$ and $j$,
not on the path. We can therefore define
\begin{equation*}
 r_j^{\rm fd}(\Lambda_\xi)
 :=\rank\!\left(\Lambda_{(-1)^j\xi}^{(j)}\right)
 \qquad(0\le j\le s,\ \xi=\pm1),
\end{equation*}
and we set $r_j^{\rm fd}(\Lambda_\xi)=0$ for $j>s$. These ranks do not depend on
the choice of symbol representatives; see \eqref{eq:symbol-first-rank-entries}.
Note that $r_0^{\rm fd}(\Lambda)=\rank\Lambda$.

At a first descent, the partition $\lambda[a]$ at $[a]\ne[\pm1]$
becomes $\lambda[a]_-$ at $[-a]$. All nonempty partitions at the next step arise
in this way. Let $\mathfrak t$ be the type of $\pi$. For $j\ge0$ we define
\begin{equation}\label{eq:first-descent-inherited-rank}
\begin{aligned}
 r_j^{\rm fd}(\mathfrak t)
 \mathrel{:=}{}&\sum_{[a]\ne[\pm1]}d([a])
       \sum_{i\ge1}\max\{\lambda[a]_i-j,0\}\\
 &+r_j^{\rm fd}(\Lambda_1)+r_j^{\rm fd}(\Lambda_{-1}).
\end{aligned}
\end{equation}
For $0\le j\le s$, the representation $\pi^{(j)}$ belongs to
$\Irr(\Sp_{2r_j^{\rm fd}(\mathfrak t)}(\Fq))$. Thus
$r_j^{\rm fd}(\mathfrak t)$ is the rank of the group after $j$ successive
first descents. This rank depends only on $\mathfrak t$, even though the
constituents $\pi^{(j)}$ need not be unique.
The index $j$ counts descent steps, including steps with $\ell=0$.

\subsection{Comparison with arbitrary descents}
Now consider an arbitrary descent path
\[
 \mathcal P:\quad
 \pi=\pi^{(0)}\searrow_{2\ell_1,a_1}\pi^{(1)}
 \searrow\cdots\searrow_{2\ell_r,a_r}\pi^{(r)}=\mathbf{1}_{\Sp_0}
\]
with $\pi^{(j)}\in\Irr(\Sp_{2N_j}(\Fq))$ and initial
type $\mathfrak t$. Write $\lambda[b]^{(j)}$,
$\Lambda_1^{(j)},\Lambda_{-1}^{(j)}$ for its Lusztig parameters at
stage $j$. If $[b]$ does not occur in the semisimple parameter
of $\pi^{(j)}$, set $\lambda[b]^{(j)}=\varnothing$.
Write $\Lambda^{(j)}_{\pm1}=(\Lambda^{(j)}_1,\Lambda^{(j)}_{-1})$.
At the $j$th step, by \eqref{eq:GGP-primary-numerical},
$\lambda[a]^{(j-1)}$ is paired with $\lambda[-a]^{(j)}$, and
$\Lambda^{(j-1)}_\xi$ is paired with $\Lambda^{(j)}_{-\xi}$ for $\xi=\pm1$
in the factor $\mathbf b_t(\Lambda^{(j-1)}_{\pm1};\Lambda^{(j)}_{\pm1})$
of \eqref{eq:symbol-pairing-factor}. We say that the following parameters of
$\pi^{(j)}$ are \emph{inherited from $\pi$}:
\begin{equation}\label{eq:inherited-primary-partition}
 \lambda[(-1)^ja]^{(j)}\quad([a]\ne[\pm1]\text{ occurring in }\pi),
 \qquad \Lambda_{(-1)^j\xi}^{(j)}\quad(\xi=\pm1).
\end{equation}
We call
\begin{equation}\label{eq:inherited-rank}
 r_j^{\rm inh}
 :=\sum_{\substack{[a]\ne[\pm1]\\ [a]\text{ occurs in }\pi}}
 d([a])\,\bigl|\lambda[(-1)^ja]^{(j)}\bigr|
 +\rank\Lambda_1^{(j)}+\rank\Lambda_{-1}^{(j)}
\end{equation}
the \emph{inherited rank} after $j$ steps; it is the part of $N_j$
occupied by parameters inherited from $\pi$.

\begin{proposition}
\label{prop:stagewise}
For every $0\le j\le r$, we have
\begin{equation}\label{eq:stagewise-minimality}
 r_j^{\rm inh}\ge r_j^{\rm fd}(\mathfrak t).
\end{equation}
Equality holds if the first $j$ steps are successive first
descents.
\end{proposition}
\begin{proof}
By \eqref{eq:GGP-primary-numerical} and
\eqref{eq:exact-primary-pairings}, for the partitions $\lambda[a]$ we have
\[
 \lambda[(-1)^{h+1}a]^{(h+1)}_i
 \ge\max\{\lambda[(-1)^ha]^{(h)}_i-1,0\}.
\]
Indeed, a nonzero factor $b_\kappa(\lambda,\lambda')$ in
\eqref{eq:GGP-primary-numerical} forces $\lambda$ and $\lambda'$ to be close.
Then
\begin{equation}\label{eq:primary-stagewise}
 \lambda[(-1)^ja]^{(j)}_i
 \ge\max\{\lambda[a]^{(0)}_i-j,0\}
 \qquad(i\ge1).
\end{equation}
At a first descent, the partition on the left is obtained by deleting the
first column, so equality holds at every step.

For the symbols, apply Lemma~\ref{lem:symbol-first-descent-comparison} in
\ref{app:symbol} to obtain the one-step comparison
\begin{equation*}
 r_k^{\rm fd}\!\left(\Lambda_{(-1)^{h+1}\xi}^{(h+1)}\right)
 \ge r_{k+1}^{\rm fd}\!\left(\Lambda_{(-1)^h\xi}^{(h)}\right)
 \qquad(k\ge0).
\end{equation*}
By induction on $j$, we obtain
\begin{equation}\label{eq:symbol-stagewise}
 r_k^{\rm fd}\!\left(\Lambda_{(-1)^j\xi}^{(j)}\right)
 \ge r_{k+j}^{\rm fd}\!\left(\Lambda_\xi^{(0)}\right)
 \qquad(k\ge0,\ \xi=\pm1).
\end{equation}
Equality holds at every first descent by the same lemma.

Finally, multiply \eqref{eq:primary-stagewise} by $d([a])$, sum over the rows
and the initial classes $[a]$, and add
\eqref{eq:symbol-stagewise} for $k=0$ and $\xi=\pm1$.
By \eqref{eq:first-descent-inherited-rank} and
\eqref{eq:inherited-rank}, we obtain \eqref{eq:stagewise-minimality}, with
equality along successive first descents.
\end{proof}

The remaining classes $[a]$ at stage $j$ account for rank
\begin{equation}\label{eq:rank-split-inherited-new}
 \sum_{\substack{[b]\ne[\pm1]\\
        [(-1)^jb]\text{ absent from }\pi}}
 d([b])\,|\lambda[b]^{(j)}|
 =N_j-r_j^{\rm inh}
 \le N_j-r_j^{\rm fd}(\mathfrak t).
\end{equation}
In Section~\ref{subsec:budget-sign-average} we apply
\eqref{eq:rank-split-inherited-new} at successive steps to bound the degree
of the counting polynomials.

\section{\texorpdfstring{Calculation of $m_\Gamma(\pi)$ and its degree in $q$}{Calculation of m\_Gamma(pi) and its degree in q}}\label{sec:poly}

\subsection{Grouping descent paths}\label{subsec:type-transitions}
\par
Fix a type $\mathfrak t=(\{(d_i,\kappa_i,\lambda_i)\}_{i\in I};\Lambda_1,\Lambda_{-1})$
of rank $N$, as in \eqref{eq:qtype}, and a descent sequence index
\[
 \Gamma=\llbracket(2\ell_1,a_1),\ldots,(2\ell_r,a_r)\rrbracket,
 \qquad \sum_{j=1}^r\ell_j=N.
\]
We prove that, for fixed $q$, the multiplicity $m_\Gamma(\pi)$ takes the same value
for all $\pi\in\Irr_{\mathfrak t}(G)$, and that these multiplicities are polynomial in
$q$ as $q$ varies.

Let $\pi=\pi^s_{\rho,\Lambda_1,\Lambda_{-1}}\in\Irr_{\mathfrak t}(G)$ where
$
 \rho=\bigotimes_{i\in I}\pi[b_i],$ and
$[b_i]$ are the distinct classes $\ne[\pm1]$ occurring in the
semisimple parameter $s$, indexed so that
\[
 d([b_i])=d_i,\qquad \kappa([b_i])=\kappa_i,\qquad \lambda[b_i]=\lambda_i
 \qquad(i\in I).
\]
Thus $\mathfrak t$ records the degrees, kinds and partitions at the classes $[b_i]$,
together with the symbols $\Lambda_1,\Lambda_{-1}$, but not the classes $[b_i]$
themselves.
\par

\par
For $\pi$ as above, put $N_1=N-\ell_1$ and
$\Gamma_{>1}:=\llbracket(2\ell_2,a_2),\ldots,(2\ell_r,a_r)\rrbracket$. Writing
$\cD^{\FJ}_{2\ell_1,a_1}(\pi)$ as a sum of irreducible representations, we have, by
\eqref{eq:cGamma-iterated} and \eqref{eq:GGP-mult},
\begin{equation}\label{eq:cGamma-first-step-type}
\begin{aligned}
 m_\Gamma(\pi)
 &=\sum_{\substack{\pi'\in\Irr(\Sp_{2N_1}(\Fq))\\ \pi\searrow_{2\ell_1,a_1}\pi'}}
   m_{2\ell_1,a_1}(\pi,\pi')\,m_{\Gamma_{>1}}(\pi')\\
 &=\sum_{\mathfrak t'}\ \sum_{\substack{\pi'\in\Irr_{\mathfrak t'}(\Sp_{2N_1}(\Fq))\\ \pi\searrow_{2\ell_1,a_1}\pi'}}
   m_{2\ell_1,a_1}(\pi,\pi')\,m_{\Gamma_{>1}}(\pi'),
\end{aligned}
\end{equation}
where $\mathfrak t'$ runs over the types of rank $N_1$.
We group the representations $\pi'$ with $\pi\searrow_{2\ell_1,a_1}\pi'$, that is, the
irreducible constituents of the descent $\cD^{\FJ}_{2\ell_1,a_1}(\pi)$, according to the
data inherited from $\pi$ along this descent, namely their Lusztig parameters at the classes matched
with those of $\pi$ and their symbols.
\par

Let $\lambda'[-b_i]$ be the partition attached to $[-b_i]$ in the
Lusztig parameter of $\pi'$, with $\lambda'[-b_i]=\varnothing$ if $[-b_i]$
does not occur. We call

\[
 \mathfrak d=(\{\lambda_i^{\mathfrak d}\}_{i\in I};\Lambda_1^{\mathfrak d},\Lambda_{-1}^{\mathfrak d}),
 \qquad \lambda_i^{\mathfrak d}=\lambda'[-b_i],\quad
 \Lambda_{\pm1}^{\mathfrak d}=\Lambda'_{\pm1},
\]
the inherited data of $\pi'$ relative to $\pi$. 
Only the $\pi'$ with $m_{2\ell_1,a_1}(\pi,\pi')\ne0$ contribute to
\eqref{eq:cGamma-first-step-type}. By \eqref{eq:GGP-primary-numerical}, we may therefore require that
$b_{\kappa_i}(\lambda_i,\lambda_i^{\mathfrak d})\ne0$ for every $i$
and that
$\mathbf b_t(\Lambda_{\pm1};\Lambda^{\mathfrak d}_{\pm1})\ne0$
for at least one value of $t$; the sign $t$ in \eqref{eq:GGP-rational-sign} also depends on the
new classes of $\pi'$ introduced below, which are not recorded in $\mathfrak d$.
Define
\begin{equation}\label{eq:H-inherited-rank}
 r(\mathfrak d)=\sum_i d_i|\lambda_i^{\mathfrak d}|+\rank\Lambda_1^{\mathfrak d}+\rank\Lambda_{-1}^{\mathfrak d},
 \qquad m(\mathfrak d)=N_1-r(\mathfrak d).
\end{equation}
There are only finitely many such data $\mathfrak d$ with $r(\mathfrak d)\le N_1$.

By \eqref{eq:GGP-one-sided-column}, if $m_{2\ell_1,a_1}(\pi,\pi')\ne0$, then, for $[x]\ne[\pm1]$, $\lambda'[x]=(1^e)$ for some $e\ge1$ whenever
$[x]\notin\{[-b_i]:i\in I\}$ and $\lambda'[x]\ne\varnothing$.
We call these classes $[x]$ the \emph{new classes} of $\pi'$.
For $d,e\ge1$ and $\kappa\in\{\GL,\U\}$, let $n_{d,\kappa,e}$ count the new
classes of degree $d$ and kind $\kappa$ carrying $(1^e)$.
Let $B$ be the set of new classes. By Remark~\ref{rem:generic-new-classes}, the corresponding $F$-stable Levi subgroup $\mathbf L$ of $\bSp_{2N_1}$ has

\[
 \mathbf L^F\simeq\left(\prod_{[x]\in B}H_x\right)\times\Sp_{2r(\mathfrak d)}(\Fq),
\]

where, for a new class $[x]$ of degree $d$ carrying $(1^e)$, the factor $H_x$ is $\GL_e(\mathbb F_{q^d})$ or $\U_e(\mathbb F_{q^d})$ according to its kind. The symplectic factor carries the representation $\sigma$ of Remark~\ref{rem:generic-new-classes}.
The tuple $\mathbf{n}=(n_{d,\kappa,e})$ satisfies
\begin{equation}\label{eq:target-profile-rank-split}
 \sum_{d,\kappa,e}de\,n_{d,\kappa,e}=m(\mathfrak d).
\end{equation}
We call $(\mathfrak d,\mathbf n)$ the \emph{step data} of $\pi'$ relative to
$\pi$; that is, $\mathfrak d$ records the data of $\pi'$ inherited from $\pi$, and
$\mathbf n$ records its new classes.
The type $\mathfrak t(\mathfrak d,\mathbf n)$ of $\pi'$ consists of
the nonempty triples $(d_i,\kappa_i,\lambda_i^{\mathfrak d})$,
$n_{d,\kappa,e}$ copies of $(d,\kappa,(1^e))$, and the symbols
$\Lambda_1^{\mathfrak d},\Lambda_{-1}^{\mathfrak d}$.
Note that the step data $(\mathfrak d,\mathbf n)$ of $\pi'$ carry more information than
the type of $\pi'$: besides determining $\mathfrak t(\mathfrak d,\mathbf n)$, they fix the
parameters of $\pi'$ \textup{(}see \eqref{eq:qtype-primary}\textup{)} at the classes $[-b_i]$.

When $\pi$ varies in $\Irr_{\mathfrak t}(G)$, its semisimple
parameter $s$ varies, and so do the classes $[b_i]$, $i\in I$. The classes of $\pi'$ attached to the indices $i\in I$ change accordingly, being
always the classes $[-b_i]$. This change does not affect the multiplicity:
$m_{2\ell_1,a_1}(\pi,\pi')$ in \eqref{eq:GGP-primary-numerical} depends only on $\mathfrak t$,
on $(2\ell_1,a_1)$ and on the step data
$(\mathfrak d,\mathbf n)$ of $\pi'$ relative to $\pi$, as the computation below shows; we denote it by
$m_{2\ell_1,a_1}(\mathfrak t;\mathfrak d,\mathbf n)$.

To compute it, let $\rho_{\mathfrak d}=\bigotimes_{i\in I}\pi'[-b_i]$ be the part of $\rho'$
at the classes $[-b_i]$. By \eqref{eq:epsilon-rho}, $\epsilon_{\rho'}=\epsilon_{\rho_{\mathfrak d}}\eta(\mathbf n)$, where
\[
 \epsilon_{\rho_{\mathfrak d}}=(-1)^{\sum_{i:\kappa_i=\U}|\lambda_i^{\mathfrak d}|},
 \qquad \eta(\mathbf{n})=(-1)^{\sum_{d,e}e n_{d,\U,e}}.
\]
Let $\delta_{-1},\delta'_{-1}$ be the discriminants of the even orthogonal factors
labelled by $\Lambda_{-1}$ and $\Lambda_{-1}^{\mathfrak d}$, respectively,
and let
\[
 t=-a_1\epsilon_\rho\epsilon_{\rho_{\mathfrak d}}
       \eta(\mathbf n)\delta_{-1}\delta'_{-1}.
\]
For every $\pi'$ with the prescribed step data,
\begin{equation}\label{eq:MH-definition}
\begin{aligned}
 &m_{2\ell_1,a_1}(\mathfrak t;\mathfrak d,\mathbf n)
   =m_{2\ell_1,a_1}(\pi,\pi')\\
 &\qquad=\prod_{i\in I}
 b_{\kappa_i}(\lambda_i,\lambda_i^{\mathfrak d})\\
 &\qquad\quad\times
 b_{\delta'_{-1}}\!\left(\Lambda_1,(\Lambda_{-1}^{\mathfrak d})^{\langle\delta t\rangle}\right)
 b_{\delta_{-1}}\!\left(\Lambda_1^{\mathfrak d},\Lambda_{-1}^{\langle t\rangle}\right).
\end{aligned}
\end{equation}
The factor
$\mathbf b_t(\Lambda_{\pm1};\Lambda^{\mathfrak d}_{\pm1})$ may be nonzero for
both values of $t$, so
$\mathfrak d$ alone need not determine $\eta(\mathbf n)$.
By Theorem~\ref{thm:FJ-branching}, we may also use
\eqref{eq:MH-definition} when $\ell_1=0$.

We now count the representations $\pi'$ with prescribed step data.
Put
\[
 P_d^\kappa(q):=\#\{[a]\ne[\pm1]:d([a])=d,\ \kappa([a])=\kappa\},
 \qquad d\ge1,\quad\kappa\in\{\GL,\U\}.
\]
Note that $d([a])=\#[a]/2$, and $\kappa([a])=\U$ precisely when the
Frobenius orbit of $a$ contains $a^{-1}$.

\begin{proposition}\label{prop:primary-counting}
The functions $P_d^\kappa(q)$ are polynomials in $q$ and satisfy
\begin{equation}\label{eq:AB-leading}
 P_d^{\GL}(q)=\frac{q^d}{2d}+O(q^{d-1}),\qquad
 P_d^{\U}(q)=\frac{q^d}{2d}+O(q^{d-1}).
\end{equation}
\end{proposition}
\begin{proof}
Let $I_d(q)$ be the number of Frobenius orbits in $\overline{\Fq}^\times$ of
length exactly $d$, that is, the number of monic irreducible polynomials of
degree $d$ in $\Fq[T]$ other than $T$. Since
$\sum_{e\mid d}eI_e(q)=|\mathbb F_{q^d}^\times|=q^d-1$, M\"obius inversion
gives
\begin{equation*}
I_d(q)=\frac{1}d\sum_{r\mid d}\mu_{\rm mob}(r)(q^{d/r}-1),
\end{equation*}
where $\mu_{\rm mob}$ is the M\"obius function.
If a Frobenius orbit of length $2d$ is stable under inversion, then inversion acts on it as $F^d$ (it commutes with $F$ and has no fixed points), that is, $a^{-1}=a^{q^d}$ for every $a$ in the orbit. Hence the orbit lies in $\U_1(\mathbb F_{q^d})$, which has order
$q^d+1$. Removing $\pm1$ and applying M\"obius inversion over the odd divisors of $d$, we get
\begin{equation*}
P_d^{\U}(q)=\frac{1}{2d}
 \sum_{\substack{r\mid d\\ r\text{ odd}}}
 \mu_{\rm mob}(r)(q^{d/r}-1).
\end{equation*}
If $\kappa([a])=\GL$, the Frobenius orbit $C$ of $a$ does not contain $a^{-1}$, so
$[a]=C\sqcup C^{-1}$ is the union of the Frobenius orbits of $a$ and $a^{-1}$. Then
\begin{equation*}
P_d^{\GL}(q)=
\begin{cases}
 \dfrac{q-3}{2},&d=1,\\[
4pt]
 \dfrac12 I_d(q),&d>1\text{ odd},\\[4pt]
 \dfrac12\bigl(I_d(q)-P_{d/2}^{\U}(q)\bigr),&d\text{ even}.
\end{cases}
\end{equation*}
It is now clear that $P_d^{\GL}(q)$ and $P_d^{\U}(q)$ both have leading term $q^d/(2d)$.
\end{proof}
Set
\begin{equation*}
 h_{d,\kappa}:=\#\{i\in I:(d_i,\kappa_i)=(d,\kappa)\},
 \qquad K_{d,\kappa}(\mathbf{n}):=\sum_e n_{d,\kappa,e}.
\end{equation*}
Fix step data $(\mathfrak d,\mathbf n)$. By the discussion above, $\mathfrak d$ fixes the
parameters of $\pi'$ at the classes $[-b_i]$ and its symbols, and each new class carries some
$(1^e)$. Hence $\pi'$ is determined by assigning to each of the $n_{d,\kappa,e}$ new classes
carrying $(1^e)$ a class of degree $d$ and kind $\kappa$ outside $\{[-b_i]:i\in I\}$, distinct
new classes receiving distinct classes. There are $P_d^\kappa(q)-h_{d,\kappa}$ classes of degree
$d$ and kind $\kappa$ outside $\{[-b_i]:i\in I\}$, so the number of such $\pi'$ is
$N_{\mathfrak d,\mathbf n}(q)$, where
\begin{equation}\label{eq:profile-count-product}
 N_{\mathfrak d,\mathbf{n}}(X)=\prod_{d,\kappa}
 \frac{(P_d^\kappa(X)-h_{d,\kappa})_{\underline{K_{d,\kappa}(\mathbf n)}}}
      {\prod_e n_{d,\kappa,e}!}.
\end{equation}
Here $(x)_{\underline K}=x(x-1)\cdots(x-K+1)$ denotes the falling factorial
and $(x)_{\underline0}=1$. We divide by $\prod_e n_{d,\kappa,e}!$ because new classes with the same $e$ are not ordered, whereas the partitions $\lambda_i^{\mathfrak d}$ keep their indices $i\in I$. By Proposition~\ref{prop:primary-counting}, and since $d$, $\kappa$ and $e$ do not depend on $q$, $N_{\mathfrak d,\mathbf n}(X)$ is a polynomial in $X$.
Each of these representations $\pi'$ occurs in $\cD^{\FJ}_{2\ell_1,a_1}(\pi)$ with multiplicity $m_{2\ell_1,a_1}(\mathfrak t;\mathfrak d,\mathbf n)$.
Each $\pi'$ occurs in $\cD^{\FJ}_{2\ell_1,a_1}(\pi)$ precisely
when $m_{2\ell_1,a_1}(\mathfrak t;\mathfrak d,\mathbf n)\ne0$. Several pairs $(\mathfrak d,\mathbf n)$ may give the same type $\mathfrak t'$. Summing their contributions, we obtain
\begin{equation}\label{eq:branching-between-types}
\begin{aligned}
 &\sum_{\pi'\in\Irr_{\mathfrak t'}(\Sp_{2N_1}(\Fq))}
       m_{2\ell_1,a_1}(\pi,\pi')\\
 &\qquad=\sum_{\mathfrak t(\mathfrak d,\mathbf n)=\mathfrak t'}
       {m_{2\ell_1,a_1}(\mathfrak t;\mathfrak d,\mathbf n)}\,
       N_{\mathfrak d,\mathbf n}(q).
\end{aligned}
\end{equation}
Here the first factor is the multiplicity of each such
$\pi'$, and the second is the number of such $\pi'$.
\begin{proposition}\label{prop:exact-poly}
There is a polynomial
$m_{\Gamma,\mathfrak t}(X)\in\mathbb Q[X]$ such that, for every
sufficiently large $q$ and every
$\pi\in\Irr_{\mathfrak t}(\Sp_{2N}(\Fq))$,
\begin{equation}\label{eq:type-polynomial-value}
 m_\Gamma(\pi)=m_{\Gamma,\mathfrak t}(q).
\end{equation}
In particular, this polynomial is independent of the realization of $\mathfrak t$ and of $\varepsilon_q(-1)$.
\end{proposition}
\begin{proof}
Fix $\delta\in\{\pm1\}$ and consider first $q$ with $\varepsilon_q(-1)=\delta$. We argue by induction on the length of $\Gamma$, for all types simultaneously; for the
rank-zero type, $m_{\varnothing,\mathfrak t_\varnothing,\delta}=1$. In \eqref{eq:cGamma-first-step-type},
the induction hypothesis gives $m_{\Gamma_{>1}}(\pi')=m_{\Gamma_{>1},\mathfrak t(\mathfrak d,\mathbf n),\delta}(q)$,
where $(\mathfrak d,\mathbf n)$ are the step data of $\pi'$ relative to $\pi$. Grouping the terms by
$(\mathfrak d,\mathbf n)$ and using \eqref{eq:branching-between-types}, we get $m_\Gamma(\pi)=m_{\Gamma,\mathfrak t,\delta}(q)$ for
\begin{equation}\label{eq:cGamma-expanded-pattern-recursion}
 m_{\Gamma,\mathfrak t,\delta}(X)
 =\sum_{\mathfrak d,\mathbf n}
 m_{2\ell_1,a_1}(\mathfrak t;\mathfrak d,\mathbf n)
 N_{\mathfrak d,\mathbf n}(X)\,
 m_{\Gamma_{>1},\mathfrak t(\mathfrak d,\mathbf n),\delta}(X).
\end{equation}
The sum is finite, since $r(\mathfrak d)\le N_1$ and \eqref{eq:target-profile-rank-split}
bound the step data. The step data, the integers $m_{2\ell_1,a_1}(\mathfrak t;\mathfrak d,\mathbf n)$
\textup{(}see \eqref{eq:MH-definition}\textup{)}, the polynomials $N_{\mathfrak d,\mathbf n}(X)$
\textup{(}see \eqref{eq:profile-count-product} and Proposition~\ref{prop:primary-counting}\textup{)} and,
by induction, the polynomials $m_{\Gamma_{>1},\mathfrak t(\mathfrak d,\mathbf n),\delta}(X)$ depend only on $\mathfrak t$ and $\delta$, not on $q$ or $\pi$.
To see that $m_{\Gamma,\mathfrak t,\delta}(X)$ is independent of $\delta$, for $e\in\{\pm1\}$ replace $\Lambda_{-1}$ in $\mathfrak t$ by $\Lambda_{-1}^{\langle e\rangle}$ to form $\mathfrak t^{\langle e\rangle}$, and set $\Gamma^{(e)}:=\llbracket(2\ell_1,ea_1),(2\ell_2,-ea_2),\ldots,(2\ell_r,(-1)^{r-1}ea_r)\rrbracket$. We check by induction on $r$ that $m_{\Gamma,\mathfrak t,\delta}(X)=m_{\Gamma^{(e)},\mathfrak t^{\langle e\rangle},-\delta}(X)$. For $r=0$ both sides are $1$. For $r\ge1$, replacing $\Lambda_{-1}^{\mathfrak d}$ by $(\Lambda_{-1}^{\mathfrak d})^{\langle -e\rangle}$ gives a bijection $(\mathfrak d,\mathbf n)\mapsto(\mathfrak d^{\langle -e\rangle},\mathbf n)$ preserving $N_{\mathfrak d,\mathbf n}(X)$ and $\delta_{-1},\delta'_{-1}$, with $t'=et$ in \eqref{eq:MH-definition}; since $(-e)(-\delta)t'=\delta t$ and $et'=t$, we have $m_{2\ell_1,ea_1}(\mathfrak t^{\langle e\rangle};\mathfrak d^{\langle -e\rangle},\mathbf n)=m_{2\ell_1,a_1}(\mathfrak t;\mathfrak d,\mathbf n)$. Since $(\Gamma^{(e)})_{>1}=(\Gamma_{>1})^{(-e)}$, \eqref{eq:cGamma-expanded-pattern-recursion} completes the induction. Setting $e=+1$ gives $m_{\Gamma,\mathfrak t,+1}(X)=m_{\Gamma,\mathfrak t,-1}(X)=:m_{\Gamma,\mathfrak t}(X)$.
\end{proof}
By iterating \eqref{eq:cGamma-expanded-pattern-recursion}, we express
$m_{\Gamma,\mathfrak t}(X)$ as a sum over paths of types. For each type
$\mathfrak u=(\{(d_i,\kappa_i,\lambda_i)\}_{i\in I_{\mathfrak u}};\Lambda_1,\Lambda_{-1})$ we fix the index
set $I_{\mathfrak u}$, and for each irreducible $\sigma$ with type $\mathfrak u$ we write its classes as $[b_i]$,
$i\in I_{\mathfrak u}$, as in \eqref{eq:qtype-primary}.
A \emph{labelled path of types} for $\Gamma$, not to be confused with a descent path
$\mathcal P$ of representations, is a sequence
\[
 \mathfrak p:\quad
 \mathfrak t^{(0)}
 \xrightarrow{(\mathfrak d_1,\mathbf n_1)}\mathfrak t^{(1)}
 \xrightarrow{(\mathfrak d_2,\mathbf n_2)}\cdots
 \xrightarrow{(\mathfrak d_r,\mathbf n_r)}\mathfrak t^{(r)}=\mathfrak t_\varnothing,
\]
in which $\mathfrak t^{(0)}=\mathfrak t$ and each
$\mathfrak t^{(j)}=\mathfrak t(\mathfrak d_j,\mathbf n_j)$ is obtained from
$\mathfrak t^{(j-1)}$ by the step data $(\mathfrak d_j,\mathbf n_j)$. The edge labels are part of $\mathfrak p$: two paths with the same vertex
types are distinct if their indexed step data differ. For a descent path
\[
 \mathcal P:\quad \pi=\pi^{(0)}\searrow_{2\ell_1,a_1}\pi^{(1)}\searrow\cdots
 \searrow_{2\ell_r,a_r}\pi^{(r)}=\mathbf 1_{\Sp_0}
\]
as in Section~\ref{subsec:good-indexes}, write $\mathfrak p(\mathcal P)$ for its labelled path of
types. Iterating \eqref{eq:cGamma-first-step-type} and grouping the descent paths by
$\mathfrak p(\mathcal P)$, we obtain
\begin{equation}\label{eq:grouped-path-weight}
 m_\Gamma(\pi)=\sum_{\mathfrak p}T_{\mathfrak p}(\pi),\qquad
 T_{\mathfrak p}(\pi):=\sum_{\mathfrak p(\mathcal P)=\mathfrak p}\ \prod_{j=1}^r
 m_{2\ell_j,a_j}\bigl(\pi^{(j-1)},\pi^{(j)}\bigr).
\end{equation}
By \eqref{eq:MH-definition} and \eqref{eq:profile-count-product}, for every realization
$\sigma\in\Irr_{\mathfrak t^{(j-1)}}(\Sp_{2N_{j-1}}(\Fq))$,
we have
\begin{equation}\label{eq:fixed-step-data-weight}
\begin{aligned}
 &\sum_{\substack{\sigma'\in\Irr(\Sp_{2N_j}(\Fq))\\
       \text{with step data }(\mathfrak d_j,\mathbf n_j)
           \text{ relative to }\sigma}}
       m_{2\ell_j,a_j}(\sigma,\sigma')\\
 &\qquad= {m_{2\ell_j,a_j}\bigl(\mathfrak t^{(j-1)};\mathfrak d_j,\mathbf n_j\bigr)}
       N_{\mathfrak d_j,\mathbf n_j}(q).
\end{aligned}
\end{equation}
The sum depends on $\sigma$ only through its type.
Summing in \eqref{eq:grouped-path-weight} over $\pi^{(r)},\pi^{(r-1)},\ldots,\pi^{(1)}$ in
turn and applying \eqref{eq:fixed-step-data-weight} at each step, we get
$T_{\mathfrak p}(\pi)=T_{\mathfrak p}(q)$, where
\begin{equation}\label{eq:weighted-descent-paths}
 \begin{gathered}
 T_{\mathfrak p}(X):=\prod_{j=1}^r
 m_{2\ell_j,a_j}\bigl(\mathfrak t^{(j-1)};\mathfrak d_j,\mathbf n_j\bigr)
 N_{\mathfrak d_j,\mathbf n_j}(X),\\
 m_{\Gamma,\mathfrak t}(X)=\sum_{\mathfrak p}T_{\mathfrak p}(X);
 \end{gathered}
\end{equation}
by induction using \eqref{eq:cGamma-expanded-pattern-recursion}, we obtain the second identity.
Here $\mathfrak p$ runs over the labelled paths of types for the fixed index $\Gamma$.
By \eqref{eq:profile-count-product} and \eqref{eq:AB-leading},
$N_{\mathfrak d_j,\mathbf n_j}(X)$ has leading coefficient
$\prod_{d,\kappa}(2d)^{-K_{d,\kappa}(\mathbf n_j)}/\prod_e n^{(j)}_{d,\kappa,e}!>0$, so every nonzero
$T_{\mathfrak p}(X)$ has positive leading coefficient and the leading terms cannot cancel.
Consequently, if $m_{\Gamma,\mathfrak t}(X)\ne0$, then $\deg m_{\Gamma,\mathfrak t}(X)$ is the
maximum of $\deg T_{\mathfrak p}(X)$ over the nonzero $T_{\mathfrak p}(X)$.

To identify the terms of largest degree, we compare $\Gamma$ with the maximal descent
sequence indices of $\pi$.

Let $\Gamma^{\max}$ be maximal with respect to $\le_{\rm seq}$ among the descent sequence
indices realized by $\pi$. By \cite[Theorem~7.1]{wang}, $\Gamma^{\max}$ is good and every descent
sequence realizing it consists of first descents; hence, by \eqref{eq:first-descent-inherited-rank},
\begin{equation}\label{eq:max-index-num}
 (\Gamma^{\max})^{\rm num}=\bigl(2(r_{j-1}^{\rm fd}-r_j^{\rm fd})\bigr)_{j\ge1},\qquad
 r_j^{\rm fd}=r_j^{\rm fd}(\mathfrak t).
\end{equation}
Since each elementary downward mutation raises $\sum_j j\ell_j$ by one, strictly increasing
chains for $\le_{\rm seq}$ are finite, so every descent sequence index realized by $\pi$ lies below
some such $\Gamma^{\max}$.

\subsection{Highest-degree terms in the path expansion}\label{subsec:budget-sign-average}
We now bound the degree of each term $T_{\mathfrak p}(X)$ in
\eqref{eq:weighted-descent-paths}. 
Fix a descent sequence index $\Gamma$ of length $r$ with $\sum_{j=1}^r\ell_j=N$.
Along a descent path with index $\Gamma$, the rank after $j$ steps is $N_j$ below; we compare
it with the rank $r_j^{\rm fd}(\mathfrak t)$ after $j$ successive first descents
\textup{(}see \eqref{eq:first-descent-inherited-rank}\textup{)}:
\begin{equation}\label{eq:sequence-budget-main}
\begin{aligned}
 N_j&=N-\sum_{i\le j}\ell_i,\\
 R_j=R_j(\Gamma,\mathfrak t)&:=N_j-r_j^{\rm fd}(\mathfrak t)
                 \qquad(1\le j<r),\\
 B_{\Gamma,\mathfrak t}&:=\sum_{j<r}R_j.
\end{aligned}
\end{equation}
In Proposition~\ref{prop:degree-budget}, we show that $B_{\Gamma,\mathfrak t}$
bounds the degree  of $m_{\Gamma,\mathfrak t}(X)$.

\begin{lemma}
\label{lem:fresh-rank}
For every pair $(\mathfrak d,\mathbf{n})$,
\[
 \deg N_{\mathfrak d,\mathbf{n}}(X)=\sum_{d,\kappa,e}d n_{d,\kappa,e}\le m(\mathfrak d),
\]
with equality if and only if $n_{d,\kappa,e}=0$ for all $e>1$.
\end{lemma}
\begin{proof}
In \eqref{eq:profile-count-product}, the integers $h_{d,\kappa}$, $K_{d,\kappa}(\mathbf n)$ and
$n_{d,\kappa,e}$ are combinatorial data independent of $q$, and $P_d^\kappa(X)$ has degree $d$ by
\eqref{eq:AB-leading}. Hence $(P_d^\kappa(X)-h_{d,\kappa})_{\underline{K_{d,\kappa}(\mathbf n)}}$ has
degree $dK_{d,\kappa}(\mathbf n)$, and

$\deg N_{\mathfrak d,\mathbf n}=\sum_{d,\kappa}dK_{d,\kappa}(\mathbf n)
=\sum_{d,\kappa,e}dn_{d,\kappa,e}$. By
\eqref{eq:target-profile-rank-split}, we obtain
\[
 m(\mathfrak d)-\deg N_{\mathfrak d,\mathbf n}
 =\sum_{d,\kappa,e}d(e-1)n_{d,\kappa,e}.
\]
The difference vanishes exactly when every new class carries $(1^e)$ with
$e=1$.
\end{proof}

In the next lemma, we compare the ranks after an arbitrary descent step with those
after a first descent.

\begin{lemma}
\label{lem:minimal-inherited-tail-equality}
Let $(\mathfrak d,\mathbf n)$ be step data for a
$(2\ell,a)$-descent from a type $\mathfrak t=(\{(d_i,\kappa_i,\lambda_i)\}_{i\in I};\Lambda_{\pm1})$
to $\mathfrak t(\mathfrak d,\mathbf n)$, where $\mathfrak d=(\{\lambda_i^{\mathfrak d}\}_{i\in I};\Lambda^{\mathfrak d}_{\pm1})$, with $\ell\ge0$, and suppose
that $m_{2\ell,a}(\mathfrak t;\mathfrak d,\mathbf n)\ne0$.

\textup{(a)} We have, with $r(\mathfrak d)$ as in \eqref{eq:H-inherited-rank},
\begin{equation}\label{eq:type-one-step}
 r(\mathfrak d)\ge r_1^{\rm fd}(\mathfrak t),\qquad
 r_k^{\rm fd}(\mathfrak t(\mathfrak d,\mathbf n))\ge r_{k+1}^{\rm fd}(\mathfrak t)\quad(k\ge1).
\end{equation}

\textup{(b)} If $r(\mathfrak d)=r_1^{\rm fd}(\mathfrak t)$, then $\lambda_i^{\mathfrak d}=(\lambda_i)_-$ for all $i\in I$ \textup{(}see \eqref{eq:first-column}\textup{)}, the multiplicity $m_{2\ell,a}(\mathfrak t;\mathfrak d,\mathbf n)$ is $1$ \textup{(}cf. Proposition~\ref{cor:first-descent-rule}\textup{)}, and
\begin{equation*}
 r_j^{\rm fd}(\mathfrak t(\mathfrak d,\mathbf n))=r_{j+1}^{\rm fd}(\mathfrak t)
 \qquad(j\ge1).
\end{equation*}
\end{lemma}
\begin{proof}
\textup{(a)} As $m_{2\ell,a}(\mathfrak t;\mathfrak d,\mathbf n)\ne0$, every factor
on the right of \eqref{eq:MH-definition} is nonzero. Thus
$\lambda_i$ and $\lambda_i^{\mathfrak d}$ are close, so that each part of
$\lambda_i^{\mathfrak d}$ is at least the corresponding part of $(\lambda_i)_-$;
and $\mathbf b_t(\Lambda_{\pm1};\Lambda^{\mathfrak d}_{\pm1})\ne0$ for the value of
$t$ defined there, so that
$r_k^{\rm fd}(\Lambda^{\mathfrak d}_{-\xi})\ge r_{k+1}^{\rm fd}(\Lambda_\xi)$ for
$k\ge0$ and $\xi=\pm1$ by Lemma~\ref{lem:symbol-first-descent-comparison}.
Multiplying the partition bounds by $d_i$ and adding the symbol bounds, we obtain
the first inequality in \eqref{eq:type-one-step} \textup{(}recall
\eqref{eq:first-descent-inherited-rank}\textup{)}. The second follows from the same
bounds, because for $k\ge1$ the columns $(1^e)$ at the new classes of
$\mathfrak t(\mathfrak d,\mathbf n)$ do not contribute to $r_k^{\rm fd}(\mathfrak t(\mathfrak d,\mathbf n))$.

\textup{(b)} Remove the new classes from $\mathfrak t(\mathfrak d,\mathbf n)$ and denote the resulting
type by $\mathfrak t_{\mathfrak d}$. Let $\rho,\rho'$ be the tensor products of unipotent factors \eqref{eq:rho-def} of
representations realizing $\mathfrak t,\mathfrak t(\mathfrak d,\mathbf n)$, respectively. Then
\[
 \epsilon_{\rho'}=\epsilon_{\rho_{\mathfrak d}}\eta(\mathbf{n}),\qquad
 a_{\mathfrak d}:=a\eta(\mathbf{n}),\qquad
 a_{\mathfrak d}\epsilon_\rho\epsilon_{\rho_{\mathfrak d}}
 =a\epsilon_\rho\epsilon_{\rho'}.
\]
For each new class $[b]$ carrying $(1^e)$,
$b_{\kappa([b])}(\varnothing,(1^e))=1$ by
\eqref{eq:GGP-one-sided-column}. Deleting these classes therefore leaves
the product over $[x]\ne[\pm1]$ in \eqref{eq:GGP-primary-numerical}
unchanged. With sign $a_{\mathfrak d}$, the value of $t$, and hence
$\mathbf b_t(\Lambda_{\pm1};\Lambda^{\mathfrak d}_{\pm1})$, is also
unchanged. The remaining
parameters $(\rho_{\mathfrak d},\Lambda_1^{\mathfrak d},\Lambda_{-1}^{\mathfrak d})$
therefore give a constituent at descent index $\rank(\mathfrak t)-r(\mathfrak d)$.
Since $r(\mathfrak d)=r_1^{\rm fd}(\mathfrak t)$, this is a first
descent. By \eqref{eq:first-primary},
\[
 \lambda_i^{\mathfrak d}=(\lambda_i)_-,\qquad (1^e)_-=\varnothing.
\]
Consequently,
\[
 r_j^{\rm fd}(\mathfrak t(\mathfrak d,\mathbf n))=r_j^{\rm fd}(\mathfrak t_{\mathfrak d})
 =r_{j+1}^{\rm fd}(\mathfrak t)\qquad(j\ge1).
\]
By Proposition~\ref{cor:first-descent-rule} and
\eqref{eq:GGP-one-sided-column}, every factor in
\eqref{eq:MH-definition} equals one. Hence
${m_{2\ell,a}(\mathfrak t;\mathfrak d,\mathbf n)}=1$.
The same argument applies when $\mathfrak t$ and $\mathfrak t(\mathfrak d,\mathbf n)$ have equal rank.
\end{proof}

Accordingly, we say that step data $(\mathfrak d,\mathbf n)$ as in
Lemma~\ref{lem:minimal-inherited-tail-equality} are a \emph{first descent of the type
$\mathfrak t$} if $r(\mathfrak d)=r_1^{\rm fd}(\mathfrak t)$.

\begin{proposition}\label{prop:degree-budget}
Let $\Gamma$ be a descent sequence index with $\sum_{j=1}^r\ell_j=N$. If
$m_{\Gamma,\mathfrak t}(X)\ne0$, then
\begin{equation*}
 \deg m_{\Gamma,\mathfrak t}(X)\le B_{\Gamma,\mathfrak t}
 =\sum_{j=1}^{r-1}\bigl(N_j-r_j^{\rm fd}(\mathfrak t)\bigr).
\end{equation*}
Moreover, let $\mathfrak p:\ \mathfrak t=\mathfrak t^{(0)}\xrightarrow{(\mathfrak d_1,\mathbf n_1)}\cdots
\xrightarrow{(\mathfrak d_r,\mathbf n_r)}\mathfrak t^{(r)}=\mathfrak t_\varnothing$ be a labelled path
of types with $T_{\mathfrak p}(X)\ne0$, where
$\mathbf n_j=(n^{(j)}_{d,\kappa,e})$.
Then $\deg T_{\mathfrak p}(X)=B_{\Gamma,\mathfrak t}$ if and only if, for every $1\le j<r$,
\[
\begin{aligned}
 r(\mathfrak d_j)&=r_j^{\rm fd}(\mathfrak t),\\
 n^{(j)}_{d,\kappa,e}&=0\quad(e>1).
\end{aligned}
\]
\end{proposition}
\begin{proof}
For the rank-zero type, $m_{\varnothing,\mathfrak t_\varnothing}=1$; assume $r\ge1$.
By \eqref{eq:weighted-descent-paths}, it suffices to bound $\deg T_{\mathfrak p}(X)$ for each
nonzero $T_{\mathfrak p}(X)$. By the definition \eqref{eq:weighted-descent-paths} of $T_{\mathfrak p}(X)$,
in which the multiplicities $m_{2\ell_j,a_j}(\mathfrak t^{(j-1)};\mathfrak d_j,\mathbf n_j)$ are
constants, and by Lemma~\ref{lem:fresh-rank},
\[
\deg T_{\mathfrak p}(X)=\sum_{j<r}\deg N_{\mathfrak d_j,\mathbf n_j}(X)
 =\sum_{j<r}\sum_{d,\kappa,e}d\,n^{(j)}_{d,\kappa,e}.
\]
 By \eqref{eq:H-inherited-rank} and
\eqref{eq:target-profile-rank-split}, $N_j=r(\mathfrak d_j)+\sum_{d,\kappa,e}de\,n^{(j)}_{d,\kappa,e}$.
Hence
\begin{equation}\label{eq:expanded-degree-loss}
\begin{aligned}
 B_{\Gamma,\mathfrak t}-\deg T_{\mathfrak p}(X)
 ={}&\sum_{j<r}\bigl(r(\mathfrak d_j)-r_j^{\rm fd}(\mathfrak t)\bigr)\\
 &+\sum_{j<r}\sum_{d,\kappa,e}
       d(e-1)n^{(j)}_{d,\kappa,e}.
\end{aligned}
\end{equation}
By Lemma~\ref{lem:minimal-inherited-tail-equality}\textup{(a)},
applied to the $j$th step, we have
\[
\begin{aligned}
 r(\mathfrak d_j)&\ge r_1^{\rm fd}(\mathfrak t^{(j-1)}),\\
 r_k^{\rm fd}(\mathfrak t^{(j)})
 &\ge r_{k+1}^{\rm fd}(\mathfrak t^{(j-1)})
 \qquad(k\ge1).
\end{aligned}
\]
Iterating yields
\begin{equation}\label{eq:inherited-rank-chain}
 r(\mathfrak d_j)
 \ge r_1^{\rm fd}(\mathfrak t^{(j-1)})
 \ge\cdots\ge r_j^{\rm fd}(\mathfrak t).
\end{equation}
In particular, $r_r^{\rm fd}(\mathfrak t)=0$.
Both sums in \eqref{eq:expanded-degree-loss} are therefore termwise
nonnegative, and they vanish exactly when $r(\mathfrak d_j)=r_j^{\rm fd}(\mathfrak t)$ and
$n^{(j)}_{d,\kappa,e}=0$ for $e>1$, for every $1\le j<r$.
\end{proof}

\subsection{Leading coefficients}\label{subsec:sequence-leading-count}
By \eqref{eq:weighted-descent-paths}, $m_{\Gamma,\mathfrak t}(X)=\sum_{\mathfrak p}T_{\mathfrak p}(X)$.
By Proposition~\ref{prop:degree-budget}, the coefficient of $X^{B_{\Gamma,\mathfrak t}}$ comes from the
$\mathfrak p$ with $r(\mathfrak d_j)=r_j^{\rm fd}(\mathfrak t)$ and $n^{(j)}_{d,\kappa,e}=0$ for $e>1$
$(1\le j<r)$. Thus, writing $[X^d]P$ for the coefficient of $X^d$ in $P$,
\begin{equation}\label{eq:lc-sum-over-data}
 [X^{B_{\Gamma,\mathfrak t}}]m_{\Gamma,\mathfrak t}(X)
 =\sum_{\substack{\mathfrak p:\,T_{\mathfrak p}(X)\ne0\\
                  \deg T_{\mathfrak p}(X)=B_{\Gamma,\mathfrak t}}}
       \LC(T_{\mathfrak p}(X)).
\end{equation}
For such $\mathfrak p$, the inequalities in \eqref{eq:inherited-rank-chain} are equalities. Hence
$r(\mathfrak d_j)=r_1^{\rm fd}(\mathfrak t^{(j-1)})$, and by
Lemma~\ref{lem:minimal-inherited-tail-equality}\textup{(b)},
$m_{2\ell_j,a_j}(\mathfrak t^{(j-1)};\mathfrak d_j,\mathbf n_j)=1$
for every $1\le j\le r$.
By the definition \eqref{eq:weighted-descent-paths} of $T_{\mathfrak p}(X)$,
\eqref{eq:profile-count-product} and \eqref{eq:AB-leading}, we obtain
\begin{equation}\label{eq:lc-single-data-weight}
 \LC(T_{\mathfrak p}(X))
 =\prod_{j<r}\prod_{d,\kappa}
   \frac{1}{(2d)^{n^{(j)}_{d,\kappa,1}}\,
                 n^{(j)}_{d,\kappa,1}!}.
\end{equation}
We next determine which first-descent paths occur.

Fix $\pi\in\Irr_{\mathfrak t}(\Sp_{2N}(\Fq))$ and choose a maximal
descent sequence index $\Gamma^{\max}$ realized by $\pi$.
Appending entries $(0,a)$ as in Section~\ref{subsec:mutations}, we assume that $\Gamma^{\max}$
and the indices $\Gamma$ compared with it have a common length $r$. By
\eqref{eq:max-index-num},
\[
 (\Gamma^{\max})^{\rm num}=(2\ell_1^{\max},\ldots,2\ell_r^{\max}),
 \qquad \ell_j^{\max}=r_{j-1}^{\rm fd}-r_j^{\rm fd}.
\]
For an index
$\Gamma=\llbracket(2\ell_1,a_1),\ldots,(2\ell_r,a_r)\rrbracket$, define
\[
 \sigma_i(\Gamma):=\prod_{j\le i}a_j,\qquad \sigma_0(\Gamma):=1.
\]
We call $\sigma_i(\Gamma)$ the $i$th \emph{prefix sign} of $\Gamma$.

Since $\sum_j\ell_j=N$ and $r_i^{\rm fd}=N-\sum_{j\le i}\ell_j^{\max}$, we have
\begin{equation}\label{eq:mi-Ri}
 R_i=\sum_{j\le i}(\ell_j^{\max}-\ell_j)\qquad(1\le i<r).
\end{equation}

\begin{lemma}\label{lem:prefix-criterion}
Let $\Gamma$ be an index of length $r$, and set
$R_0=R_r:=0$. Then
\begin{equation}\label{eq:signed-prefix-reachability}
\begin{gathered}
 \Gamma\le_{\rm seq}\Gamma^{\max}\quad\Longleftrightarrow\\
 \begin{cases}
 R_i\ge0,&1\le i<r,\\
 \sigma_i(\Gamma)=\sigma_i(\Gamma^{\max}),&R_i=0,\ 1\le i\le r.
 \end{cases}
\end{gathered}
\end{equation}
\end{lemma}
\begin{proof}
This follows directly from \eqref{eq:mi-Ri} and \eqref{eq:elementary-seq-mutation}.
\end{proof}

Assume now that $R_i\ge0$ for $i<r$, and let
\begin{equation}\label{eq:S-maxnum-Gamma}
 \mathscr S_\Gamma^{\pi}:=
 \{\Gamma':(\Gamma')^{\rm num}=(\Gamma^{\max})^{\rm num},
                     \ \Gamma\le_{\rm seq}\Gamma'\}.
\end{equation}
By Lemma~\ref{lem:prefix-criterion},
\begin{equation}\label{eq:S-sign-parametrization}
 \mathscr S_\Gamma^{\pi}\simeq\{\pm1\}^{h},\qquad
 \Gamma'\longmapsto(\sigma_i(\Gamma'))_{R_i>0},
 \qquad h=\#\{i<r:R_i>0\}.
\end{equation}

Let $\Gamma'\in\mathscr S_\Gamma^\pi$. Then $B_{\Gamma',\mathfrak t}=0$, so by
Proposition~\ref{prop:degree-budget} every $\mathfrak p$ along $\Gamma'$ with
$T_{\mathfrak p}(X)\ne0$ consists of first descents without new classes. Hence
$N_{\mathfrak d_j,\mathbf n_j}(X)=1$ and, by Lemma~\ref{lem:minimal-inherited-tail-equality}\textup{(b)},
$T_{\mathfrak p}(X)=1$; by \eqref{eq:grouped-path-weight}, $\mathfrak p=\mathfrak p(\mathcal P)$ for
exactly one descent path $\mathcal P$ from $\pi$. Thus $m_{\Gamma',\mathfrak t}(X)$ is the number of
descent paths from $\pi$ to $\mathbf 1_{\Sp_0}$ with signed index $\Gamma'$. Let
$\mathscr T_{\Gamma,\mathfrak t}$ be the set of these paths for all $\Gamma'\in\mathscr S_\Gamma^\pi$, so that
\[
 |\mathscr T_{\Gamma,\mathfrak t}|
 =\sum_{\Gamma'\in\mathscr S_\Gamma^\pi}m_{\Gamma',\mathfrak t}(X).
\]

\begin{proposition}
\label{cor:cGamma-lc-explicit}
Assume $R_i\ge0$ for $i<r$ and $r_r^{\rm fd}=0$,
that is, $r$ is at least the length of a first-descent path from $\pi$ to
$\mathbf{1}_{\Sp_0}$. Then
\begin{equation}\label{eq:budget-chain-count}
 [X^{B_{\Gamma,\mathfrak t}}]m_{\Gamma,\mathfrak t}(X)
 =\frac1{2^h}\sum_{\Gamma'\in\mathscr S_\Gamma^\pi}m_{\Gamma',\mathfrak t}(X)
 =2^{-h}|\mathscr T_{\Gamma,\mathfrak t}|.
\end{equation}
Moreover,
\begin{equation*}
 m_{\Gamma,\mathfrak t}(X)\ne0
 \quad\Longleftrightarrow\quad\mathscr T_{\Gamma,\mathfrak t}\ne\varnothing.
\end{equation*}
If the polynomial is nonzero, its degree is $B_{\Gamma,\mathfrak t}$, and its
leading coefficient is given by \eqref{eq:budget-chain-count}.
\end{proposition}

\begin{proof}
\leavevmode By \eqref{eq:lc-sum-over-data},
\[
 [X^{B_{\Gamma,\mathfrak t}}]m_{\Gamma,\mathfrak t}(X)
 =\sum_{\substack{\mathfrak p:\,T_{\mathfrak p}(X)\ne0\\
                  \deg T_{\mathfrak p}(X)=B_{\Gamma,\mathfrak t}}}
       \LC(T_{\mathfrak p}(X)).
\]
We count the labelled paths $\mathfrak p$ in this sum with
weights $\LC(T_{\mathfrak p}(X))$. Consider a descent path
\[
 \mathcal P:\quad
 \pi=\pi^{(0)}\searrow_{2\ell_1,a_1}\pi^{(1)}
 \searrow_{2\ell_2,a_2}\cdots
 \searrow_{2\ell_r,a_r}\pi^{(r)}=\mathbf 1_{\Sp_0},
\]
whose labelled path of types $\mathfrak p=\mathfrak p(\mathcal P)$ satisfies
$T_{\mathfrak p}(X)\ne0$ and $\deg T_{\mathfrak p}(X)=B_{\Gamma,\mathfrak t}$.
At the $j$th step, $\pi^{(j)}\in\Irr(\Sp_{2N_j}(\Fq))$ is an irreducible
constituent of $\cD^{\FJ}_{2\ell_j,a_j}(\pi^{(j-1)})$, where
$N_j=N_{j-1}-\ell_j$. For $1\le j<r$, let $S_j$ be the set of classes
$[x]\ne[\pm1]$ occurring in the parameter of $\pi^{(j)}$ such that $[-x]$ does
not occur in the parameter of $\pi^{(j-1)}$. Thus $S_j$ consists of the new
classes at the $j$th step, as in Section~\ref{subsec:type-transitions}. Set
$S_0=S_r=\varnothing$. Recall from \eqref{eq:sequence-budget-main} that $R_j=N_j-r_j^{\rm fd}(\mathfrak t)$.
By Proposition~\ref{prop:degree-budget} and \eqref{eq:target-profile-rank-split},
every class in $S_j$ carries $(1)$ and
\[
 \sum_{[x]\in S_j}d([x])=R_j\qquad(1\le j<r).
\]
Let $\mathbf n_j=(n^{(j)}_{d,\kappa,e})$ record the new classes in $S_j$. Using
$\eta(\mathbf n)$ from Section~\ref{subsec:type-transitions}, set
\[
 \eta_j:=\eta(\mathbf n_j)
 =(-1)^{\#\{[x]\in S_j:\kappa([x])=\U\}}
 \qquad(1\le j<r).
\]
For the empty sets $S_0$ and $S_r$, set $\eta_0=\eta_r=1$.

\leavevmode For $0\le j\le r$, let $\pi_{\rm fd}^{(j)}$ be obtained from
$\pi^{(j)}$ by deleting $S_j$, and write $\rho^{(j)}$ and $\rho_{\rm fd}^{(j)}$
for their respective factors \eqref{eq:rho-def}. Their $[\pm1]$-symbols agree, and
\[
 \pi_{\rm fd}^{(j)}\in\Irr(\Sp_{2r_j^{\rm fd}}(\Fq)),\qquad
 \epsilon_{\rho_{\rm fd}^{(j)}}=\eta_j\epsilon_{\rho^{(j)}}.
\]
Put $\Gamma'=\llbracket(2\ell_1^{\max},b_1),\ldots,(2\ell_r^{\max},b_r)\rrbracket$,
where
\begin{equation*}
 b_j:=a_j\eta_{j-1}\eta_j,\qquad
 \sigma_j(\Gamma')=\sigma_j(\Gamma)\eta_j\quad(1\le j\le r).
\end{equation*}
By \eqref{eq:GGP-rational-sign}, $t$ is unchanged. By Lemma~\ref{lem:minimal-inherited-tail-equality}\textup{(b)}, the classes in $S_j$ disappear at the next step, and
\[
 \mathcal P_{\rm fd}:\quad
 \pi=\pi_{\rm fd}^{(0)}\fdsearrow_{2\ell_1^{\max},b_1}\pi_{\rm fd}^{(1)}
 \fdsearrow\cdots\fdsearrow_{2\ell_r^{\max},b_r}\pi_{\rm fd}^{(r)}=\mathbf 1_{\Sp_0}.
\]
Since $R_j=0$ implies $\eta_j=1$, and $\eta_r=1$, we have
\[
 \sigma_j(\Gamma')=\sigma_j(\Gamma)\quad(R_j=0,\ 1\le j<r),\qquad
 \sigma_r(\Gamma')=\sigma_r(\Gamma).
\]
By Lemma~\ref{lem:prefix-criterion}, $\Gamma'\in\mathscr S_\Gamma^\pi$.
Thus every descent path $\mathcal P$ whose labelled type path
$\mathfrak p(\mathcal P)$ indexes a term in the sum on the right-hand side of
\eqref{eq:lc-sum-over-data} determines a first-descent path
$\mathcal P_{\rm fd}\in\mathscr T_{\Gamma,\mathfrak t}$.

\leavevmode Conversely, every
$\mathcal P_{\rm fd}\in\mathscr T_{\Gamma,\mathfrak t}$ arises from such a descent
path $\mathcal P$ for all sufficiently large $q$.
Fix $\mathcal P_{\rm fd}\in\mathscr T_{\Gamma,\mathfrak t}$ with signed
index $\Gamma'$, and put $\eta_j:=\sigma_j(\Gamma)\sigma_j(\Gamma')$ for
$0\le j\le r$. By Lemma~\ref{lem:prefix-criterion},
\[
 \eta_0=\eta_r=1,\qquad \eta_j=1\quad\text{when }R_j=0\quad(1\le j<r).
\]
Starting with $\pi^{(0)}=\pi$, we construct $\pi^{(j)}$ successively for
$1\le j<r$. Let $E_j$ be the non-$[\pm1]$ classes occurring in
$\pi^{(0)},\ldots,\pi^{(j-1)}$ and their negatives. Thus $|E_j|\le2rN$. Choose a
set $S_j$ of non-$[\pm1]$ classes such that
\[
 S_j\cap E_j=\varnothing,\qquad
 \sum_{[a]\in S_j}d([a])=R_j,\qquad
 (-1)^{\#\{[a]\in S_j:\kappa([a])=\U\}}=\eta_j.
\]
By Lemma~\ref{lem:singleton-count}\textup{(b)}, such sets exist for sufficiently
large $q$, and excluding $E_j$ does not affect the leading coefficient. Let
$\pi^{(j)}$ be obtained from $\pi_{\rm fd}^{(j)}$ by adjoining the classes in
$S_j$, each carrying $(1)$, and set $\pi^{(r)}=\mathbf1_{\Sp_0}$. We have
\[
 b_\kappa((1),\varnothing)=b_\kappa(\varnothing,(1))=1
 \quad(\kappa\in\{\GL,\U\}),\qquad a_j=b_j\eta_{j-1}\eta_j.
\]
Thus, by \eqref{eq:GGP-rational-sign} and \eqref{eq:GGP-primary-numerical}, we have
$\pi^{(j-1)}\searrow_{2\ell_j,a_j}\pi^{(j)}$ for $1\le j\le r$, and the resulting
path has index $\Gamma$ and first-descent path $\mathcal P_{\rm fd}$. It is uniquely
determined by the sets $S_j$.
By \eqref{eq:GGP-primary-numerical}, with $\ell_j=0$, each appended
entry $(0,a_j)$ is counted once with its prescribed sign.

\leavevmode We now count the descent paths $\mathcal P$
whose labelled type paths index the sum in \eqref{eq:lc-sum-over-data} and whose
associated first-descent path is the fixed $\mathcal P_{\rm fd}$.
Let $K_j$ record the degrees and kinds of the classes added back at
stage $j$. Write $|K_j|$ for the total rank, $\eta(K_j)$ for the parity
sign, and $w(K_j)$ for the factor at stage $j$ in \eqref{eq:lc-single-data-weight}.
Let $\alpha_R(\eta)$ be the sum of the weights over tuples of rank $R$ and parity
sign $\eta$. By \eqref{eq:lc-single-data-weight}, the total contribution to
the sum in \eqref{eq:lc-sum-over-data} from the labelled type paths
$\mathfrak p(\mathcal P)$ associated with the fixed $\mathcal P_{\rm fd}$ is
\begin{equation*}
 \sum_{\substack{|K_j|=R_j,\ \eta(K_j)=\eta_j\\1\le j<r}}
       \prod_{j<r}w(K_j)
 =\prod_{j<r}\alpha_{R_j}(\eta_j)
 =2^{-h},
\end{equation*}
where we obtain the last equality by Lemma~\ref{lem:singleton-count}\textup{(a)}.
By Lemma~\ref{lem:singleton-count}\textup{(b)}, the number of
such descent paths $\mathcal P$ corresponding to this fixed $\mathcal P_{\rm fd}$ is
$2^{-h}q^{B_{\Gamma,\mathfrak t}}+O(q^{B_{\Gamma,\mathfrak t}-1})$.
There are $m_{\Gamma',\mathfrak t}(X)$ first-descent paths with signed index $\Gamma'$, and
adding up we obtain
\begin{equation}\label{eq:direct-primary-sign-average}
\begin{aligned}
 [X^{B_{\Gamma,\mathfrak t}}]m_{\Gamma,\mathfrak t}(X)
 &=\sum_{\Gamma'\in\mathscr S_\Gamma^\pi}m_{\Gamma',\mathfrak t}(X)
    \prod_{j<r}\alpha_{R_j}
       \bigl(\sigma_j(\Gamma)\sigma_j(\Gamma')\bigr)\\
 &=2^{-h}\sum_{\Gamma'\in\mathscr S_\Gamma^\pi}m_{\Gamma',\mathfrak t}(X),
\end{aligned}
\end{equation}
\leavevmode If $m_{\Gamma,\mathfrak t}(X)\ne0$, then $\Gamma$ is realized
for sufficiently large $q$. By the discussion around
\eqref{eq:max-index-num}, $\Gamma$ lies below a realized maximal index
$\Gamma'\in\mathscr S_\Gamma^\pi$, which is realized by a first-descent path.
Thus $\mathscr T_{\Gamma,\mathfrak t}\ne\varnothing$. Conversely, if
$\mathscr T_{\Gamma,\mathfrak t}\ne\varnothing$, then
by \eqref{eq:direct-primary-sign-average} and Proposition~\ref{prop:degree-budget},
we have $\deg m_{\Gamma,\mathfrak t}(X)=B_{\Gamma,\mathfrak t}$.
\end{proof}

\begingroup\samepage
\begin{lemma}\label{lem:singleton-count}
For a tuple $K=(K_{d,\kappa})_{d\ge1,\,\kappa\in\{\GL,\U\}}$ of nonnegative integers, put
\begin{equation*}
\begin{gathered}
 |K|:=\sum_{d,\kappa}dK_{d,\kappa},\qquad
 \eta(K):=(-1)^{\sum_dK_{d,\U}},\\
 w(K):=\prod_{d,\kappa}
       \frac{1}{(2d)^{K_{d,\kappa}}K_{d,\kappa}!},
\end{gathered}
\end{equation*}
and $\alpha_R(\eta):=\sum_{|K|=R,\,\eta(K)=\eta}w(K)$ for $R\ge0$ and $\eta=\pm1$.

\textup{(a)} We have
\begin{equation}\label{eq:singleton-parity-mass}
 \alpha_R(\eta)
 =\begin{cases}
  1,&R=0,\ \eta=+1,\\
  0,&R=0,\ \eta=-1,\\
  \tfrac12,&R>0.
 \end{cases}
\end{equation}

\textup{(b)} Let $E=E(q)$ be a set of classes $[x]\ne[\pm1]$ with $|E|\le M$, where $M$ is independent of $q$. The number of sets $S$ of classes $[x]\ne[\pm1]$, $[x]\notin E$, with
$\sum_{[x]\in S}d([x])=R$ and $(-1)^{\#\{[x]\in S:\,\kappa([x])=\U\}}=\eta$ is
\begin{equation}\label{eq:FR-leading}
 F_R^\eta(q)=\alpha_R(\eta)\,q^R+O_{R,M}(q^{R-1}).
\end{equation}
\end{lemma}
\endgroup
\begin{proof}
(a) As formal power series in $x$,
\[
 \sum_{K}w(K)\,\eta(K)^{\epsilon}x^{|K|}
 =\prod_{d\ge1}\exp\Bigl(\frac{x^d}{2d}\Bigr)
   \exp\Bigl((-1)^{\epsilon}\frac{x^d}{2d}\Bigr)
 =\begin{cases}(1-x)^{-1},&\epsilon=0,\\ 1,&\epsilon=1.\end{cases}
\]
Comparing coefficients gives $\sum_{|K|=R}w(K)=1$ for every $R\ge0$ and
$\sum_{|K|=R}\eta(K)w(K)=0$ for $R>0$, whence \eqref{eq:singleton-parity-mass}.

(b) Let $e_{d,\kappa}$ be the number of classes in $E$ of degree $d$ and kind $\kappa$. Since $e_{d,\kappa}\le M$, the estimates below are uniform in $E$. The sets
$S$ with exactly $K_{d,\kappa}$ classes of degree $d$ and kind $\kappa$ number
$\prod_{d,\kappa}\binom{P_d^\kappa(q)-e_{d,\kappa}}{K_{d,\kappa}}$. By \eqref{eq:AB-leading}, each
factor is $\frac1{K_{d,\kappa}!}\bigl(\frac{q^d}{2d}\bigr)^{K_{d,\kappa}}+O(q^{dK_{d,\kappa}-1})$, so the
product is $w(K)\,q^{|K|}+O(q^{|K|-1})$. Summing over the finitely many $K$ with $|K|=R$ and
$\eta(K)=\eta$ and applying \textup{(a)} gives \eqref{eq:FR-leading}.
\end{proof}

\subsection{Behaviour under elementary mutations}\label{subsec:local-transfer}
We use  Proposition~\ref{cor:cGamma-lc-explicit} to compare the leading coefficients before and after the mutation.

\leavevmode Let $\Gamma=\llbracket(2\ell_1,a_1),\ldots,(2\ell_r,a_r)\rrbracket$
and fix $1\le k<r$ with $\ell_k>0$. Put $\Gamma_k^+=\Gamma$. For
$\varepsilon\in\{\pm1\}$, let $\Gamma_k^-$ and $\Gamma_{k,\varepsilon}^0$ agree
with $\Gamma$ outside positions $k,k+1$, where their respective entries are
\[
\begin{aligned}
 \Gamma_k^-&:\quad\bigl((2\ell_k,-a_k),(2\ell_{k+1},-a_{k+1})\bigr),\\
 \Gamma_{k,\varepsilon}^0&:\quad
 \bigl((2\ell_k-2,\varepsilon a_k),(2\ell_{k+1}+2,\varepsilon a_{k+1})\bigr).
\end{aligned}
\]
Since $\varepsilon^2=1$, both choices preserve the product $a_ka_{k+1}$ of the
two signs, as in \eqref{eq:elementary-seq-mutation}. By \eqref{eq:sequence-budget-main}, $\Gamma_k^+$ and $\Gamma_k^-$ have
the same $R_i$, while $\Gamma_{k,\varepsilon}^0$ increases $R_k$ by one and
leaves the other $R_i$ unchanged. Hence
\begin{equation}\label{eq:budget-mutation}
 B_{\Gamma_{k,\varepsilon}^0,\mathfrak t}=B_{\Gamma_k^+,\mathfrak t}+1
 =B_{\Gamma_k^-,\mathfrak t}+1.
\end{equation}
\begin{lemma}
\label{thm:even-transfer}
\leavevmode For $1\le k<r$ with $\ell_k>0$ and $\varepsilon\in\{\pm1\}$,
let $\Gamma_k^+$, $\Gamma_k^-$ and $\Gamma_{k,\varepsilon}^0$ be as above.
Suppose that all $R_i(\Gamma,\mathfrak t)$ are nonnegative. Then
\begin{equation}\label{eq:two-step-budget-average}
 \begin{aligned}
 [X^{B_{\Gamma_{k,\varepsilon}^0,\mathfrak t}}]
 m_{\Gamma_{k,\varepsilon}^0,\mathfrak t}(X)
 &=\tfrac12[X^{B_{\Gamma_k^+,\mathfrak t}}]m_{\Gamma_k^+,\mathfrak t}(X)\\
 &\quad+\tfrac12[X^{B_{\Gamma_k^-,\mathfrak t}}]m_{\Gamma_k^-,\mathfrak t}(X).
 \end{aligned}
\end{equation}
If at least one of $m_{\Gamma_k^+,\mathfrak t}(X)$ and
$m_{\Gamma_k^-,\mathfrak t}(X)$ has degree $B_{\Gamma_k^+,\mathfrak t}$, then
$m_{\Gamma_{k,\varepsilon}^0,\mathfrak t}(X)$ has degree
$B_{\Gamma_{k,\varepsilon}^0,\mathfrak t}$, and
\begin{equation}\label{eq:even-transfer}
 m_{\Gamma_{k,\varepsilon}^0,\mathfrak t}(q)
 =\frac q2\bigl(m_{\Gamma_k^+,\mathfrak t}(q)
                  +m_{\Gamma_k^-,\mathfrak t}(q)\bigr)
                  +O(q^{B_{\Gamma_k^+,\mathfrak t}}).
\end{equation}
\end{lemma}

\begin{proof}
\leavevmode Recall from Section~\ref{subsec:sequence-leading-count} that
$\sigma_i(\Gamma)=\prod_{j\le i}a_j$. By \eqref{eq:mi-Ri} and the definition of
the three indices,
\[
\begin{aligned}
 R_i(\Gamma_{k,\varepsilon}^0,\mathfrak t)
 &=R_i(\Gamma_k^\pm,\mathfrak t)+\mathbf 1_{\{i=k\}}\qquad(1\le i<r),\\
 \sigma_i(\Gamma_{k,\varepsilon}^0)&=\sigma_i(\Gamma_k^+)=\sigma_i(\Gamma_k^-)
 \qquad(1\le i\le r,\ i\ne k),\\
 \sigma_k(\Gamma_{k,\varepsilon}^0)&=\varepsilon\sigma_k(\Gamma_k^+),\qquad
 \sigma_k(\Gamma_k^-)=-\sigma_k(\Gamma_k^+).
\end{aligned}
\]
If $R_k(\Gamma,\mathfrak t)>0$, then by Lemma~\ref{lem:prefix-criterion}, we have
\[
 \mathscr S_{\Gamma_{k,\varepsilon}^0}^\pi
 =\mathscr S_{\Gamma_k^+}^\pi=\mathscr S_{\Gamma_k^-}^\pi,
\]
where the sets $\mathscr S_\Theta^\pi$ are defined in \eqref{eq:S-maxnum-Gamma}.
If $R_k(\Gamma,\mathfrak t)=0$, then by the same lemma, we have
\begin{equation*}
 \begin{gathered}
 \mathscr S_{\Gamma_{k,\varepsilon}^0}^\pi
 =\mathscr S_{\Gamma_k^+}^\pi\sqcup\mathscr S_{\Gamma_k^-}^\pi,\\
 |\mathscr S_{\Gamma_k^+}^\pi|=|\mathscr S_{\Gamma_k^-}^\pi|
 =\tfrac12|\mathscr S_{\Gamma_{k,\varepsilon}^0}^\pi|.
 \end{gathered}
\end{equation*}
In either case, by \eqref{eq:S-sign-parametrization} and
Proposition~\ref{cor:cGamma-lc-explicit}, we obtain \eqref{eq:two-step-budget-average}.

If $\deg m_{\Gamma_k^\xi,\mathfrak t}(X)=B_{\Gamma_k^+,\mathfrak t}$ for
some $\xi\in\{+,-\}$, then by \eqref{eq:two-step-budget-average} and
Proposition~\ref{cor:cGamma-lc-explicit},
$[X^{B_{\Gamma_{k,\varepsilon}^0,\mathfrak t}}]
m_{\Gamma_{k,\varepsilon}^0,\mathfrak t}(X)>0$. By
Proposition~\ref{prop:degree-budget}, $\deg m_{\Gamma_{k,\varepsilon}^0,\mathfrak t}(X)
=B_{\Gamma_{k,\varepsilon}^0,\mathfrak t}$. Finally, by \eqref{eq:two-step-budget-average} and
\eqref{eq:budget-mutation}, the polynomial
\[
 m_{\Gamma_{k,\varepsilon}^0,\mathfrak t}(X)
 -\frac X2\bigl(m_{\Gamma_k^+,\mathfrak t}(X)+m_{\Gamma_k^-,\mathfrak t}(X)\bigr)
\]
has degree at most $B_{\Gamma_k^+,\mathfrak t}$. This proves \eqref{eq:even-transfer}.
\end{proof}

For descent sequence indices $\Gamma'\le_{\rm seq}\Gamma$, write their numerical sequences as $(2\ell'_i)_i$ and $(2\ell_i)_i$. After padding the shorter sequence by zeros, define
\[
s_{\rm seq}(\Gamma,\Gamma'):=\sum_{i\ge1}i(\ell'_i-\ell_i).
\]
\begin{lemma}\label{lem:sseq-length}
Let $\Gamma'\le_{\rm seq}\Gamma$. Then every chain of elementary downward sequence mutations
from $\Gamma$ to $\Gamma'$ has length $s_{\rm seq}(\Gamma,\Gamma')$.
\end{lemma}
\begin{proof}
Every elementary downward mutation increases $\sum_i i\ell_i$ by one.
\end{proof}

\begin{corollary}\label{cor:cGamma-degree}
If $\Gamma\le_{\rm seq}\Gamma^{\max}$, where $\Gamma^{\max}$ is maximal with respect to
$\le_{\rm seq}$ among the indices realized by $\pi\in\Irr_{\mathfrak t}(\Sp_{2N}(\Fq))$ \textup{(}by
\eqref{eq:max-index-num}, its numerical sequence depends only on $\mathfrak t$\textup{)}, then
\begin{equation}\label{eq:cGamma-distance}
 \deg m_{\Gamma,\mathfrak t}(X)
 =B_{\Gamma,\mathfrak t}=s_{\rm seq}(\Gamma^{\max},\Gamma).
\end{equation}
\end{corollary}
\begin{proof}
\leavevmode By
Proposition~\ref{cor:cGamma-lc-explicit},
$\deg m_{\Gamma,\mathfrak t}(X)=B_{\Gamma,\mathfrak t}$, since the realized
$\Gamma^{\max}$ lies in $\mathscr S_\Gamma^\pi$.
Then, by \eqref{eq:mi-Ri} and $\sum_j\ell_j=\sum_j\ell_j^{\max}=N$, we have
\[
\begin{aligned}
 B_{\Gamma,\mathfrak t}
 &=\sum_{i<r}\sum_{j\le i}(\ell_j^{\max}-\ell_j)\\
 &=\sum_{j=1}^r j(\ell_j-\ell_j^{\max})
  =s_{\rm seq}(\Gamma^{\max},\Gamma).
\end{aligned}
\]
\end{proof}
{}

\subsection{Sequence mutations and stable orbit codimensions}\label{sec:geometry}
Fix $\pi\in\Irr_{\mathfrak t}(\Sp_{2N}(\Fq))$, and let
$\Gamma^{\max}$ be a maximal descent sequence index of $\pi$. For
$\Gamma\le_{\rm seq}\Gamma^{\max}$, the degree of
$m_{\Gamma,\mathfrak t}(X)$ is the length of any chain of elementary
sequence mutations from $\Gamma^{\max}$ to $\Gamma$, by
\eqref{eq:cGamma-distance} and Lemma~\ref{lem:sseq-length}. We now relate
this length to stable orbit dimensions. For good indices
$\Gamma'\le_{\rm seq}\Gamma$, we will prove
\[
 s_{\rm seq}(\Gamma,\Gamma')
 =\tfrac12\bigl(\dim\CO(\Gamma)^{\rm st}-\dim\CO(\Gamma')^{\rm st}\bigr),
\]
where $\CO(\Gamma)$ is the rational nilpotent orbit attached to $\Gamma$
in Section~\ref{subsec:good-indexes}.

Write
\[
 \Gamma^{\rm num}=(2\ell_1,\ldots,2\ell_r),
 \qquad
 (\Gamma')^{\rm num}=(2\ell'_1,\ldots,2\ell'_{r'}),
\]
append zeros to the shorter list, and let $\nu_\Gamma,\nu_{\Gamma'}$ be the partitions of the
corresponding stable orbits.

For a  partition $\nu$ of $2N$, choose $X_\nu\in(\Ost{\nu})^F$,
where $\Ost{\nu}$ is the stable nilpotent orbit with partition $\nu$
(Section~\ref{subsec:rational-orbit-order}).
By \cite[Corollary~6.1.4]{CM},
\begin{equation}\label{eq:sp-centralizer}
\begin{aligned}
 \dim C_{\bSp_{2N}}(X_\nu)
 &=\frac{1}{2}\left(\sum_j(\nu_j^t)^2+\#\{j:\nu_j\text{ is odd}\}\right),\\
 \sum_j(\nu_j^t)^2&=\sum_i(2i-1)\nu_i.
\end{aligned}
\end{equation}

\begin{lemma}\label{lem:good-index-dimension}
For a good index $\Gamma$ with $\Gamma^{\rm num}=(2\ell_1,\ldots,2\ell_r)$
and $\sum_i\ell_i=N$,
\begin{equation*}
 \dim\CO(\Gamma)^{\rm st}=\dim\bSp_{2N}-\sum_{i\ge1}(2i-1)\ell_i.
\end{equation*}
\end{lemma}
\begin{proof}
By \eqref{eq:sp-centralizer},
\[
 2\dim C_{\bSp_{2N}}(X_{\nu_\Gamma})
 =\sum_{j\ge1}(2j-1)(\nu_\Gamma)_j
  +\#\{j:(\nu_\Gamma)_j\text{ is odd}\}.
\]
For each $1\le i<r$ with $\ell_i=\ell_{i+1}-1$, the entries
$2\ell_i,2\ell_{i+1}$ of $\Gamma^{\rm num}$ give the odd parts in positions
$i,i+1$ of $\nu_\Gamma$, namely
$(\nu_\Gamma)_i=(\nu_\Gamma)_{i+1}=2\ell_i+1=2\ell_{i+1}-1$
(Section~\ref{subsec:good-indexes}). In the displayed formula, the terms with
$j=i,i+1$ contribute $(2i-1)(\nu_\Gamma)_i+(2i+1)(\nu_\Gamma)_{i+1}$ to the sum, while
$(\nu_\Gamma)_i,(\nu_\Gamma)_{i+1}$ contribute $2$ to
$\#\{j:(\nu_\Gamma)_j\text{ is odd}\}$. Thus
\[
 (2i-1)(2\ell_i+1)+(2i+1)(2\ell_i+1)+2
 =(2i-1)(2\ell_i)+(2i+1)(2\ell_{i+1}).
\]
Repeating this calculation for all odd pairs of $\nu_\Gamma$ gives
\begin{equation*}
 \dim C_{\bSp_{2N}}(X_{\nu_\Gamma})
 =\frac12\sum_{j\ge1}(2j-1)(2\ell_j)
 =\sum_{j\ge1}(2j-1)\ell_j.
\end{equation*}
This proves the lemma, since
$\dim\CO(\Gamma)^{\rm st}=\dim\bSp_{2N}-\dim C_{\bSp_{2N}}(X_{\nu_\Gamma})$.
\end{proof}

\begin{proposition}\label{prop:endpoint-codim}
For good indices $\Gamma'\le_{\rm seq}\Gamma$,
\begin{equation*}
s_{\rm seq}(\Gamma,\Gamma')
 =\frac{1}{2}\left(
 \dim\CO(\Gamma)^{\rm st}
 -\dim\CO(\Gamma')^{\rm st}
 \right).
\end{equation*}
\end{proposition}
\begin{proof}
Every sequence mutation preserves $\sum_i\ell_i$. Thus, by
Lemma~\ref{lem:good-index-dimension},
\[
 \dim\CO(\Gamma)^{\rm st}-\dim\CO(\Gamma')^{\rm st}
 =\sum_{i\ge1}(2i-1)(\ell'_i-\ell_i)
 =2\sum_{i\ge1}i(\ell'_i-\ell_i)
 =2s_{\rm seq}(\Gamma,\Gamma').
\]
\end{proof}

Let $\Gamma$ be a good index realized by $\pi$. As noted after \eqref{eq:max-index-num}, there is an index $\Gamma^{\max}$, maximal
among those realized by $\pi$, with $\Gamma\le_{\rm seq}\Gamma^{\max}$.
Combining Corollary~\ref{cor:cGamma-degree}, the prefix formula
\eqref{eq:mi-Ri}, and Proposition~\ref{prop:endpoint-codim}, we obtain
\begin{equation}\label{eq:prefix-budget-codim}
\begin{aligned}
 \deg m_{\Gamma,\mathfrak t}(X)
 &=B_{\Gamma,\mathfrak t}=\sum_{i<r}R_i\\
 &=s_{\rm seq}(\Gamma^{\max},\Gamma)\\
 &=\tfrac12\bigl(\dim\CO_\pi^{\rm st}-\dim\CO(\Gamma)^{\rm st}\bigr).
\end{aligned}
\end{equation}
Here $\CO(\Gamma^{\max})^{\rm st}=\CO_\pi^{\rm st}$ by
\cite[Theorem~7.1]{wang}.

\begin{remark}
The number $s_{\rm seq}(\Gamma,\Gamma')$ counts the sequence mutations.
Even when the rational orbits $\CO(\Gamma)$ and $\CO(\Gamma')$ differ by one orbit mutation,
this number can exceed one.
\end{remark}

\begin{example}\label{ex:4222}
In the stable orbit order for $\bSp_{10}$, consider the interval from
$\Ost{2^5}$ to $\Ost{4,2,2,2}$.
Vertex subscripts give dimensions; edge labels give codimensions.
\[
\begin{tikzpicture}[
x=2.8cm,y=1.9cm,
every node/.style={font=\small,inner sep=1pt},
cover/.style={draw,line width=0.45pt},
descent/.style={draw,line width=1.35pt}
]
\node (T) at (0,0) {$\Ost{4,2,2,2}_{\,38}$};
\node (L) at (-1,-1) {$\Ost{4,2,2,1,1}_{\,36}$};
\node (R) at (1,-1) {$\Ost{3,3,2,2}_{\,36}$};
\node (M) at (0,-2) {$\Ost{3,3,2,1,1}_{\,34}$};
\node (B) at (0,-3) {$\Ost{2^5}_{\,30}$};
\draw[cover] (T)--node[pos=.48,above left,font=\scriptsize] {$2$} (L);
\draw[descent] (T)--node[pos=.48,above right,font=\scriptsize] {$2$} (R);
\draw[cover] (L)--node[pos=.48,below left,font=\scriptsize] {$2$} (M);
\draw[descent] (R)--node[pos=.48,below right,font=\scriptsize] {$2$} (M);
\draw[descent] (M)--node[pos=.5,right,font=\scriptsize] {$4$} (B);
\end{tikzpicture}
\]
The bold edges form a possible downward path in the stable orbit order.
The good indices $\Gamma$ and $\Gamma'$ representing the top
and bottom orbits have $(\ell_i)_i=(2,1,1,1,0)$ and
$(\ell'_i)_i=(1,1,1,1,1)$, respectively.
Four elementary sequence mutations move one unit from the first to the
fifth position. The stable orbit dimensions are $38$ and $30$, and
\[
 s_{\rm seq}(\Gamma,\Gamma')=4=\tfrac12(38-30).
\]
The signed descent sequence indices along the branch through
$\Ost{3,3,2,2}$ are as follows:
\[
\begin{aligned}
&\llbracket(4,+),(2,\delta),(2,+),(2,\delta),(0,+)\rrbracket\\
&\quad\rightsquigarrow
 \llbracket(2,\delta),(4,+),(2,+),(2,\delta),(0,+)\rrbracket\\
&\quad\rightsquigarrow
 \llbracket(2,\delta),(4,+),(2,+),(0,\delta),(2,+)\rrbracket\\
&\quad\rightsquigarrow
 \llbracket(2,\delta),(2,+),(4,+),(0,\delta),(2,+)\rrbracket\\
&\quad\rightsquigarrow
 \llbracket(2,\delta),(2,+),(2,\delta),(2,+),(2,+)\rrbracket .
\end{aligned}
\]
The fourth index is not good. The other four indices give, in order,
the stable orbits at the vertices of the bold path:
\[
 \Ost{4,2,2,2},\quad \Ost{3,3,2,2},\quad
 \Ost{3,3,2,1,1},\quad \Ost{2^5}.
\]

Along this chain, the values of
$\dim\bSp_{10}-\sum_{i=1}^5(2i-1)\ell_i=55-\sum_{i=1}^5(2i-1)\ell_i$ are
\[
 38\longrightarrow36\longrightarrow34\longrightarrow32\longrightarrow30.
\]
At each good index this is the dimension of the underlying stable orbit
(Lemma~\ref{lem:good-index-dimension}).
The index giving $32$ is not good and has no associated rational orbit.
Thus the covering relation $\Ost{3,3,2,1,1}>\Ost{2^5}$, of codimension
$4$, requires two sequence mutations.

There are fourteen rational orbits lying over the stable orbits in
this interval. Their Hasse diagram is shown below, with equal even
rows combined as in \eqref{eq:equal-even-combination} and the following
chain highlighted:
\[
 [(4,+),(2^3,\delta)]
 >_{\rm rat}[3^2,(2^2,+)]
 >_{\rm rat}[3^2,(2,+),1^2]
 >_{\rm rat}[(2^5,+)].
\]

\begin{center}
\resizebox{0.94\textwidth}{!}{%
\begin{tikzpicture}[
x=0.75cm,y=1.05cm,
every node/.style={inner sep=1.5pt,outer sep=1pt,font=\small,fill=white},
cover/.style={draw,line width=0.45pt},
descent/.style={draw,line width=1.35pt}
]
\node (App) at (-6,6) {$\CO_{[(4,+),(2^3,\delta)]}$};
\node (Apm) at (-2,6) {$\CO_{[(4,+),(2^3,-\delta)]}$};
\node (Amp) at ( 2,6) {$\CO_{[(4,-),(2^3,-\delta)]}$};
\node (Amm) at ( 6,6) {$\CO_{[(4,-),(2^3,\delta)]}$};

\node (Bpp) at (-9,4) {$\CO_{[(4,+),(2^2,+),1^2]}$};
\node (Bpm) at (-5.4,4) {$\CO_{[(4,+),(2^2,-),1^2]}$};
\node (Bmp) at ( 5.4,4) {$\CO_{[(4,-),(2^2,-),1^2]}$};
\node (Bmm) at ( 9,4) {$\CO_{[(4,-),(2^2,+),1^2]}$};
\node (Cp) at (-1.6,4) {$\CO_{[3^2,(2^2,+)]}$};
\node (Cm) at ( 1.6,4) {$\CO_{[3^2,(2^2,-)]}$};

\node (Dp) at (-1.8,2) {$\CO_{[3^2,(2,+),1^2]}$};
\node (Dm) at ( 1.8,2) {$\CO_{[3^2,(2,-),1^2]}$};

\node (Ep) at (-1.8,0) {$\CO_{[(2^5,+)]}$};
\node (Em) at ( 1.8,0) {$\CO_{[(2^5,-)]}$};

\begin{scope}[on background layer]
\draw[cover] (Bpp) -- (App);
\draw[cover] (Bpp) -- (Apm);
\draw[cover] (Bpm) -- (App);
\draw[cover] (Bpm) -- (Apm);
\draw[cover] (Bmp) -- (Amp);
\draw[cover] (Bmp) -- (Amm);
\draw[cover] (Bmm) -- (Amp);
\draw[cover] (Bmm) -- (Amm);

\draw[descent] (Cp) to[bend left=11] (App);
\draw[cover] (Cp) to[bend right=9] (Amp);
\draw[cover] (Cm) to[bend left=9] (Apm);
\draw[cover] (Cm) to[bend right=11] (Amm);

\draw[cover] (Dp) to[bend left=18] (Bpp);
\draw[cover] (Dp) to[bend right=13] (Bmp);
\draw[cover] (Dm) to[bend left=13] (Bpm);
\draw[cover] (Dm) to[bend right=18] (Bmm);

\draw[descent] (Dp) -- (Cp);
\draw[cover] (Dp) -- (Cm);
\draw[cover] (Dm) -- (Cp);
\draw[cover] (Dm) -- (Cm);

\draw[descent] (Ep) -- (Dp);
\draw[cover] (Em) -- (Dm);
\end{scope}
\end{tikzpicture}%
}
\end{center}
\end{example}

\section{Whittaker multiplicities for symplectic groups}\label{sec:orbit-multiplicities}
Fix a $q$-system $\boldsymbol\pi=\{\pi_q\}$ of type $\mathfrak t$.
For a rational nilpotent orbit $\CO$ in $\mathfrak{sp}_{2N}(\Fq)$, define
\begin{equation*}
 m_{\CO,\boldsymbol\pi}(q):=\dim\Wh_\CO(\pi_q).
\end{equation*}
\subsection{Composition formulas and descent paths}
We compute $m_\CO(\pi)=\dim\Wh_\CO(\pi)$ by removing the rows of
$\CO$ from largest to smallest. If the largest part is even, say $2\ell$, write
$\CO=\CO_{2\ell,a}*\CO'$, where $\CO'$ is obtained by removing the signed row
$(2\ell,a)$. By \eqref{eq:even-composition-law},
\[
 \dim\Wh_\CO(\pi)
 =\dim\Wh_{\CO'}\!\left(\cD^{\FJ}_{2\ell,a}(\pi)\right).
\]
If the largest part is odd, say $2\ell+1$, write
$\CO=\CO_{(2\ell+1)^2}*\CO'$, where $\CO'$ is obtained by removing the pair
$(2\ell+1)^2$. By \eqref{eq:odd-general-input}, for every $a\in\{\pm1\}$, we have
\[
 \begin{aligned}
 2\dim\Wh_\CO(\pi)
 ={}&\dim\Wh_{\CO'}\!\left(
       \cD^{\FJ}_{2\ell+2,\delta a}
       \cD^{\FJ}_{2\ell,a}(\pi)\right)\\
 &+\frac{q+1}{2}\dim\Wh_{\CO'}\!\left(
       \cD^{\FJ}_{2\ell,\delta a}
       \cD^{\FJ}_{2\ell+2,a}(\pi)\right)\\
 &+\frac{q-1}{2}\dim\Wh_{\CO'}\!\left(
       \cD^{\FJ}_{2\ell,-\delta a}
       \cD^{\FJ}_{2\ell+2,-a}(\pi)\right).
 \end{aligned}
\]

Write
\begin{equation*}
 \CO=[(\lambda_1^{n_1},a_1),\ldots,(\lambda_t^{n_t},a_t)]
\end{equation*}
with underlying partition $\mu=(\mu_1\ge\cdots\ge\mu_L)$. For each
even block $(\lambda_i^{n_i},a_i)$, choose signs
$\beta_{i,1},\ldots,\beta_{i,n_i}\in\{\pm1\}$ such that
\begin{equation}\label{eq:even-block-singleton-signs}
 a_i=\delta^{n_i(n_i-1)/2}\prod_{k=1}^{n_i}\beta_{i,k}.
\end{equation}
By \eqref{eq:equal-even-combination}, this splits the block into
the signed rows $(\lambda_i,\beta_{i,k})$.
For the odd rows, set
\[
 \begin{gathered}
 s_\CO:=\tfrac12\#\{h:\mu_h\text{ odd}\},\qquad
 \mu_{y_j}=\mu_{y_j+1}=2r_j+1,\\
 \#\{h<y_j:\mu_h\text{ odd}\}=2j-2
 \qquad(1\le j\le s_\CO).
 \end{gathered}
\]
Choose $b_j\in\{\pm1\}$ for $1\le j\le s_\CO$, and recall
$\delta=\varepsilon_q(-1)$. Let
\begin{equation}\label{eq:odd-channel-table}
 \begin{aligned}
 \mathbf I_j^0&:=\bigl((2r_j,b_j),(2r_j+2,\delta b_j)\bigr),
 &w_X(0)&:=1,\\
 \mathbf I_j^+&:=\bigl((2r_j+2,b_j),(2r_j,\delta b_j)\bigr),
 &w_X(+)&:=\frac{X+1}{2},\\
 \mathbf I_j^-&:=\bigl((2r_j+2,-b_j),(2r_j,-\delta b_j)\bigr),
 &w_X(-)&:=\frac{X-1}{2}.
\end{aligned}
\end{equation}
For $\mathbf c=(c_1,\ldots,c_{s_\CO})\in
\mathscr R(\CO):=\{0,+,-\}^{s_\CO}$, define the descent sequence index
$\Gamma_{\CO,\mathbf c}$  by the
following replacements, preserving the order of the parts of $\mu$:
\[
 \begin{aligned}
 (\lambda_i^{n_i},a_i)&\longmapsto
 \bigl((\lambda_i,\beta_{i,1}),\ldots,(\lambda_i,\beta_{i,n_i})\bigr)
 &&(\lambda_i\text{ even}),\\
 (\mu_{y_j},\mu_{y_j+1})&\longmapsto\mathbf I_j^{c_j}
 &&(1\le j\le s_\CO).
 \end{aligned}
\]
Put
\begin{equation}\label{eq:resolution-state-weight}
 w_X(\mathbf c):=\prod_{j=1}^{s_\CO}w_X(c_j).
\end{equation}
Define
\begin{equation}\label{eq:uniform-orbit-expansion-C}
 m_{\CO,\mathfrak t}(X)
 :=2^{-s_\CO}\sum_{\mathbf c\in\mathscr R(\CO)}
 w_X(\mathbf c)\,m_{\Gamma_{\CO,\mathbf c},\mathfrak t}(X).
\end{equation}

\begin{proposition}
\label{prop:uniform-orbit-expansion}
For every finite-dimensional representation $\pi$ of $\Sp_{2N}(\Fq)$,
\begin{equation}\label{eq:uniform-orbit-expansion}
 2^{s_{\CO}}\dim\Wh_{\CO}(\pi)
 =\sum_{\mathbf c\in\mathscr R(\CO)}
{w_q(\mathbf c)}\,m_{\Gamma_{\CO,\mathbf c}}(\pi).
\end{equation}
For every $q$-system $\boldsymbol\pi$ of type $\mathfrak t$ and every
sufficiently large $q$,
\begin{equation*}
m_{\CO,\boldsymbol\pi}(q)=m_{\CO,\mathfrak t}(q).
\end{equation*}
\end{proposition}

\begin{proof}
\begingroup
For $\mathbf c\in\mathscr R(\CO)$, write
\[
 \Gamma_{\CO,\mathbf c}
 =\llbracket(2\ell_1(\mathbf c),\alpha_1(\mathbf c)),\ldots,
 (2\ell_L(\mathbf c),\alpha_L(\mathbf c))\rrbracket.
\]
By \eqref{eq:even-composition-law} and \eqref{eq:odd-general-input}, applied
successively with the choices in \eqref{eq:odd-channel-table}, we have
\[
\begin{aligned}
 2^{s_\CO}\dim\Wh_\CO(\pi)
 &=\sum_{\mathbf c\in\mathscr R(\CO)}
   \left(\prod_{j=1}^{s_\CO}w_q(c_j)\right)
   \dim\!\left(\cD^{\FJ}_{2\ell_L(\mathbf c),\alpha_L(\mathbf c)}\cdots
        \cD^{\FJ}_{2\ell_1(\mathbf c),\alpha_1(\mathbf c)}(\pi)\right)\\
 &=\sum_{\mathbf c\in\mathscr R(\CO)}
   w_q(\mathbf c)\,m_{\Gamma_{\CO,\mathbf c}}(\pi),
\end{aligned}
\]
by \eqref{eq:cGamma-iterated} and \eqref{eq:resolution-state-weight}.
To see that $m_{\CO,\mathfrak t}(X)$ in \eqref{eq:uniform-orbit-expansion-C} is
independent of $\delta=\varepsilon_q(-1)$, let $\CO'$ be the orbit for $q'$ identified
with $\CO$ under Section~\ref{subsec:rational-orbit-order}. By the compatibility of
$\Gamma\mapsto\CO(\Gamma)$ with the cross-$q$ identification of
Section~\ref{subsec:sequence-coefficients}, multiplying the $m$th sign of
$\Gamma_{\CO,\mathbf0}$ by $(\varepsilon_q(-1)\varepsilon_{q'}(-1))^{m-1}$ yields a
valid splitting $\Gamma_{\CO',\mathbf0}$ of $\CO'$ for $q'$ (matching
\eqref{eq:even-block-singleton-signs} on each even block). Moreover, on the $j$th odd
pair (at consecutive positions $y_j,y_j+1$), multiplying the two signs by
$(\varepsilon_q(-1)\varepsilon_{q'}(-1))^{y_j-1}$ and
$(\varepsilon_q(-1)\varepsilon_{q'}(-1))^{y_j}$ changes their relative sign by
$\varepsilon_q(-1)\varepsilon_{q'}(-1)$, which maps each channel $\mathbf I_j^{c_j}$ in
\eqref{eq:odd-channel-table} for $q$ to the corresponding channel for $q'$ with the
same $c_j$ and the same weight $w_X(c_j)$. Hence $\CO$ and $\CO'$ yield the same
identified indices $\Gamma_{\CO,\mathbf c}$ and weights $w_X(\mathbf c)$ in
\eqref{eq:uniform-orbit-expansion-C}; by Proposition~\ref{prop:exact-poly},
$m_{\CO,\mathfrak t}(X)$ is independent of $\delta$.
By \eqref{eq:uniform-orbit-expansion}, \eqref{eq:type-polynomial-value} and
\eqref{eq:uniform-orbit-expansion-C},
\begin{align*}
 m_{\CO,\boldsymbol\pi}(q)
 &=2^{-s_\CO}\sum_{\mathbf c\in\mathscr R(\CO)}
   w_q(\mathbf c)\,m_{\Gamma_{\CO,\mathbf c}}(\pi_q)\\
 &=2^{-s_\CO}\sum_{\mathbf c\in\mathscr R(\CO)}
   w_q(\mathbf c)\,m_{\Gamma_{\CO,\mathbf c},\mathfrak t}(q)
  =m_{\CO,\mathfrak t}(q).\qedhere
\end{align*}
\endgroup
\end{proof}

For $\mathbf c\in\mathscr R(\CO)$, set
$u(\mathbf c):=\#\{j:c_j\ne0\}$. Both
$\mathbf I_j^+\rightsquigarrow\mathbf I_j^0$ and
$\mathbf I_j^-\rightsquigarrow\mathbf I_j^0$ are elementary downward sequence
mutations \eqref{eq:elementary-seq-mutation}. Hence $\Gamma_{\CO,\mathbf 0}$ is
obtained from $\Gamma_{\CO,\mathbf c}$ by $u(\mathbf c)$ elementary downward mutations, one at the
positions $y_j,y_j+1$ for each $j$ with $c_j\ne0$. By \eqref{eq:budget-mutation}, we have
\begin{equation}\label{eq:channel-budget-drop}
 B_{\Gamma_{\CO,\mathbf c},\mathfrak t}
 =B_{\Gamma_{\CO,\mathbf 0},\mathfrak t}-u(\mathbf c).
\end{equation}

\begin{example}\label{ex:orbit-expansion}
Take the rational nilpotent orbit
$\CO=[(4,+),3^2,(2,+)]$ under $\Sp_{12}(\Fq)$.
Then $s_\CO=1$, and choosing $b_1=+$ gives
\[
\begin{aligned}
 \Gamma_{\CO,(0)}
 &=\llbracket(4,+),(2,+),(4,\delta),(2,+)\rrbracket,\\
 \Gamma_{\CO,(+)}
 &=\llbracket(4,+),(4,+),(2,\delta),(2,+)\rrbracket,\\
 \Gamma_{\CO,(-)}
 &=\llbracket(4,+),(4,-),(2,-\delta),(2,+)\rrbracket.
\end{aligned}
\]
The first and last entries come from the even rows $(4,+)$ and $(2,+)$;
the two middle entries replace the odd pair $3^2$.
\end{example}

\subsection{Degrees of Whittaker multiplicities for rational orbits}

\begin{theorem}\label{prop:orbit-degree}
Let $\pi\in\Irr_{\mathfrak t}(\Sp_{2N}(\Fq))$ and $\CO\in\ScO(\pi)$. Then
\begin{equation}\label{eq:orbit-degree}
 \deg m_{\CO,\mathfrak t}(X)
 =\frac{1}{2}\left(\dim\CO_\pi^{\rm st}-\dim\CO^{\rm st}\right).
\end{equation}
\end{theorem}
\begin{proof}
\begingroup
By \cite[Corollary~5.2]{wang}, choose a good index $\Gamma\in\mathfrak G(\CO)$
realized by $\pi$, and choose the signs in \eqref{eq:even-block-singleton-signs}
and \eqref{eq:odd-channel-table} so that $\Gamma_{\CO,\mathbf0}=\Gamma$.
By \eqref{eq:prefix-budget-codim}, we have
\[
 \deg m_{\Gamma,\mathfrak t}(X)=B_{\Gamma,\mathfrak t}.
\]
By \eqref{eq:odd-channel-table}, \eqref{eq:resolution-state-weight},
Proposition~\ref{prop:degree-budget} and \eqref{eq:channel-budget-drop}, for
$\mathbf c\in\mathscr R(\CO)$ with
$m_{\Gamma_{\CO,\mathbf c},\mathfrak t}(X)\ne0$, we have
\[
\begin{aligned}
 \deg\!\left(w_X(\mathbf c)\,m_{\Gamma_{\CO,\mathbf c},\mathfrak t}(X)\right)
 &\le u(\mathbf c)+B_{\Gamma_{\CO,\mathbf c},\mathfrak t}\\
 &=B_{\Gamma,\mathfrak t}.
\end{aligned}
\]
Every nonzero $m_{\Gamma_{\CO,\mathbf c},\mathfrak t}(X)$ has positive leading
coefficient \textup{(}see after \eqref{eq:weighted-descent-paths}\textup{)}, and
$\LC(w_X(\mathbf c))=2^{-u(\mathbf c)}>0$. Since $w_X(\mathbf0)=1$, the
$\mathbf c=\mathbf0$ term attains the degree of the sum in
\eqref{eq:uniform-orbit-expansion-C}. Using
$\CO(\Gamma)=\CO$ in \eqref{eq:prefix-budget-codim}, we obtain
\begin{equation}\label{eq:orbit-degree-equals-B}
 \deg m_{\CO,\mathfrak t}(X)=B_{\Gamma,\mathfrak t}{}=\tfrac12\bigl(\dim\CO_\pi^{\rm st}-\dim\CO^{\rm st}\bigr).
\end{equation}
\endgroup
\end{proof}
In particular, $m_{\CO,\mathfrak t}(X)$ has degree zero at every $\CO\in\ScO(\pi)^{\max}$, that is,
$m_\CO(\pi)$ does not depend on $q$; compare \cite[Theorem~1.9]{wang}.

\subsection{Leading coefficients and maximal rational orbits}\label{subsec:rational-LC}
For $\pi\in\Irr_{\mathfrak t}(\Sp_{2N}(\Fq))$ and $\CO\in\ScO(\pi)$, the expansion \eqref{eq:uniform-orbit-expansion-C} expresses $m_{\CO,\mathfrak t}(X)$ in terms of the polynomials $m_{\Gamma_{\CO,\mathbf c},\mathfrak t}(X)$, and for each descent sequence index $\Gamma$ with $m_{\Gamma,\mathfrak t}(X)\ne0$, Proposition~\ref{cor:cGamma-lc-explicit} expresses the leading coefficient of $m_{\Gamma,\mathfrak t}(X)$ in terms of the first-descent paths from $\pi$ to $\mathbf{1}_{\Sp_0}$ with index $\Gamma'\ge_{\rm seq}\Gamma$. To compute the leading coefficient of $m_{\CO,\mathfrak t}(X)$, we reduce the leading-term sum to a single good index $\Gamma\in\mathfrak G(\CO)$ and group the resulting first-descent paths according to the maximal rational orbits $\CO(\Gamma')\in\ScO(\pi)^{\max}$ attached to their indices $\Gamma'$.

For the partition
$
 \nu=(d_1^{m_1},\ldots,d_u^{m_u}),
$
of a nilpotent orbit of a symplectic group, define
\[
M_v(\nu):=\sum_{a=1}^v m_a,\qquad M_0(\nu):=0.
\]
Each odd part $d_v$ has even multiplicity $m_v$. Define
\begin{align*}
J_{\rm end}(\nu)
 &:=\{M_v(\nu):d_v\text{ is even}\}
 \cup\left\{M_{v-1}(\nu)+2b:
 \substack{d_v\text{ is odd}\\1\le b\le m_v/2}\right\},\\
J_{\rm odd}(\nu)
 &:=\left\{M_{v-1}(\nu)+2b-1:
 \substack{d_v\text{ is odd}\\1\le b\le m_v/2}\right\}.
\end{align*}
Thus $J_{\rm end}(\nu)$ consists of the last row of each even block and
the second row of each odd pair, while $J_{\rm odd}(\nu)$ consists of
the first row of each odd pair.
For a rational nilpotent orbit $\CO$ and a stable nilpotent orbit ${\CO'}^{\rm st}$ with
$\CO^{\rm st}\le{\CO'}^{\rm st}$, let $\lambda=(\lambda_i)_i$ and $\lambda'=(\lambda_i')_i$ be their
partitions and define
\[
S_i(\lambda):=\sum_{j=1}^i\lambda_j,
 \qquad
S_i(\lambda'):=\sum_{j=1}^i\lambda'_j,
\]
with zeros appended when necessary. Define
\begin{equation}\label{eq:Delta-stable-partitions}
\begin{aligned}
 \Delta(\CO,{\CO'}^{\rm st})
 &:={\#}\{i\in J_{\rm end}(\lambda'):S_i(\lambda)<S_i(\lambda')\}\\
 &\quad-{\#}\{i\in J_{\rm odd}(\lambda'):S_i(\lambda)<S_i(\lambda')\}.
\end{aligned}
\end{equation}

For an irreducible representation $\pi$ and $\CO\in\ScO(\pi)$, define
\begin{equation*}
Q_\pi(\CO)
 :=\#\{\CO'\in\ScO(\pi)^{\max}:\CO\le_{\rm rat}\CO'\}.
\end{equation*}
Thus $Q_\pi(\CO)$ is the number of maximal rational wavefront orbits of $\pi$ lying above
$\CO$. Recall from Section~\ref{sec:intro} that, for $\pi\in\cE(G,s)$, $\pi^{\rm J}=\CL_s(\pi)$,
$\CO_{\pi^{\rm J}}^{\rm st}$ denotes the stable wavefront orbit of $\pi^{\rm J}$ in the Lie algebra
of $C_{\bG^*}(s)$ over $\overline{\Fq}$, $A(\CO_\pi^{\rm st})$ is the component group of $\CO_\pi^{\rm st}$, and $\bA(\CO_{\pi^{\rm J}}^{\rm st})$ denotes Lusztig's canonical quotient.

\begin{theorem}\label{prop:orbit-leading-coefficient}
Let $\pi\in\Irr_{\mathfrak t}(\Sp_{2N}(\Fq))$ and $\CO\in\ScO(\pi)$.
Then
\begin{equation}\label{eq:orbit-leading-intrinsic}
 \LC(m_{\CO,\mathfrak t}(X))
 =2^{-\Delta(\CO,\CO_\pi^{\rm st})}
 \sum_{\substack{\CO''\in\ScO(\pi)^{\max}\\\CO\le_{\rm rat}\CO''}}
 \dim\Wh_{\CO''}(\pi).
\end{equation}
Moreover, for $\pi\in\cE(G,s)$,
\begin{equation}\label{eq:orbit-leading-canonical}
 \LC(m_{\CO,\mathfrak t}(X))
 =2^{-\Delta(\CO,\CO_\pi^{\rm st})}\,
  \frac{Q_\pi(\CO)\,|A(\CO_\pi^{\rm st})|}{|\ScO(\pi)^{\max}|\,|\bA(\CO_{\pi^{\rm J}}^{\rm st})|}.
\end{equation}
\end{theorem}

Before proving Theorem~\ref{prop:orbit-leading-coefficient}, we rewrite $\Delta(\CO(\Gamma),\CO(\Gamma')^{\rm st})$ in terms of good descent sequence indices.

\begin{lemma}\label{lem:Delta-prefix}
Let $\Gamma\le_{\rm seq}\Gamma'$ be good descent sequence indices, extended by zero entries
to a common length $r$, with numerical sequences $(2\ell_i)_i$ and $(2\ell'_i)_i$. Put
$R_i(\Gamma,\Gamma'):=\sum_{j\le i}(\ell'_j-\ell_j)$, and let $\lambda$ and $\lambda'$ be the
partitions of $\CO(\Gamma)^{\rm st}$ and $\CO(\Gamma')^{\rm st}$. Set
\[
\begin{aligned}
k&:=\#\{i\in J_{\rm end}(\lambda'):R_i(\Gamma,\Gamma')>0\},
 \qquad s':=|J_{\rm odd}(\lambda')|,\\
z&:=\#\{1\le i<r:R_i(\Gamma,\Gamma')=0\text{ and }i\text{ begins an odd pair of }\CO(\Gamma)\}.
\end{aligned}
\]
Then $\Delta(\CO(\Gamma),\CO(\Gamma')^{\rm st})=k+z-s'$.
\end{lemma}
\begin{proof}
Write $R_i:=R_i(\Gamma,\Gamma')$. For a partition $\nu$ of a nilpotent orbit of a symplectic group, let $\epsilon_i(\nu)\in\{0,1\}$ be $1$ exactly when
$i\in J_{\rm odd}(\nu)$. By the definition of $\CO(\Gamma)$ in Section~\ref{subsec:good-indexes},
\[
S_i(\lambda)-\epsilon_i(\lambda)=2\sum_{j\le i}\ell_j,
 \qquad
S_i(\lambda')-\epsilon_i(\lambda')=2\sum_{j\le i}\ell'_j.
\]
Then
\begin{equation}\label{eq:prefix-R-relation}
S_i(\lambda')-S_i(\lambda)
 =2R_i+\epsilon_i(\lambda')-\epsilon_i(\lambda),
\end{equation}
whose left-hand side is nonnegative because $\lambda\le\lambda'$ by \eqref{eq:good-order-compatibility}.
Since $\epsilon_i(\lambda')$ is $0$ on $J_{\rm end}(\lambda')$ and $1$ on $J_{\rm odd}(\lambda')$,
for the indices $i$ occurring in \eqref{eq:Delta-stable-partitions} \eqref{eq:prefix-R-relation} becomes
\[
 S_i(\lambda')-S_i(\lambda)
 =\begin{cases}
  2R_i-\epsilon_i(\lambda),&i\in J_{\rm end}(\lambda'),\\
  2R_i+1-\epsilon_i(\lambda),&i\in J_{\rm odd}(\lambda').
 \end{cases}
\]
For $i\in J_{\rm end}(\lambda')$, the difference $S_i(\lambda')-S_i(\lambda)$ is positive precisely when $R_i>0$;
for $i\in J_{\rm odd}(\lambda')$, it vanishes precisely when $R_i=0$ and $\epsilon_i(\lambda)=1$.
Conversely, these two conditions force $\epsilon_i(\lambda')=1$. Thus exactly $z$ indices in
$J_{\rm odd}(\lambda')$ satisfy $S_i(\lambda')=S_i(\lambda)$; the endpoint $r$ begins no odd
pair. The two counts in \eqref{eq:Delta-stable-partitions} are consequently $k$ and $s'-z$,
respectively. This proves the lemma.
\end{proof}

\begin{proof}[Proof of Theorem~\ref{prop:orbit-leading-coefficient}]
\noindent\textit{Step 1: the contribution of odd pairs.}
Let $\lambda^{\max}$ be the partition of $\CO_\pi^{\rm st}$. As in the proof of Theorem~\ref{prop:orbit-degree},
choose a realized good index $\Gamma\in\mathfrak G(\CO)$ and signs in
\eqref{eq:even-block-singleton-signs} and \eqref{eq:odd-channel-table} so that
$\Gamma_{\CO,\mathbf 0}=\Gamma$; then $m_{\Gamma,\mathfrak t}(X)\ne0$. Extending
$\Gamma$ and all $\Gamma_{\CO,\mathbf c}$ ($\mathbf c\in\mathscr R(\CO)$) by zero entries to a common length $r$ with $r_r^{\rm fd}(\mathfrak t)=0$, we write
$R_i:=R_i(\Gamma,\mathfrak t)$ (as in \eqref{eq:sequence-budget-main}) and define
\[
\begin{aligned}
h&:=\#\{1\le i<r:R_i>0\},\\
z&:=\#\{1\le i<r:R_i=0\text{ and }i\text{ begins an odd pair of }\CO(\Gamma)\}.
\end{aligned}
\]
By \eqref{eq:channel-budget-drop} and Proposition~\ref{prop:degree-budget},
$\deg m_{\Gamma_{\CO,\mathbf c},\mathfrak t}(X)\le B_{\Gamma_{\CO,\mathbf c},\mathfrak t}=B_{\Gamma,\mathfrak t}-u(\mathbf c)$
for every $\mathbf c\in\mathscr R(\CO)$.
Recall from \eqref{eq:odd-channel-table} and \eqref{eq:resolution-state-weight} that $w_X(\mathbf c)=\prod_{j=1}^{s_\CO}w_X(c_j)$ with $w_X(0)=1$ and $w_X(\pm)=(X\pm1)/2$, so $\deg w_X(\mathbf c)=u(\mathbf c)$ and $\LC(w_X(\mathbf c))=2^{-u(\mathbf c)}$. It therefore follows from \eqref{eq:uniform-orbit-expansion-C} that
\begin{equation}\label{eq:odd-budget-weighted-sum}
 [X^{B_{\Gamma,\mathfrak t}}]m_{\CO,\mathfrak t}(X)
 =2^{-s_{\CO}}\sum_{\mathbf c\in\mathscr R(\CO)}
   2^{-u(\mathbf c)}[X^{B_{\Gamma_{\CO,\mathbf c},\mathfrak t}}]m_{\Gamma_{\CO,\mathbf c},\mathfrak t}(X).
\end{equation}

For $\mathbf c\in\mathscr R(\CO)$, $1\le j\le s_\CO$ and $e\in\{0,+,-\}$, let $\mathbf c^{(j,e)}$ be
obtained from $\mathbf c$ by replacing its $j$th entry by $e$:
\[
 \mathbf c^{(j,e)}=(c_1,\ldots,c_{j-1},e,c_{j+1},\ldots,c_{s_\CO}),
\]
and let $i=y_j$ be the starting position of the $j$th odd pair of $\CO$ (as in Section~\ref{sec:orbit-multiplicities}).
We claim that
\begin{equation}\label{eq:local-odd-average-top}
\begin{aligned}
&\frac12[X^{B_{\Gamma_{\CO,\mathbf c^{(j,+)},\mathfrak t}}}]m_{\Gamma_{\CO,\mathbf c^{(j,+)},\mathfrak t}}(X)
 +\frac12[X^{B_{\Gamma_{\CO,\mathbf c^{(j,-)},\mathfrak t}}}]m_{\Gamma_{\CO,\mathbf c^{(j,-)},\mathfrak t}}(X)\\
&\qquad=
\begin{cases}
 0,&R_i=0,\\
 [X^{B_{\Gamma_{\CO,\mathbf c^{(j,0)},\mathfrak t}}}]m_{\Gamma_{\CO,\mathbf c^{(j,0)},\mathfrak t}}(X),&R_i>0.
\end{cases}
\end{aligned}
\end{equation}
Since the odd pairs occupy disjoint positions, we have $R_i(\Gamma_{\CO,\mathbf c^{(j,0)}},\mathfrak t)=R_i$. Moreover, $\Gamma_{\CO,\mathbf c^{(j,0)}}$ is obtained from $\Gamma_{\CO,\mathbf c^{(j,\pm)}}$ by an elementary
downward mutation at positions $i,i+1$, so (see Section~\ref{subsec:local-transfer})
\[
 R_k(\Gamma_{\CO,\mathbf c^{(j,\pm)}},\mathfrak t)=R_k(\Gamma_{\CO,\mathbf c^{(j,0)}},\mathfrak t)
 \quad(k\ne i),\qquad
 R_i(\Gamma_{\CO,\mathbf c^{(j,\pm)}},\mathfrak t)=R_i-1.
\]
We consider three cases:
\begin{itemize}[leftmargin=1.8em]
 \item \emph{Case 1: $R_k(\Gamma_{\CO,\mathbf c^{(j,0)}},\mathfrak t)<0$ for some $1\le k<r$ with $k\ne i$.} Then all three terms in \eqref{eq:local-odd-average-top} vanish by Proposition~\ref{prop:stagewise}, since no descent path has an index with a negative $R_k$.
 \item \emph{Case 2: $R_k(\Gamma_{\CO,\mathbf c^{(j,0)}},\mathfrak t)\ge0$ for all $k\ne i$, and $R_i=0$.} Then $R_i(\Gamma_{\CO,\mathbf c^{(j,\pm)}},\mathfrak t)=-1<0$, so the two terms on the left-hand side of \eqref{eq:local-odd-average-top} vanish, again by Proposition~\ref{prop:stagewise}.
 \item \emph{Case 3: $R_k(\Gamma_{\CO,\mathbf c^{(j,0)}},\mathfrak t)\ge0$ for all $k\ne i$, and $R_i>0$.} Then all $R_m(\Gamma_{\CO,\mathbf c^{(j,+)}},\mathfrak t)$ ($1\le m<r$) are nonnegative. Applying the notation of Section~\ref{subsec:local-transfer} to $\Gamma_{\CO,\mathbf c^{(j,+)}}$ with $k=i$ and $\varepsilon=+1$ gives $\Gamma_i^+=\Gamma_{\CO,\mathbf c^{(j,+)}}$, $\Gamma_i^-=\Gamma_{\CO,\mathbf c^{(j,-)}}$ and $\Gamma_{i,+1}^0=\Gamma_{\CO,\mathbf c^{(j,0)}}$, so \eqref{eq:two-step-budget-average} yields \eqref{eq:local-odd-average-top}.
\end{itemize}

In the sum \eqref{eq:odd-budget-weighted-sum}, fix $j\in\{1,\ldots,s_\CO\}$ and the remaining entries of $\mathbf c$. Since $u(\mathbf c^{(j,\pm)})=u(\mathbf c^{(j,0)})+1$, summing the corresponding three terms of \eqref{eq:odd-budget-weighted-sum} over $c_j=e\in\{0,+,-\}$ and applying \eqref{eq:local-odd-average-top} gives
\[
\begin{aligned}
 &\sum_{e\in\{0,+,-\}}2^{-u(\mathbf c^{(j,e)})}
  [X^{B_{\Gamma_{\CO,\mathbf c^{(j,e)}},\mathfrak t}}]m_{\Gamma_{\CO,\mathbf c^{(j,e)}},\mathfrak t}(X)\\
 &\qquad=
 \begin{cases}
  2^{-u(\mathbf c^{(j,0)})}[X^{B_{\Gamma_{\CO,\mathbf c^{(j,0)}},\mathfrak t}}]m_{\Gamma_{\CO,\mathbf c^{(j,0)}},\mathfrak t}(X), & R_{y_j}=0,\\
  2\cdot 2^{-u(\mathbf c^{(j,0)})}[X^{B_{\Gamma_{\CO,\mathbf c^{(j,0)}},\mathfrak t}}]m_{\Gamma_{\CO,\mathbf c^{(j,0)}},\mathfrak t}(X), & R_{y_j}>0.
 \end{cases}
\end{aligned}
\]
Applying this identity successively for $j=1,\ldots,s_\CO$ reduces the sum over all $\mathbf c\in\mathscr R(\CO)$ in \eqref{eq:odd-budget-weighted-sum} to the single term $\mathbf c=\mathbf0=(0,\ldots,0)$: each of the $s_\CO-z$ odd pairs with $R_{y_j}>0$ contributes a factor of $2$, while each of the $z$ odd pairs with $R_{y_j}=0$ contributes $1$. Hence
\[
 [X^{B_{\Gamma,\mathfrak t}}]m_{\CO,\mathfrak t}(X)
 =2^{-s_\CO}2^{s_\CO-z}[X^{B_{\Gamma,\mathfrak t}}]m_{\Gamma,\mathfrak t}(X)
 =2^{-z}[X^{B_{\Gamma,\mathfrak t}}]m_{\Gamma,\mathfrak t}(X).
\]
Since both $m_{\CO,\mathfrak t}(X)$ and $m_{\Gamma,\mathfrak t}(X)$ have degree $B_{\Gamma,\mathfrak t}$ by
\eqref{eq:orbit-degree-equals-B},
\begin{equation}\label{eq:orbit-lc-vs-C0}
\begin{aligned}
 \LC(m_{\CO,\mathfrak t}(X))
 &=2^{-z}\LC(m_{\Gamma,\mathfrak t}(X))\\
 &=2^{-z-h}\sum_{\Gamma'\in\mathscr S_\Gamma^\pi}
                 m_{\Gamma'}(\pi),
\end{aligned}
\end{equation}
where the last equality is \eqref{eq:budget-chain-count} applied to $\Gamma$, together with
$m_{\Gamma',\mathfrak t}(X)=m_{\Gamma'}(\pi)$ (by Proposition~\ref{prop:exact-poly}, since $m_{\Gamma',\mathfrak t}(X)$ is constant for $\Gamma'\in\mathscr S_\Gamma^\pi$). In this way, the sum over $\mathbf c\in\mathscr R(\CO)$ in \eqref{eq:odd-budget-weighted-sum} becomes a sum over $\mathscr S_\Gamma^\pi$.

\medskip
\noindent\textit{Step 2: computing $\Delta$ and counting good indices.}
Since $\Gamma$ is realized by $\pi$, there exists a maximal descent sequence index $\Gamma^{\max}$ realized by $\pi$ with $\Gamma\le_{\rm seq}\Gamma^{\max}$ (see Section~\ref{subsec:type-transitions}). Then
$\CO(\Gamma^{\max})^{\rm st}=\CO_\pi^{\rm st}$ by \cite[Theorem~7.1]{wang}, and
$R_i=R_i(\Gamma,\mathfrak t)=\sum_{j\le i}(\ell_j^{\max}-\ell_j)=R_i(\Gamma,\Gamma^{\max})$ by \eqref{eq:mi-Ri}, where $(2\ell_i)_i$ and $(2\ell_i^{\max})_i$ are the numerical sequences of $\Gamma$ and $\Gamma^{\max}$ (extended by zero entries to length $r$). Hence $R_i$ and $z$ from Step~1 agree with $R_i(\Gamma,\Gamma^{\max})$ and $z$ in Lemma~\ref{lem:Delta-prefix}. Applying Lemma~\ref{lem:Delta-prefix} to $\Gamma\le_{\rm seq}\Gamma^{\max}$ with $k:=\#\{i\in J_{\rm end}(\lambda^{\max}):R_i>0\}$ and $s_{\max}:=|J_{\rm odd}(\lambda^{\max})|$, we obtain
\begin{equation}\label{eq:Delta-k-z}
 \Delta(\CO,\CO_\pi^{\rm st})=k+z-s_{\max}.
\end{equation}

Fix $\CO''\in\ScO(\pi)^{\max}$ with $\CO\le_{\rm rat}\CO''$. Write
\[
 \CO''=[(d_1^{m_1},a_1),\ldots,(d_u^{m_u},a_u)]
\]
and let $M_v=\sum_{w\le v}m_w$ and $M_0=0$.
By \eqref{eq:S-sign-parametrization}, the map $\Gamma'\mapsto(\sigma_i(\Gamma'))_{R_i>0}$ of prefix signs identifies $\mathscr S_\Gamma^\pi$ with $\{\pm1\}^{h}$.
Since every $\Gamma'\in\mathscr S_\Gamma^\pi$ has numerical sequence $(\Gamma^{\max})^{\rm num}$ (which determines the underlying partition $\lambda^{\max}$ of $\CO''$), an index $\Gamma'\in\mathscr S_\Gamma^\pi$ lies in $\mathfrak G(\CO'')$ if and only if its prefix signs $\sigma_i:=\sigma_i(\Gamma')$ obey the sign rule for good indices \eqref{eq:good-cases} together with \eqref{eq:equal-even-combination}. Explicitly,
\begin{equation}\label{eq:fixed-orbit-prefix-constraints}
 \begin{aligned}
  \sigma_{M_v}
   &=\delta^{m_v(m_v-1)/2}a_v
       \sigma_{M_{v-1}},
       &&d_v\text{ even},\\
  \sigma_{M_{v-1}+2b}
   &=\delta^b\sigma_{M_{v-1}},
       &&d_v\text{ odd},\quad1\le b\le m_v/2.
 \end{aligned}
\end{equation}
Starting from $\sigma_0=1$, the relations \eqref{eq:fixed-orbit-prefix-constraints} successively determine
the prefix signs $\sigma_i(\Gamma')$ for $i\in J_{\rm end}(\lambda^{\max})$, while leaving all remaining prefix signs $\sigma_i(\Gamma')$ ($i\notin J_{\rm end}(\lambda^{\max})$) unconstrained.
Moreover, by \eqref{eq:orbit-to-sequence-order}, there exists $\Gamma'\in\mathfrak G(\CO'')$ with $\Gamma\le_{\rm seq}\Gamma'$; since $\CO''$ has underlying stable orbit $\CO_\pi^{\rm st}$, any such $\Gamma'$  has numerical sequence $(\Gamma^{\max})^{\rm num}$ and hence lies in $\mathscr S_\Gamma^\pi$. Thus $\mathscr S_\Gamma^\pi\cap\mathfrak G(\CO'')\ne\varnothing$, so the conditions \eqref{eq:fixed-orbit-prefix-constraints} are compatible with the prefix signs fixed by $\Gamma$ at indices where $R_i=0$.
Exactly $k$ of the $h$ prefix signs $\sigma_i(\Gamma')$ with $R_i>0$ are therefore prescribed independently, giving
\begin{equation}\label{eq:beta-Gamma-Opp}
 \#\bigl(\mathscr S_\Gamma^{\pi}\cap\mathfrak G(\CO'')\bigr)
 =2^{h-k}.
\end{equation}
By \cite[Corollary~5.2]{wang}, all $2^{h-k}$ indices in $\mathscr S_\Gamma^\pi\cap\mathfrak G(\CO'')$ are realized by $\pi$.

\medskip
\noindent\textit{Step 3: sum over maximal rational orbits.}
By \cite[Theorem~7.1]{wang} and
\eqref{eq:signed-prefix-reachability}, every realized index $\Gamma'\in\mathscr S_\Gamma^\pi$ is maximal and good. Its associated rational orbit $\CO''=\CO(\Gamma')$ is maximal in $\ScO(\pi)$ by \cite[Corollary~7.2 and its proof]{wang}, and satisfies $\CO\le_{\rm rat}\CO''$ by \eqref{eq:good-order-compatibility}.
To compute $m_{\Gamma'}(\pi)$ for such a realized index $\Gamma'$, choose the even-row and odd-pair
signs in \eqref{eq:uniform-orbit-expansion} for $\CO''$ so that
$\Gamma_{\CO'',\mathbf 0}=\Gamma'$. Since $\CO''$ has $s_{\max}$ odd pairs, \eqref{eq:uniform-orbit-expansion} reads
\[
 2^{s_{\max}}\dim\Wh_{\CO''}(\pi)
 =\sum_{\mathbf c\in\mathscr R(\CO'')} w_q(\mathbf c)\,m_{\Gamma_{\CO'',\mathbf c}}(\pi).
\]
Every $\mathbf c\in\mathscr R(\CO'')$ with $\mathbf c\ne\mathbf0$ gives an index $\Gamma_{\CO'',\mathbf c}$ strictly larger than $\Gamma'$ with respect to $\le_{\rm seq}$, which is not realized by $\pi$ (as $\Gamma'$ is maximal), so $m_{\Gamma_{\CO'',\mathbf c}}(\pi)=0$. Thus only the $\mathbf c=\mathbf0$ term (with $w_q(\mathbf0)=1$) survives, yielding
\[
 m_{\Gamma'}(\pi)=2^{s_{\max}}\dim\Wh_{\CO''}(\pi).
\]
By \eqref{eq:cGamma}, an index $\Gamma'\in\mathscr S_\Gamma^\pi$ satisfies $m_{\Gamma'}(\pi)\ne0$ if and only if it is realized by $\pi$. Combining Step~2 with the discussion above, we see that the indices $\Gamma'\in\mathscr S_\Gamma^\pi$ with $m_{\Gamma'}(\pi)\ne0$ decompose as the disjoint union
\[
 \{\Gamma'\in\mathscr S_\Gamma^\pi:m_{\Gamma'}(\pi)\ne0\}
 =\bigsqcup_{\substack{\CO''\in\ScO(\pi)^{\max}\\\CO\le_{\rm rat}\CO''}}
   \bigl(\mathscr S_\Gamma^\pi\cap\mathfrak G(\CO'')\bigr).
\]
Grouping the nonzero terms of $\sum_{\Gamma'\in\mathscr S_\Gamma^\pi}m_{\Gamma'}(\pi)$ accordingly and using \eqref{eq:beta-Gamma-Opp} gives
\begin{equation}\label{eq:T-intrinsic-count}
 \sum_{\Gamma'\in\mathscr S_\Gamma^\pi}m_{\Gamma'}(\pi)
 =2^{h-k+s_{\max}}
   \sum_{\substack{\CO''\in\ScO(\pi)^{\max}\\
                   \CO\le_{\rm rat}\CO''}}
       \dim\Wh_{\CO''}(\pi).
\end{equation}
Combining \eqref{eq:orbit-lc-vs-C0}, \eqref{eq:T-intrinsic-count}
and \eqref{eq:Delta-k-z}, we obtain
\[
 \LC(m_{\CO,\mathfrak t}(X))
 =2^{-\Delta(\CO,\CO_\pi^{\rm st})}
   \sum_{\substack{\CO''\in\ScO(\pi)^{\max}\\
                   \CO\le_{\rm rat}\CO''}}
       \dim\Wh_{\CO''}(\pi),
\]
which is \eqref{eq:orbit-leading-intrinsic}.
By \cite[Theorem~1.9]{wang}, the maximal Whittaker dimensions $\dim\Wh_{\CO''}(\pi)$ ($\CO''\in\ScO(\pi)^{\max}$) are
independent of $\CO''$ and satisfy
\[
 \dim\Wh_{\CO''}(\pi)
 =\frac{|A(\CO_\pi^{\rm st})|}{|\ScO(\pi)^{\max}|\,|\bA(\CO_{\pi^{\rm J}}^{\rm st})|}.
\]
There are $Q_\pi(\CO)$ equal summands in \eqref{eq:orbit-leading-intrinsic}, which gives \eqref{eq:orbit-leading-canonical}.
\end{proof}

\begin{remark}
\label{rem:geck-taylor-endpoints}
For the rational orbit $\CO=\{0\}$, we have $m_{\{0\},\mathfrak t}(X)=D_{\mathfrak t}(X)$, since $\Gamma_{\{0\}}$ is the regular
representation of $G$. For $G=\Sp_{2N}(\Fq)$ and $\pi\in\Irr_{\mathfrak t}(G)$,
by \eqref{eq:orbit-degree} and \eqref{eq:orbit-leading-canonical}, we have
\[
 \deg D_{\mathfrak t}(X)=\tfrac12\dim\CO_\pi^{\rm st},\qquad
 {\LC(D_{\mathfrak t}(X))=\frac1{n_\pi},\qquad
 n_\pi=|\bA(\CO_{\pi^{\rm J}}^{\rm st})|.}
\]
For the second identity, note that
$Q_\pi(\{0\})=|\ScO(\pi)^{\max}|$, and that
$\Delta(\{0\},\CO_\pi^{\rm st})$ is the number of distinct even parts of the
partition of $\CO_\pi^{\rm st}$. Hence
$2^{\Delta(\{0\},\CO_\pi^{\rm st})}=|A(\CO_\pi^{\rm st})|$
\cite[Section~6.1]{CM}.
For the standard rational representatives, the sum of the
GGGRs over the rational orbits contained in $\CO_\pi^{\rm st}$ is Geck's
stable GGGR \cite{Geck99}. The multiplicity of $\pi$ in this sum is
$|A(\CO_\pi^{\rm st})|/n_\pi$. Thus
\eqref{eq:orbit-leading-canonical} agrees with the generic-degree and
GGGR expressions in \cite[(2.3),(2.8)]{Geck99},
\cite[Theorem~3.1 and \S8]{Taylor13}, and \cite[Lemma~14.15]{Taylor16}.
\end{remark}

\section{Orthogonal groups and theta correspondence}\label{sec:orth-theta}
By the finite-field Gomez--Zhu theorem \cite[Theorem~8.6]{wang} (see also \cite{GZ14}), the generalized Whittaker models of a representation $\pi$ of an orthogonal group can be transferred to those of its theta lift to a symplectic group. In this section we combine this transfer with the symplectic results of Sections~\ref{sec:minimality-degree}--\ref{sec:orbit-multiplicities} to prove Theorems~\ref{thm:mainA} and~\ref{thm:mainB} for full orthogonal groups.

\subsection{Theta lifts and orbit transfer}\label{subsec:theta-eta-choice}
Let $V$ be a nondegenerate quadratic space of dimension $n$ over $\Fq$ with symmetric bilinear form $B_V$,
and let $W$ be a symplectic space of dimension $2n$ with symplectic form $B_W$. Put
$G=\OO(V)$ and $G'=\Sp(W)$,
so that $(G,G')$ is a finite reductive dual pair in $\Sp(V\otimes W)$.

Following \cite[Section~8]{wang}, we use the theta correspondence
to pass from $\OO(V)$ to $\Sp(W)$. Although, the Coxeter number of $\Sp(W)$
is $2n$, roughly twice that of $\OO(V)$. As Wang observes in
\cite[v2, Section~1.3(3)]{wang}, the assumption
on $q$ for $\OO_n$ also suffices for applying symplectic results of Sections~\ref{sec:minimality-degree}--\ref{sec:orbit-multiplicities}.

Let $\omega_{G,G'}$ be the
restriction to $G\times G'$ of the standard Weil representation of
$\Sp(V\otimes W)$ attached to the fixed additive character $\psi$
\cite{Ger}.
For $\sigma\in\Irr(G)$, the \emph{big theta lift} of $\sigma$ to $G'$ is defined by
\begin{equation*}
\begin{aligned}
 \Theta_{n,2n}(\sigma)
 &:=\Hom_{G}(\sigma,\omega_{G,G'}),\\
 \omega_{G,G'}
 &\simeq\bigoplus_{\sigma\in\Irr(G)}
          \sigma\otimes\Theta_{n,2n}(\sigma),
\end{aligned}
\end{equation*}
where $G'$ acts via its action on $\omega_{G,G'}$.
For $T\in\Hom(V,W)$, its adjoint $T^*\in\Hom(W,V)$ is defined by
\[
 B_W(Tv,w)=B_V(v,T^*w)\qquad(v\in V,\ w\in W).
\]
The moment maps (see \cite[Section~8]{wang}) are morphisms of affine spaces over $\Fq$:
\[
\xymatrix@C=3em@R=1.2em{
 \mathfrak o(V)
 & \Hom(V,W)\ar[l]_{\varphi}\ar[r]^{\varphi'}
 & \mathfrak{sp}(W)\\
 T^*T
 & T\ar@{|->}[l]\ar@{|->}[r]
 & TT^*
}
\]

Let
\[
 \CO=[(d_1^{r_1},a_1),\ldots,(d_u^{r_u},a_u)],
 \qquad d_1>\cdots>d_u>0,
\]
be a rational nilpotent orbit under $G=\OO(V)$, with partition
$\lambda=(\lambda_1,\ldots,\lambda_\ell)$ of $n=\dim V$ (where $\ell=\ell(\lambda)\le n$). Adding a first column of length $n$ produces a partition of $2n=\dim W$ in which every odd part has even multiplicity (note that $n-\ell=\sum_i r_i(d_i-1)$ is even), namely
\begin{equation}\label{eq:partition-first-column-lift}
 \lambda^{\uparrow}
 :=(\lambda_1+1,\ldots,\lambda_\ell+1,1^{n-\ell}),
\end{equation}
and a corresponding rational nilpotent orbit $\CO^\uparrow$ in $\mathfrak{sp}(W)\simeq\mathfrak{sp}_{2n}(\Fq)$:
\begin{equation}\label{eq:orth-first-column-lift}
 \CO^{\uparrow}
 :=[((d_1+1)^{r_1},a_1),\ldots,((d_u+1)^{r_u},a_u),1^{n-\ell}].
\end{equation}
By \cite[Propositions~8.3--8.4]{wang}, there exists an injective linear map $T\in\Hom(V,W)$ (that is, $\ker T=0$) with $T^*T\in\CO$, and for every such $T$ the $G'$-orbit of $TT^*$ is $\CO^\uparrow$.
For a stable nilpotent orbit $\Ost{\nu}$ of $\bO(V)$, we similarly define $\Ost{\nu}^\uparrow:=\Ost{\nu^\uparrow}$ in $\mathfrak{sp}(W)$ (with $\nu^\uparrow$ as in \eqref{eq:partition-first-column-lift}); thus $(\CO^\uparrow)^{\rm st}=(\CO^{\rm st})^\uparrow$.
By \cite[Proposition~8.4]{wang}, $\CO^\uparrow$ is the unique generalized lift of $\CO$ in $\operatorname{im}\varphi'$ (in the sense of \cite[Section~8.1]{wang}) whose partition has $n$ rows.
Moreover, by the same proposition, for any two rational nilpotent orbits $\CO_1,\CO_2\subset\mathfrak o(V)$ we have
\begin{equation}\label{eq:first-column-order}
 \CO_1\le_{\rm rat}\CO_2
 \quad\Longleftrightarrow\quad
 \CO_1^\uparrow\le_{\rm rat}\CO_2^\uparrow.
\end{equation}
Since $\CO\mapsto\CO^\uparrow$ preserves $r_i$ and $a_i$ and appends an even number $n-\ell$ of $1$'s, it is compatible with the orbit identification of Section~\ref{subsec:rational-orbit-order}.

\begin{lemma}\label{lem:first-column-dimension}
With the notation above (where $\dim V=n$ and $\dim W=2n$), for every rational nilpotent orbit $\CO\subset\mathfrak o(V)$ we have
\begin{equation}\label{eq:first-column-dimension}
 \dim(\CO^\uparrow)^{\rm st}=\dim\CO^{\rm st}+n(n+1).
\end{equation}
\end{lemma}
\begin{proof}
Let $\lambda$ be the partition of $\CO$. Since $\lambda^\uparrow$ is obtained from $\lambda$ by adding a first column of length $n$, we have $(\lambda^\uparrow)^t=(n,\lambda_1^t,\lambda_2^t,\ldots)$, and the number of odd parts of $\lambda^\uparrow$ is $n-\#\{i:\lambda_i\text{ is odd}\}$.
Choose representatives $X\in\CO\subset\mathfrak o(V)$ (with partition $\lambda$) and $Y\in\CO^\uparrow\subset\mathfrak{sp}(W)$ (with partition $\lambda^\uparrow$). By the centralizer dimension formulas \cite[Corollary~6.1.4]{CM},
\[
 \dim C_{\bO_n}(X)
 =\frac12\left(\sum_j(\lambda_j^t)^2-\#\{i:\lambda_i\text{ is odd}\}\right),
\]
and
\[
\begin{aligned}
 \dim C_{\bSp_{2n}}(Y)
 &=\frac12\left(n^2+\sum_j(\lambda_j^t)^2+n-\#\{i:\lambda_i\text{ is odd}\}\right)\\
 &=\dim C_{\bO_n}(X)+\frac{n(n+1)}2.
\end{aligned}
\]
Since $\dim\bSp_{2n}=n(2n+1)$ and $\dim\bO_n=n(n-1)/2$, subtracting the centralizer dimensions yields \eqref{eq:first-column-dimension}.
\end{proof}

Fix $\pi\in\Irr(G)$, and write $\pi^\vee$ for its contragredient. Define $\tau\in\Irr(G)$ by
\[
 \tau:=\begin{cases}
  \pi^\vee, & n\text{ even},\\
  \text{the unique element of }\{\pi^\vee,\pi^\vee\otimes\det\nolimits_G\}\text{ with }\delta_{-1}=+1, & n\text{ odd},
 \end{cases}
\]
where, for each irreducible constituent $\sigma'\in\cE(G',s')$ of $\Theta_{n,2n}(\tau)$, $\delta_{-1}\in\{\pm1\}$ is the discriminant of the quadratic space defining the $[-1]$-primary factor $(G')^*_{[-1]}(s')\simeq\OO_{2r_{-1}}^{\delta_{-1}}(\Fq)$ as in Section~\ref{sec:classification} (by \cite[Theorem~6.1(2)]{WangJordanCorrected}, $\delta_{-1}$ is independent of $\sigma'$ and is reversed when $\tau$ is twisted by $\det\nolimits_G$). Since $\det\nolimits_G$ is trivial on unipotent subgroups, twisting by $\det\nolimits_G$ does not change Whittaker multiplicities.

The following proposition summarizes the structure of $\Theta_{n,2n}(\tau)$ needed in the next subsection.

\begin{proposition}\label{prop:theta-distinguished}
Let $G=\OO(V)$, let $\pi\in\Irr(G)$, and let $\tau\in\Irr(G)$ be chosen as above. Then the big theta lift $\Theta_{n,2n}(\tau)$ decomposes as
\begin{equation}\label{eq:selected-theta-is-eta}
 \Theta_{n,2n}(\tau)=\pi^\theta\oplus\bigoplus_{j\in J}\pi_j,
\end{equation}
where $\pi^\theta$ is the irreducible constituent of $\Theta_{n,2n}(\tau)$ chosen in \cite[Theorems~8.7 and~8.12]{wang}, and $\pi_j$ ($j\in J$, counted with multiplicity) are the remaining irreducible constituents, such that:
\begin{enumerate}[label=\textup{(\roman*)},leftmargin=2em]
 \item the maximal rational wavefront orbits of $\pi^\theta$ are given by
 \begin{equation}\label{eq:theta-top-separation}
  \ScO(\pi^\theta)^{\max}
  =\{(\CO'')^\uparrow:\CO''\in\ScO(\pi)^{\max}\};
\end{equation}
 \item the stable wavefront orbits of all remaining constituents have strictly smaller dimension:
 \begin{equation}\label{eq:theta-strict-separation}
  \dim\CO_{\pi_j}^{\rm st}<\dim\CO_{\pi^\theta}^{\rm st}
  \qquad(j\in J).
\end{equation}
\end{enumerate}
\end{proposition}
\begin{proof}
The existence of the distinguished constituent $\pi^\theta$ in \eqref{eq:selected-theta-is-eta} and the comparison of maximal rational wavefront orbits \eqref{eq:theta-top-separation} are established in \cite[Proposition~8.4, proof of Theorem~8.7, and Lemma~8.10]{wang} for even $n$ and in \cite[proof of Theorem~8.12]{wang} for odd $n$. The strict dimension inequality \eqref{eq:theta-strict-separation} follows from \cite[Lemma~8.10 and proof of Theorem~8.12]{wang}.
\end{proof}

We also recall the finite-field Gomez--Zhu theorem \cite[Theorem~8.6]{wang}.
Let $\CO'$ be a rational nilpotent orbit contained in $\operatorname{im}\varphi'$, and let
$\CO=\nabla^{\rm gen}_{W,V}(\CO')$ be its generalized descent, realized by a map $T\in\Hom(V,W)$ with $TT^*\in\CO'$ and $T^*T\in\CO$.
Let $(L,L')$ be the residual finite dual pair of \cite[Definition~8.5]{wang}, formed by the isometry groups of the residual spaces attached to $\ker T$ and $\ker T^*$ (in our application below to $\CO'=\CO^\uparrow$ with $\ker T=0$, it is simply $(L,L')=(1,\Sp_{n-\ell}(\Fq))$, where $\ell$ is the length of the partition of $\CO$).
For a representation $\pi'$ of $G'$ and $\tau'\in\Irr(L')$, define
\[
 \Wh_{\CO',\tau'}(\pi')
 :=\Hom_{L'}(\tau',\Wh_{\CO'}(\pi')),
\]
and define $\Wh_{\CO,\tau_0}(\widetilde\pi)$ similarly for representations
$\widetilde\pi$ of $G$ and $\tau_0\in\Irr(L)$.
For $\sigma\in\Irr(G)$, \cite[Theorem~8.6(i)]{wang} gives
\begin{equation}\label{eq:orth-gomez-zhu}
 \dim\Wh_{\CO',\tau'}(\Theta_{n,2n}(\sigma))
 =\dim\Wh_{\CO,\Theta(\tau')^\vee}(\sigma^\vee),
\end{equation}
where $\Theta(\tau')$ is the theta lift from $L'$ to $L$.
Moreover, by \cite[Theorem~8.6(ii)]{wang}, every rational orbit in the rational wavefront set of $\Theta_{n,2n}(\sigma)$ lies in $\operatorname{im}\varphi'$:
\[
 \CO'\in\ScO\bigl(\Theta_{n,2n}(\sigma)\bigr)
 \quad\Longrightarrow\quad
 \CO'\subseteq\operatorname{im}\varphi'.
\]

\subsection{The theta sum and its highest-degree term}
Let $\boldsymbol\pi=\{\pi_q\}_q$ be a $q$-system of type $\mathfrak t$ for $G=\OO(V)$. For each $\pi_q$, choose $\tau_q$ as in Section~\ref{subsec:theta-eta-choice} and decompose its theta lift as in \eqref{eq:selected-theta-is-eta}, with irreducible constituents $\pi_q^\theta$ and $\pi_{q,j}$ ($j\in J$).

To apply the symplectic polynomiality result (Proposition~\ref{prop:uniform-orbit-expansion}) to each constituent $\pi_q^\theta$ and $\pi_{q,j}$ in \eqref{eq:selected-theta-is-eta}, we first observe that their types (together with their multiplicities) depend only on $\mathfrak t$ and $\delta$. Indeed, the type of $\pi_q^\vee$ is determined by $\mathfrak t$. By the description of the theta correspondence in terms of Lusztig's Jordan decomposition \cite[Theorem~2.2 and proofs of Theorems~8.7 and~8.12]{wang} (under the conventions of \cite{WangJordanCorrected}):
\begin{itemize}[leftmargin=1.8em]
 \item every irreducible constituent of $\Theta_{n,2n}(\tau_q)$ has the same non-$[\pm1]$ partitions as $\tau_q$;
 \item the symbols $\Lambda_{\pm1}$ in the Lusztig parameters of the irreducible constituents of $\Theta_{n,2n}(\tau_q)$ (and their multiplicities in $\Theta_{n,2n}(\tau_q)$) are determined combinatorially by the symbols in the Lusztig parameters of $\tau_q$, $\delta$, and (when $n$ is odd) the discriminant $\delta_{-1}=+1$ fixed by our choice of $\tau_q$.
\end{itemize}
Consequently, the types of $\pi_q^\theta$ and $\pi_{q,j}$ form a finite list, denoted $\mathfrak t^\theta$ and $\mathfrak t_j$ ($j\in J$, with repetitions recording multiplicities), which is determined uniquely by $\mathfrak t$ and $\delta$ and is independent of $q$.

\begin{proposition}\label{prop:orth-theta-polynomial}
For a rational nilpotent orbit $\CO\subset\CN^F$ under $G$, there is a polynomial $m_{\CO,\mathfrak t}(X)\in\mathbb Q[X]$ such that, for every sufficiently large $q$ and every $\pi_q\in\Irr_{\mathfrak t}(\bG(\Fq))$,
\begin{equation*}
 m_{\CO}(\pi_q)=m_{\CO,\mathfrak t}(q).
\end{equation*}
Moreover,
\begin{equation}\label{eq:orth-polynomial-theta-sum}
 m_{\CO,\mathfrak t}(X)
 =m_{\CO^\uparrow,\mathfrak t^\theta}(X)
   +\sum_{j\in J}m_{\CO^\uparrow,\mathfrak t_j}(X).
\end{equation}
In particular, $m_\CO(\pi_q)$ does not depend on the choice of realization of $\mathfrak t$.
\end{proposition}
\begin{proof}
Fix $\delta\in\{\pm1\}$ and consider first $q$ with $\varepsilon_q(-1)=\delta$, setting $m_{\CO,\mathfrak t,\delta}(X):=m_{\CO^\uparrow,\mathfrak t^\theta}(X)+\sum_{j\in J}m_{\CO^\uparrow,\mathfrak t_j}(X)$ (with $\mathfrak t^\theta$ and $\mathfrak t_j$, $j\in J$, depending on $\delta$). Let $\ell$ be the number of parts of the partition of $\CO$.
By \cite[Propositions~8.3--8.4]{wang}, the generalized descent of $\CO^\uparrow$ is $\CO$, and any map $T\in\Hom(V,W)$ with $TT^*\in\CO^\uparrow$ and $T^*T\in\CO$ satisfies $\ker T=0$. Since the multiplicity space of the part $1$ in $\CO^\uparrow$ has dimension $n-\ell$ (by \eqref{eq:orth-first-column-lift}), the residual dual pair attached to $\CO^\uparrow$ and $\CO$ is
\[
 (L,L')=(1,\Sp_{n-\ell}(\Fq)).
\]
Since $L=1$, the Weil representation of $(L,L')$ is trivial, so $\Theta(\mathbf{1}_{L'})=\mathbf{1}_L$ and $\Theta(\tau')=0$ for every nontrivial $\tau'\in\Irr(L')$. Applying the Gomez--Zhu identity \eqref{eq:orth-gomez-zhu} with $\CO'=\CO^\uparrow$ and $\sigma=\tau_q$, we see that $\Wh_{\CO^\uparrow,\tau'}(\Theta_{n,2n}(\tau_q))=0$ for every nontrivial $\tau'\in\Irr(L')$; hence $L'$ acts trivially on $\Wh_{\CO^\uparrow}(\Theta_{n,2n}(\tau_q))$. Evaluating \eqref{eq:orth-gomez-zhu} at $\tau'=\mathbf{1}_{L'}$ therefore gives
\[
\begin{aligned}
 \dim\Wh_{\CO^\uparrow}(\Theta_{n,2n}(\tau_q))
 &=\dim\Wh_{\CO^\uparrow,\mathbf{1}_{L'}}(\Theta_{n,2n}(\tau_q))\\
 &=\dim\Wh_{\CO,\mathbf{1}_L}(\tau_q^\vee)
  =\dim\Wh_{\CO}(\tau_q^\vee)
  =\dim\Wh_{\CO}(\pi_q),
\end{aligned}
\]
where the last equality holds because $\tau_q^\vee\in\{\pi_q,\pi_q\otimes\det\nolimits_G\}$ and $\det\nolimits_G$ is trivial on unipotent subgroups.

Decomposing $\Theta_{n,2n}(\tau_q)$ via Proposition~\ref{prop:theta-distinguished} and applying Proposition~\ref{prop:uniform-orbit-expansion} to each irreducible constituent, we obtain
\begin{equation}\label{eq:orth-mult-expansion}
\begin{aligned}
 m_{\CO}(\pi_q)
 &=m_{\CO^\uparrow}(\Theta_{n,2n}(\tau_q))\\
 &=m_{\CO^\uparrow}(\pi_q^\theta)
   +\sum_{j\in J}m_{\CO^\uparrow}(\pi_{q,j})\\
 &=m_{\CO^\uparrow,\mathfrak t^\theta}(q)
   +\sum_{j\in J}m_{\CO^\uparrow,\mathfrak t_j}(q)
  =m_{\CO,\mathfrak t,\delta}(q).
\end{aligned}
\end{equation}

We now verify $m_{\CO,\mathfrak t,+1}(X)=m_{\CO,\mathfrak t,-1}(X)$. When $n$ is even, $G'=\Sp_{2n}(\Fq)$ has even rank $n$, so $\omega_{G,G'}$ coincides with the modified Weil representation of \cite[Section~4.2]{WangJordanCorrected}; by \cite[Theorem~6.1(3)]{WangJordanCorrected}, the multiset of constituent types of $\Theta_{n,2n}(\pi_q^\vee)$ depends only on $\mathfrak t$, and Proposition~\ref{prop:uniform-orbit-expansion} gives the equality. When $n$ is odd, the Gomez--Zhu argument above applies equally to $\pi_q^\vee$ and $\pi_q^\vee\otimes\det_G$, giving
\[
 m_\CO(\pi_q)
 =m_{\CO^\uparrow}\bigl(\Theta_{n,2n}(\pi_q^\vee)\bigr)
 =m_{\CO^\uparrow}\bigl(\Theta_{n,2n}(\pi_q^\vee\otimes\det\nolimits_G)\bigr),
\]
and hence
\[
 m_\CO(\pi_q)
 =\frac12\,m_{\CO^\uparrow}\!\Bigl(\Theta_{n,2n}(\pi_q^\vee)\oplus\Theta_{n,2n}(\pi_q^\vee\otimes\det\nolimits_G)\Bigr).
\]
Since $\Theta_{n,2n}(\pi_q^\vee)\oplus\Theta_{n,2n}(\pi_q^\vee\otimes\det_G)$ is invariant under twisting the Weil representation by $\det_G$, it agrees with the corresponding sum for the modified Weil representation of \cite[Section~4.2]{WangJordanCorrected}; moreover, by \cite[Theorem~6.1(2)]{WangJordanCorrected}, summing over the two twists eliminates the sign condition on the action of $-I\in G$, so the resulting multiset of constituent types depends only on $\mathfrak t$. Proposition~\ref{prop:uniform-orbit-expansion} again gives $m_{\CO,\mathfrak t,+1}(X)=m_{\CO,\mathfrak t,-1}(X)=:m_{\CO,\mathfrak t}(X)$.
\end{proof}

\begin{corollary}\label{cor:orth-whittaker-type-invariance}
Let $G$ be symplectic or orthogonal. Under our standing assumptions, for sufficiently large $q$
and $\pi,\widetilde\pi\in\Irr_{\mathfrak t}(\bG(\Fq))$,
\[
 m_\CO(\pi)=m_\CO(\widetilde\pi)
 \qquad\text{for every rational nilpotent orbit }\CO.
\]
Consequently,
\[
 \ScO(\pi)=\ScO(\widetilde\pi),\qquad
 \ScO(\pi)^{\max}=\ScO(\widetilde\pi)^{\max}.
\]
\end{corollary}
\begin{proof}
By Propositions~\ref{prop:uniform-orbit-expansion}
and~\ref{prop:orth-theta-polynomial}, both $m_\CO(\pi)$ and
$m_\CO(\widetilde\pi)$ are the value at $q$ of a polynomial depending only
on $\CO$ and $\mathfrak t$.
Thus the rational wavefront sets, and hence their maximal elements, coincide.
\end{proof}

\begin{proposition}\label{thm:orthogonal-transfer}
Assume $m_{\CO,\mathfrak t}(X)\ne0$. Then
\begin{align}
 \deg m_{\CO,\mathfrak t}(X)
 &=\deg m_{\CO^\uparrow,\mathfrak t^\theta}(X),
 \label{eq:orth-degree-transfer}\\
 \LC(m_{\CO,\mathfrak t}(X))
 &=\LC(m_{\CO^\uparrow,\mathfrak t^\theta}(X)).
 \label{eq:orth-lc-transfer}
\end{align}
Consequently, for every sufficiently large $q$ and every
$\pi_q\in\Irr_{\mathfrak t}(\bG(\Fq))$,
\begin{equation}\label{eq:orth-degree-final}
 \deg m_{\CO,\mathfrak t}(X)
 =\frac12\left(\dim\CO_{\pi_q}^{\rm st}-\dim\CO^{\rm st}\right).
\end{equation}
\end{proposition}
\begin{proof}
Fix a sufficiently large $q$ and $\pi_q\in\Irr_{\mathfrak t}(\bG(\Fq))$.
We first show that $m_{\CO^\uparrow,\mathfrak t^\theta}(X)\ne0$.
By Proposition~\ref{prop:orth-theta-polynomial}, $\CO\in\ScO(\pi_q)$.
Thus, by \cite[Theorem~1.1]{wang}, we have
$\CO\le_{\rm rat}\CO''$ for some $\CO''\in\ScO(\pi_q)^{\max}$.
By \eqref{eq:first-column-order} and \eqref{eq:theta-top-separation},
we have
\[
 \CO^\uparrow\le_{\rm rat}(\CO'')^\uparrow
 \in\ScO(\pi_q^\theta)^{\max}.
\]
By \cite[Corollary~7.2(ii)]{wang}, $\CO^\uparrow\in\ScO(\pi_q^\theta)$.
In particular, $m_{\CO^\uparrow,\mathfrak t^\theta}(X)\ne0$.

For each $j\in J$ with $m_{\CO^\uparrow,\mathfrak t_j}(X)\ne0$, we obtain
from Theorem~\ref{prop:orbit-degree} and
\eqref{eq:theta-strict-separation}
\[
\begin{aligned}
 \deg m_{\CO^\uparrow,\mathfrak t_j}(X)
 &=\frac12\left(\dim\CO_{\pi_{q,j}}^{\rm st}
                         -\dim(\CO^\uparrow)^{\rm st}\right)\\
 &<\frac12\left(\dim\CO_{\pi_q^\theta}^{\rm st}
                         -\dim(\CO^\uparrow)^{\rm st}\right)
 =\deg m_{\CO^\uparrow,\mathfrak t^\theta}(X).
\end{aligned}
\]
Thus $m_{\CO^\uparrow,\mathfrak t^\theta}(X)$ is the unique summand of
highest degree in \eqref{eq:orth-polynomial-theta-sum}.
This proves \eqref{eq:orth-degree-transfer} and \eqref{eq:orth-lc-transfer}.

Finally, combining \eqref{eq:theta-top-separation} with
Lemma~\ref{lem:first-column-dimension}, we obtain
\begin{equation*}
 \dim\CO_{\pi_q^\theta}^{\rm st}-\dim(\CO^\uparrow)^{\rm st}
 =\dim\CO_{\pi_q}^{\rm st}-\dim\CO^{\rm st},
\end{equation*}
which proves \eqref{eq:orth-degree-final}.
\end{proof}

To translate the leading-coefficient formula of Theorem~\ref{prop:orbit-leading-coefficient} from $G'=\Sp(W)$ to $G=\OO(V)$, we define the analogue of $\Delta$ from \eqref{eq:Delta-stable-partitions} for orthogonal groups. Let $\nu=(\nu_1,\ldots,\nu_\ell)$ be the partition of a nilpotent orbit of $\bO(V)$, where $\ell=\ell(\nu)\le n$. Appending $n-\ell$ zeros, write
\[
 (\nu_1,\ldots,\nu_\ell,0^{n-\ell})
 =(d_1^{m_1},\ldots,d_u^{m_u}),\qquad d_1>\cdots>d_u\ge0,
\]
and set $M_v(\nu):=\sum_{a=1}^v m_a$ with $M_0(\nu):=0$. Since every even part $d_v\ge0$ has even multiplicity $m_v$ (including $n-\ell$ when $d_u=0$), we may define
\begin{align*}
 J_{\rm end}^{\rm O}(\nu)
 &:=\{M_v(\nu):d_v\text{ is odd}\}
   \cup\left\{M_{v-1}(\nu)+2b:
       \substack{d_v\text{ is even}\\1\le b\le m_v/2}\right\},\\
 J_{\rm even}^{\rm O}(\nu)
 &:=\left\{M_{v-1}(\nu)+2b-1:
       \substack{d_v\text{ is even}\\1\le b\le m_v/2}\right\}.
\end{align*}
For a rational nilpotent orbit $\CO$ and a stable nilpotent orbit ${\CO'}^{\rm st}$ of $\bO(V)$ with $\CO^{\rm st}\le{\CO'}^{\rm st}$, let $\lambda$ and $\lambda'$ be their respective partitions, put $S_i(\lambda):=\sum_{j=1}^i\lambda_j$ (padding with zeros), and define
\begin{equation}\label{eq:Delta-orth-combinatorial}
\begin{aligned}
 \Delta_{\OO}(\CO,{\CO'}^{\rm st})
 &:=\#\{i\in J_{\rm end}^{\rm O}(\lambda'):
                     S_i(\lambda)<S_i(\lambda')\}\\
 &\quad-\#\{i\in J_{\rm even}^{\rm O}(\lambda'):
                     S_i(\lambda)<S_i(\lambda')\}.
\end{aligned}
\end{equation}
Adding $1$ to each of the $n$ parts of $(\lambda'_1,\ldots,\lambda'_{\ell(\lambda')},0^{n-\ell(\lambda')})$ yields $(\lambda')^\uparrow$ (see \eqref{eq:partition-first-column-lift}), reversing the parity of every part, preserving multiplicities, and adding $i$ to the $i$-th prefix sum; hence $J_{\rm end}^{\rm O}(\lambda')=J_{\rm end}((\lambda')^\uparrow)$, $J_{\rm even}^{\rm O}(\lambda')=J_{\rm odd}((\lambda')^\uparrow)$, and
\begin{equation}\label{eq:Delta-orth}
 \Delta_{\OO}(\CO,{\CO'}^{\rm st})
 =\Delta\bigl(\CO^\uparrow,({\CO'}^{\rm st})^\uparrow\bigr).
\end{equation}
Recall from Section~\ref{sec:intro} that $\pi^{\rm J}=\CL_s(\pi)$ for $\pi\in\cE(G,s)$, and that $c_G=1$ when $n$ is even and $c_G=2$ when $n$ is odd.

\begin{corollary}\label{cor:orth-leading}
For sufficiently large $q$, let $\pi\in\Irr_{\mathfrak t}(G)$ and $\CO\in\ScO(\pi)$.
Then
\begin{equation}\label{eq:orth-leading-intrinsic}
 \LC(m_{\CO,\mathfrak t}(X))
 =2^{-\Delta_{\OO}(\CO,\CO_\pi^{\rm st})}
 \sum_{\substack{\CO''\in\ScO(\pi)^{\max}\\\CO\le_{\rm rat}\CO''}}
 \dim\Wh_{\CO''}(\pi).
\end{equation}
Moreover, for $\pi\in\cE(G,s)$,
\begin{equation}\label{eq:orth-leading-canonical}
 \LC(m_{\CO,\mathfrak t}(X))
 =2^{-\Delta_{\OO}(\CO,\CO_\pi^{\rm st})}
 \frac{Q_\pi(\CO)|A(\CO_\pi^{\rm st})|}
      {c_G|\ScO(\pi)^{\max}|\,|\bA(\CO_{\pi^{\rm J}}^{\rm st})|}.
\end{equation}
\end{corollary}
\begin{proof}
Since $\CO\in\ScO(\pi)$, we have $m_{\CO,\mathfrak t}(X)\ne0$ and $\CO^\uparrow\in\ScO(\pi^\theta)$ (as in the proof of Proposition~\ref{thm:orthogonal-transfer}). Combining \eqref{eq:orth-lc-transfer} with Theorem~\ref{prop:orbit-leading-coefficient} applied to $\pi^\theta$ and $\CO^\uparrow$ gives
\[
 \LC(m_{\CO,\mathfrak t}(X))
 =\LC(m_{\CO^\uparrow,\mathfrak t^\theta}(X))
 =2^{-\Delta(\CO^\uparrow,\CO_{\pi^\theta}^{\rm st})}
 \sum_{\substack{\CO'''\in\ScO(\pi^\theta)^{\max}\\
                  \CO^\uparrow\le_{\rm rat}\CO'''}}
 \dim\Wh_{\CO'''}(\pi^\theta).
\]
By direct calculation, we have
\begin{align}
 \Delta(\CO^\uparrow,\CO_{\pi^\theta}^{\rm st})
 &=\Delta_{\OO}(\CO,\CO_\pi^{\rm st}),
 \label{eq:theta-transport-defect}\\
 Q_{\pi^\theta}(\CO^\uparrow)&=Q_\pi(\CO),
 \label{eq:theta-transport-Q}\\
 \dim\Wh_{(\CO'')^\uparrow}(\pi^\theta)
 &=\dim\Wh_{\CO''}(\pi)
 \qquad(\CO''\in\ScO(\pi)^{\max}).
 \label{eq:theta-transport-whittaker}
\end{align}
Indeed, \eqref{eq:theta-top-separation} implies $\CO_{\pi^\theta}^{\rm st}=(\CO_\pi^{\rm st})^\uparrow$, so \eqref{eq:Delta-orth} yields \eqref{eq:theta-transport-defect}. Moreover, by \eqref{eq:first-column-order} and \eqref{eq:theta-top-separation}, the map $\CO''\mapsto(\CO'')^\uparrow$ is a bijection from $\{\CO''\in\ScO(\pi)^{\max}:\CO\le_{\rm rat}\CO''\}$ onto $\{\CO'''\in\ScO(\pi^\theta)^{\max}:\CO^\uparrow\le_{\rm rat}\CO'''\}$, which gives \eqref{eq:theta-transport-Q}. For \eqref{eq:theta-transport-whittaker}, fix $\CO''\in\ScO(\pi)^{\max}$. Then $(\CO'')^\uparrow\in\ScO(\pi^\theta)^{\max}$ by \eqref{eq:theta-top-separation}, and hence \eqref{eq:theta-strict-separation} gives
\[
 \dim\CO_{\pi_j}^{\rm st}
 <\dim\CO_{\pi^\theta}^{\rm st}
 =\dim((\CO'')^\uparrow)^{\rm st}
 \qquad(j\in J),
\]
so $m_{(\CO'')^\uparrow}(\pi_j)=0$ for every $j\in J$. Applying \eqref{eq:orth-mult-expansion} to $\CO''$ therefore yields \eqref{eq:theta-transport-whittaker}.

Substituting \eqref{eq:theta-transport-defect} and \eqref{eq:theta-transport-whittaker} (along with the bijection $\CO''\mapsto(\CO'')^\uparrow$) into the leading-coefficient expression above gives \eqref{eq:orth-leading-intrinsic}. Finally, by \cite[Theorems~8.11(v) and~8.12(v)]{wang}, for every $\CO''\in\ScO(\pi)^{\max}$ we have
\[
 \dim\Wh_{\CO''}(\pi)
 =\frac{|A(\CO_\pi^{\rm st})|}{c_G|\ScO(\pi)^{\max}|\,|\bA(\CO_{\pi^{\rm J}}^{\rm st})|},
\]
and summing these $Q_\pi(\CO)$ equal terms in \eqref{eq:orth-leading-intrinsic} yields \eqref{eq:orth-leading-canonical}.
\end{proof}

\begin{proof}[Proof of Theorems~\ref{thm:mainA} and~\ref{thm:mainB}]
For $G=\Sp_{2N}(\Fq)$, polynomiality and independence of the realization follow from Propositions~\ref{prop:exact-poly} and~\ref{prop:uniform-orbit-expansion}, while the degree and leading-coefficient formulas are given by Theorems~\ref{prop:orbit-degree} and~\ref{prop:orbit-leading-coefficient}. For $G=\OO(V)$, the corresponding statements follow from Propositions~\ref{prop:orth-theta-polynomial} and~\ref{thm:orthogonal-transfer} and Corollary~\ref{cor:orth-leading}. In both cases, since $m_\CO(\pi)=m_{\CO,\mathfrak t}(q)$ for all sufficiently large $q$ and all $\pi\in\Irr_{\mathfrak t}(\bG(\Fq))$, and a nonzero polynomial has only finitely many roots, we have $m_{\CO,\mathfrak t}(X)\ne0$ if and only if $\CO\in\ScO(\pi)$ for all such $q$ and $\pi$.
\end{proof}

\appendix
\renewcommand{\thetheorem}{\Alph{section}.\arabic{theorem}}
\renewcommand{\theequation}{\Alph{section}.\arabic{equation}}
\section{Symbol interlacing}\label{app:symbol}
We prove the symbol inequalities used in Proposition~\ref{prop:stagewise} and the symbol part of Proposition~\ref{cor:first-descent-rule}.

\begin{lemma}\label{lem:transpose-preceq}
If $P\preccurlyeq Q$, then
\begin{equation}\label{eq:transpose-preceq}
 Q_i^t\ge P_i^t\ge Q_{i+1}^t\qquad(i\ge1).
\end{equation}
\end{lemma}
\begin{proof}
Since $Q_r-1\le P_r\le Q_r$, we have $\{r:Q_r\ge i+1\}\subseteq\{r:P_r\ge i\}\subseteq\{r:Q_r\ge i\}$ for $i\ge1$; taking cardinalities gives \eqref{eq:transpose-preceq}.
\end{proof}

\begin{lemma}\label{lem:symbol-first-descent-comparison}
Suppose $\mathbf{b}_t(\Lambda_{\pm1};\Lambda'_{\pm1})\ne0$ for some $t\in\{\pm1\}$. Then
\begin{equation}\label{eq:symbol-tail-inequality}
 r_k^{\rm fd}(\Lambda'_{-\xi})
 \ge r_{k+1}^{\rm fd}(\Lambda_\xi)
 \qquad(k\ge0,\ \xi=\pm1).
\end{equation}
For each $\xi$, equality at $k=0$ holds if and only if $\Lambda_\xi$ and $\Lambda'_{-\xi}$ satisfy the formula in Proposition~\ref{cor:first-descent-rule}, in which case
\begin{equation}\label{eq:symbol-first-tail}
 r_k^{\rm fd}(\Lambda'_{-\xi})
 =r_{k+1}^{\rm fd}(\Lambda_\xi)\qquad(k\ge0).
\end{equation}
\end{lemma}
\begin{proof}
Let $\delta_{-1},\delta'_{-1}$ be the discriminants attached to $\Lambda_{-1},\Lambda'_{-1}$. Fix $\xi\in\{\pm1\}$ and set
\[
 (\gamma,\gamma'):=
 \begin{cases}
  (\delta'_{-1},-\delta t\delta'_{-1}), & \xi=1,\\
  (t\delta_{-1},-\delta_{-1}), & \xi=-1.
 \end{cases}
\]
By \eqref{eq:GGP-symbol-rule}, \eqref{eq:AC-symbol} and Lemma~\ref{lem:transpose-preceq},
\begin{equation}\label{eq:dual-interlace-common}
\begin{aligned}
 \Upsilon({}^d\Lambda_\xi)^\gamma_i
 &\ge \Upsilon({}^d\Lambda'_{-\xi})^{\gamma'}_i
 \ge \Upsilon({}^d\Lambda_\xi)^\gamma_{i+1},\\
 \Upsilon({}^d\Lambda'_{-\xi})^{-\gamma'}_i
 &\ge \Upsilon({}^d\Lambda_\xi)^{-\gamma}_i
 \ge \Upsilon({}^d\Lambda'_{-\xi})^{-\gamma'}_{i+1},
\end{aligned}
\end{equation}
and $\gamma'\defect(\Lambda'_{-\xi})=\gamma\defect(\Lambda_\xi)-1$.

Choose representatives with $r,s\ge1$ so that the $\gamma$-row of ${}^d\Lambda_\xi$ and the $\gamma'$-row of ${}^d\Lambda'_{-\xi}$ have length $r$; then the other rows have lengths $s$ and $s+1$, respectively. Write these rows as
\[
 \begin{aligned}
  {}^d\Lambda_\xi&:\quad x_1>\cdots>x_r\text{ ($\gamma$-row)},\quad y_1>\cdots>y_s\text{ ($-\gamma$-row)},\\
  {}^d\Lambda'_{-\xi}&:\quad x'_1>\cdots>x'_r\text{ ($\gamma'$-row)},\quad y'_1>\cdots>y'_{s+1}\text{ ($-\gamma'$-row)}.
 \end{aligned}
\]
By \eqref{eq:upsilon}, \eqref{eq:dual-interlace-common} is equivalent to
\begin{equation}\label{eq:strict-row-interlace-one}
\begin{aligned}
 x_i&\ge x'_i>x_{i+1}\quad(1\le i<r), & x_r&\ge x'_r\ge0,\\
 y'_i&>y_i\ge y'_{i+1}\quad(1\le i\le s).
\end{aligned}
\end{equation}

Put $m:=r+s\ge2$ and order the entries of ${}^d\Lambda_\xi$ and ${}^d\Lambda'_{-\xi}$ as $c_1\ge\cdots\ge c_m$ and $c'_1\ge\cdots\ge c'_{m+1}$. For $v\ge0$, let $\mathcal N(v):=\#\{i:c_i\ge v\}$ and $\mathcal N'(v):=\#\{i:c'_i\ge v\}$. Since \eqref{eq:strict-row-interlace-one} maps every entry of ${}^d\Lambda_\xi$ other than $x_1$ to a strictly larger entry of ${}^d\Lambda'_{-\xi}$ ($x_{i+1}\mapsto x'_i$, $y_i\mapsto y'_i$),
\begin{equation*}
 \mathcal N'(v)\ge\mathcal N(v-1)-1\qquad(v\ge1).
\end{equation*}
Using $\sum_{i=k}^{M}d_i=\sum_{v\ge1}\max\{\#\{i:d_i\ge v\}-k+1,0\}$ and $\mathcal N(0)=m$, we obtain for $1\le k\le m$:
\begin{equation}\label{eq:ordered-tail-comparison}
\begin{aligned}
 \sum_{i=k}^{m+1}c'_i
 &=\sum_{v\ge1}\max\{\mathcal N'(v)-k+1,0\}
  \ge\sum_{v\ge1}\max\{\mathcal N(v-1)-k,0\}\\
 &=(m-k)+\sum_{w\ge1}\max\{\mathcal N(w)-k,0\}
  =(m-k)+\sum_{i=k+1}^{m}c_i.
\end{aligned}
\end{equation}

For any symbol $\Lambda$ whose dual ${}^d\Lambda$ has sorted entries $c_1\ge\cdots\ge c_m$, \eqref{eq:symbol-rank} and successive largest-entry deletions give
\begin{equation}\label{eq:symbol-first-rank-entries}
 r_j^{\rm fd}(\Lambda)=
 \begin{cases}
  \displaystyle\sum_{i=j+1}^m c_i-\Bigl\lfloor\frac{(m-j-1)^2}{4}\Bigr\rfloor, & 0\le j<m,\\[
6pt]
  0, & j\ge m,
 \end{cases}
\end{equation}
which is invariant under shifts and independent of the choice of row. Since $\lfloor(m-k+1)^2/4\rfloor-\lfloor(m-k-1)^2/4\rfloor=m-k$, applying \eqref{eq:symbol-first-rank-entries} to \eqref{eq:ordered-tail-comparison} yields $r_{k-1}^{\rm fd}(\Lambda'_{-\xi})\ge r_k^{\rm fd}(\Lambda_\xi)$ for $1\le k\le m$ (and trivially for $k>m$), proving \eqref{eq:symbol-tail-inequality}.
For equality at $k=0$, applying \eqref{eq:symbol-first-rank-entries} with $\lfloor m^2/4\rfloor-\lfloor(m-2)^2/4\rfloor=m-1=(r-1)+s$ gives
\begin{equation*}
\begin{aligned}
 r_0^{\rm fd}(\Lambda'_{-\xi})-r_1^{\rm fd}(\Lambda_\xi)
 &=(c_1-x_1)+\sum_{i=1}^{r-1}(x'_i-x_{i+1}-1)+x'_r\\
 &\quad+\sum_{i=1}^{s}(y'_i-y_i-1)+y'_{s+1}.
\end{aligned}
\end{equation*}
By \eqref{eq:strict-row-interlace-one}, every term on the right is nonnegative, so equality holds if and only if
\[
\begin{gathered}
 x_1=c_1,\qquad
 (x'_1,\ldots,x'_r)=(x_2+1,\ldots,x_r+1,0),\\
 (y'_1,\ldots,y'_{s+1})=(y_1+1,\ldots,y_s+1,0).
\end{gathered}
\]
Undoing one simultaneous shift, this means precisely ${}^d\Lambda'_{-\xi}=\bigl(\partial_\gamma({}^d\Lambda_\xi)\bigr)^{\langle\gamma\gamma'\rangle}$. Since $\gamma\gamma'=-\delta t$ for $\xi=1$ and $\gamma\gamma'=-t$ for $\xi=-1$, this is the formula in Proposition~\ref{cor:first-descent-rule}; \eqref{eq:symbol-first-tail} then follows directly from \eqref{eq:symbol-first-rank-entries}.
\end{proof}

\end{document}